\documentclass[a4paper,10pt]{article}
\usepackage{benstyle}
\usepackage{longtable}

\usepackage[normalem]{ulem}

\DefineNamedColor{named}{Purple}{cmyk}{0.45,0.86,0,0}

\title{
On the numerical stability of Fourier extensions
}
\author{Ben Adcock \\ Department of Mathematics \\ Purdue University \\ USA \\
 \and Daan Huybrechs \\ Department of Computer Science \\ Katholieke Universiteit Leuven \\ Belgium
 \and Jes\'us Mart\'in--Vaquero \\ Department of Applied Mathematics \\ E.T.S.I.I. B\'ejar,
 University of Salamanca \\Spain}
\begin{document}
\maketitle

\begin{abstract}
An effective means to approximate an analytic, nonperiodic function on a bounded interval is by using a Fourier series
on a larger domain.  When constructed appropriately, this so-called Fourier extension is known to
converge geometrically fast in the truncation parameter.  Unfortunately, computing a Fourier extension requires solving an ill-conditioned linear system, and hence one might expect such rapid convergence to be destroyed when carrying out computations in finite precision.  The purpose of this paper is to show that this is not the case.  Specifically, we show that Fourier extensions are actually numerically stable when implemented in finite arithmetic, and achieve a convergence rate that is at least superalgebraic.  Thus, in this instance, ill-conditioning of the linear system does not prohibit a good approximation.

In the second part of this paper we consider the issue of computing Fourier extensions from equispaced data.  A result of Platte, Trefethen \& Kuijlaars states that no method for this problem can be both numerically stable and exponentially convergent.  We explain how Fourier extensions relate to this theoretical barrier, and demonstrate that they are particularly well suited for this problem: namely, they obtain at least superalgebraic convergence in a numerically stable manner.
\end{abstract}

\section{Introduction}\label{s:introduction}

Let $f : [-1,1] \rightarrow \bbR$ be an analytic function.  When
periodic, an extremely effective means to approximate $f$ is via
its truncated Fourier series.  This approximation 
converges geometrically fast in the truncation
parameter $N$, and can be computed efficiently via the Fast
Fourier Transform (FFT).  Moreover, Fourier series possess high
resolution power. One requires an optimal $2$ modes per wavelength
to resolve oscillations, making Fourier methods
well suited for (most notably) PDEs with oscillatory solutions
\cite{naspec}.

For these reasons, Fourier series are extremely widely used in practice.
However, the situation changes completely when $f$ is nonperiodic.
In this case, rather than geometric convergence, one witnesses
the familiar Gibbs phenomenon near $x = \pm 1$ and only linear
pointwise convergence in $(-1,1)$.

\subsection{Fourier extensions}\label{ss:FE}
For analytic and nonperiodic functions, one way to restore the good
properties of a Fourier series expansion (in particular, geometric
convergence and high resolution power) is to approximate
$f$ with a Fourier series on an \textit{extended} domain $[-T,T]$.
Here $T>1$ is a user-determined parameter.  Thus we seek an
approximation $F_N(f)$ to $f$ from the set \bes{ \cG_{N} := \spn
\left \{ \phi_n : |n| \leq N \right \},\qquad \phi_n(x) : =
\frac{1}{\sqrt{2T}} \E^{\I \frac{n \pi}{T} x }. } Although there
are many potential ways to define $F_N(f)$, in
\cite{BoydFourCont,brunoFEP,huybrechs2010fourier} it was proposed
to compute $F_N(f)$ as the best approximation to $f$ on $[-1,1]$
in a least squares sense: \be{ \label{L2} F_N(f) : =
\underset{\phi \in \cG_N}{\operatorname{argmin}}  \| f - \phi \| .
} Here $\nm{\cdot}$ is the standard norm on $\rL^2(-1,1)$---the
space of square-integrable functions on $[-1,1]$. Henceforth, we
shall refer to $F_N(f)$ as the \textit{continuous} Fourier
extension (FE) of $f$.

In \cite{BADHFEResolution,huybrechs2010fourier} it was shown that
the continuous FE $F_N(f)$ converges geometrically fast in
$N$ and has a resolution constant (number of degrees of freedom per
wavelength required to resolve an oscillatory wave) that ranges between $2$ and $\pi$ depending on the
choice of the parameter $T$, with $T \approx 1$ giving close to the
optimal value $2$ (see \S\ref{ss:Tchoice} for a discussion).  Thus the continuous FE successfully
retains the key properties of rapid convergence and high resolution power of a standard Fourier series in the case of nonperiodic functions.

We note that one does not usually
compute the continuous FE \R{L2} in practice.  A more convenient
approach \cite{BADHFEResolution,huybrechs2010fourier} is to replace \R{L2} by the discrete least squares
\be{ \label{discL2} \tilde{F}_N(f) : = \underset{\phi \in
\cG_N}{\operatorname{argmin}} \sum_{|n| \leq N} | f(x_n) -
\phi(x_n) |^2, } for nodes $\{ x_n \}_{|n| \leq N} \subseteq
[-1,1]$.  We refer to $\tilde{F}_N(f)$ as the
\textit{discrete} Fourier extension of $f$.  When chosen suitably---in particular, as in \R{mapCheb}---such nodes ensure that the difference in approximation properties between the extensions
\R{L2} and \R{discL2} is minimal (for details, see \S
\ref{s:FEtypes}).

\subsection{Numerical convergence and stability of Fourier extensions}\label{ss:intro_summary}
The approximation properties of the continuous and discrete FEs were analyzed
in \cite{BADHFEResolution,huybrechs2010fourier}.  Therein it was also observed numerically that the condition numbers of the matrices $A$ and $\tilde{A}$ of the least squares \R{L2} and \R{discL2} are exponentially large in $N$.  We shall confirm this observation later in the paper.  Thus, if $a = (a_{-N},\ldots,a_{N})^{\top}$ is the vector of \textit{coefficients} of the continuous or discrete FE, i.e.\ $F_N(f)$ or $\tilde{F}_N(f)$ is given by $\sum_{|n| \leq N} a_n \phi_n$, one expects small perturbations in $f$ to lead to large errors in $a$.  In other words, the computation of the coefficients of the continuous or discrete FE is ill-conditioned.

Because of this ill-conditioning, it is tempting to think that
FEs will be useless in applications.  At first sight it is reasonable to expect that the good approximation properties of \textit{exact} FEs (i.e.\ those obtained in exact arithmetic) will be destroyed when computing \textit{numerical} FEs in finite precision.  However, previous numerical studies
\cite{BADHFEResolution,BoydFourCont,brunoFEP,huybrechs2010fourier,LyonFESVD,LyonFast} indicate otherwise.  Despite very large condition numbers, one typically obtains an extremely good approximation with a numerical FE, even for poorly behaved functions and in the presence of noise.

The aim of this paper is to give a full explanation of this phenomenon.  This explanation can be summarized as follows.  In computations, one's interest does not lie with the accuracy in computing the coefficient vector $a$, but rather the accuracy of the numerical FE approximation $\sum_{|n| \leq N} a_n \phi_n$.  As we show, although the mapping from a function to its coefficients is ill-conditioned, the mapping from $f$ to its numerical FE is, in fact, well-conditioned.  In other words, whilst the small singular values of $A$ (or $\tilde{A}$) have a substantial effect on $a$, they have a much less significant, and completely quantifiable, effect on the FE itself.

Although this observation explains the apparent stability of numerical FEs, it does not address their approximation properties.  In \cite{BADHFEResolution,huybrechs2010fourier} it was shown that the exact continuous and discrete FEs $F_N(f)$ and $\tilde{F}_N(f)$ converge geometrically fast in $N$.  However, the fact that there may be substantial differences between the coefficients of $F_N(f)$, $\tilde{F}_N(f)$ and those of the numerical FEs, which henceforth we denote by $G_N(f)$ and $\tilde{G}_N(f)$, suggests that geometric convergence may not be witnessed in finite arithmetic for large $N$.  As we show later, for a large class of functions, geometric convergence of $F_N(f)$ (or $\tilde{F}_N(f)$) is typically accompanied by geometric growth of the norm $\nm{a}$ of the exact (infinite-precision) coefficient vector.  Hence, whenever $N$ is sufficiently large, one expects there to be a discrepancy between the exact coefficient
vector and its numerically computed counterpart, meaning that the numerical extensions $G_N(f)$ and $\tilde{G}_N(f)$ may not exhibit the same convergence behaviour.  In the first half of this paper, besides showing stability, we also give a complete analysis and description of the convergence of $G_N(f)$ and $\tilde{G}_N(f)$, and discuss how this differs from that of $F_N(f)$ and $\tilde{F}_N(f)$.

We now summarize the main conclusions of the first half of the paper.  Concerning stability, we have:

\begin{enumerate}
\item The condition numbers of the matrices $A$ and $\tilde{A}$ of the continuous and discrete FEs are exponentially large in $N$ (see \S \ref{ss:condnumb}).
\item The condition number $\kappa(F_N)$ of the exact continuous FE mapping is exponentially large in $N$.  The condition number of the exact discrete FE mapping satisfies $\kappa(\tilde{F}_N) = 1$ for all $N$ (see \S \ref{ss:cond_numb_exact}).
\item The condition number of the numerical continuous and discrete FE mappings $G_N$  and $\tilde{G}_N$ satisfy
\bes{
\kappa(G_N) \lesssim 1/\sqrt{\epsilon},\quad \kappa(\tilde{G}_N) \lesssim 1,\qquad \forall N \in \bbN,
}
where $\epsilon = \epsilon_{\mathrm{mach}}$ is the machine precision used (see \S \ref{ss:stabfn}).
\end{enumerate}
To state our main conclusions regarding convergence, we first require some notation. Let $\cD(\rho)$, $\rho \geq 1$, be a particular one-parameter family of regions in the complex plane related to Bernstein ellipses (see \R{Mapped_Bernstein} and Definition \ref{d:Bernstein}), and define the Fourier extension \textit{constant} \cite{BADHFEResolution,huybrechs2010fourier} by
\be{
\label{cteET}
E(T) = \cot^2 \left ( \frac{\pi}{4 T} \right ).
}
We now have the following:
\begin{enumerate}
\item[4.] Suppose that $f$ is analytic in $\cD(\rho^*)$ and continuous on its boundary. Then the exact continuous and discrete FEs satisfy
\bes{
\| f - F_N(f) \|, \ \| f - \tilde{F}_{N} f \| \leq c_f \rho^{-N},
}
where $\rho = \min \left \{ \rho^*,E(T) \right \}$ and $c_f$ is proportional to $\max_{x \in \cD(\rho)} | f(x) |$ (see \S \ref{s:FEconv}).
\item[5.] For $f$ as in 4.\ the errors of the numerical continuous and discrete FEs satisfy (see \S \ref{ss:numsolnanalysis}):
\enum{
\item[(i)] For $N \leq N_0$ (continuous) or $N \leq N_1 : = 2 N_0$ (discrete), where $N_0$ is a function-independent breakpoint depending on $\epsilon$ and $T$ only, both $\| f - G_N(f) \|$ and $\| f - \tilde{G}_N f \|$ decay like $\rho^{-N}$, where $\rho$ is as in 4.
\item[(ii)] When $N = N_0$ or $N = N_1$, the errors
\bes{
\| f - G_{N_0}(f) \| \approx c_f (\sqrt{\epsilon})^{d_f},\quad \| f - \tilde{G}_{N_1}(f) \| \approx c_f \epsilon^{d_f},
}
where $c_f$ is as in 4.\ and $d_f = \frac{\log \rho}{\log E(T)} \in (0,1]$.
\item[(iii)] When $N>N_0$ or $N > N_1$, the errors decay at least superalgebraically fast down to \textit{maximal achievable accuracies} of order $\sqrt{\epsilon}$ and $\epsilon$ respectively.  In other words,
\bes{
\limsup_{N \rightarrow \infty} \| f - G_{N}(f) \| \lesssim \sqrt{\epsilon},\qquad \limsup_{N \rightarrow \infty} \| f - \tilde{G}_{N}(f) \| \lesssim \epsilon.
}
}
\end{enumerate}

\rem{
In this paper we refer to several different types of convergence of an approximation $f_N \approx f$.  We say that $f_N$ converges \textit{algebraically} fast to $f$ at rate $k$ if $\| f - f_N \| = \ord{N^{-k}}$ as $N \rightarrow \infty$.  If $\nm{f-f_N}$ decays faster than any algebraic power of $N^{-1}$ then $f_N$ is said to converge \textit{superalgebraically} fast.  We say that $f_N$ converges \textit{geometrically} fast to $f$ if there exists a $\rho > 1$ such that $\nm{f-f_N} = \ord{\rho^{-N}}$.  We shall also occasionally use the term \textit{root-exponential} to describe convergence of the form $\nm{f-f_N} = \ordu{\rho^{-\sqrt{N}}}$. }

As we explain in \S \ref{s:stability}, the reason for the disparity between the exact and numerical FEs can be traced to the fact that the system of functions $\{ \E^{\I \frac{n \pi}{T} \cdot}  \}_{n \in
\bbZ}$ forms a \textit{frame} for $\rL^2(-1,1)$.  The inherent redundancy of this frame, i.e.\ the fact that any function $f$ has infinitely many expansions in this system, leads to both the ill-conditioning in the coefficients, as well as the differing convergence between the exact and numerical approximations $F_N$, $\tilde{F}_N$ and $G_N$, $\tilde{G}_N$ respectively.

This aside, observe that conclusion 5.\ asserts that the numerical continuous FE $G_N(f)$ converges geometrically fast in the regime $N < N_0$ down to an error of order $(\sqrt{\epsilon})^{d_f}$, and then at least superalgebraically fast for $N > N_0$ down to a best achievable accuracy of order $\sqrt{\epsilon}$.  Note that $d_f = 1$ whenever $f$ is analytic in $\cD(\rho)$ with $\rho \geq E(T)$.  Thus $G_N$ approximates all sufficiently analytic functions possessing moderately small constants $c_f$ with geometric convergence down to order $\sqrt{\epsilon}$, and this is achieved at $N = N_0$.  For functions only analytic in regions $\cD(\rho)$ with $\rho < E(T)$, or possessing large constants $c_f$, this accuracy is obtained after a further regime of at least superalgebraic convergence.  Note that $c_f$ is large typically when $f$ is oscillatory or possessing boundary layers.  Hence for such functions, even though they may well be entire, one usually still sees the second 
phase of
superalgebraic convergence.

The limitation of $\sqrt{\epsilon}$ accuracy for the numerical continuous FE is undesirable.  Since $\epsilon = \epsilon_{\mathrm{mach}} \approx 10^{-16}$ in practice, this means that one cannot expect to obtain more than $7$ or $8$ digits of accuracy in general.  The condition number is also large---specifically, $\kappa(G_N) \approx 10^{8}$ (see 3.)---and hence the continuous FE has limited practical value.  This is in addition to $G_N(f)$ being difficult to compute in practice, since it requires calculation of $2N+1$ Fourier integrals of $f$ (see \S \ref{sss:continuous}).

On the other hand, conclusion 3.\ shows that the discrete FE is completely stable when implemented numerically.  Moreover, it possesses the same qualitative convergence behaviour as the continuous FE, but with two key differences.  First, the region of guaranteed geometric convergence is precisely twice as large,  $N_1 = 2 N_0$.  Second, the maximal achievable accuracy is on the order of machine precision, as opposed to its square root (see 5.).  Thus, an important conclusion of the first half of this paper is the following: it is possible to compute a numerically stable FE of any analytic function which converges at least superalgebraically fast in $N$ (in particular, geometrically fast for all small $N$), and which attains close to machine accuracy for $N$ sufficiently large.

\rem{
\label{r:numsolver}
This paper is about the discrepancy between theoretical properties of solutions to \R{L2} and \R{discL2} and their numerical solutions when computed with standard solvers.  Throughout we shall consistently use \textit{Mathematica}'s \texttt{LeastSquares} routine in our computations,  though we would like to stress that \textit{Matlab}'s command $\backslash$ gives similar results.  Occasionally, to compare theoretical and numerical properties, we shall carry out computations in additional precision to eliminate the effect of round-off error.  When done, this will be stated explicitly.  Otherwise, it is to be assumed that all computations are carried out as described in standard precision.
}

\subsection{Fourier extensions from equispaced data}\label{ss:intro_summary_equispaced}
In many applications, one is faced with the problem of recovering an analytic function $f$ to high accuracy from its values on an equispaced grid $\left \{ f \left ( \tfrac{n}{M} \right ) : n=-M,\ldots,M \right \}$.  This problem turns out to be quite challenging.  For example, the famous Runge phenomenon states that the polynomial interpolant of this data will diverge geometrically fast as $M \rightarrow \infty$ unless $f$ is analytic in a sufficiently large region.

Numerous approaches have been proposed to address this problem, and thereby `overcome' the Runge phenomenon (see \cite{BoydRunge,TrefPlatteIllCond} for a comprehensive list).  Whilst many are quite effective in practice, ill-conditioning is often an issue.  This was  recently explained by Platte, Trefethen \& Kuijlaars in \cite{TrefPlatteIllCond} (see also \S \ref{ss:PTKrelation}), wherein it was shown that any exponentially convergent method for recovering analytic functions $f$ from equispaced data must also be exponentially ill-conditioned.  As was also proved, the best possible that can be achieved by a stable method is root-exponential convergence.  This profound result, most likely the first of its kind for this type of problem, places an important theoretical benchmark against which all such methods must be measured.

As we show in the first half of this paper, the numerical discrete FE is well-conditioned and has good convergence properties.  Yet it relies on particular interpolation points \R{mapCheb} which are not equispaced.  In the second half of this paper we consider Fourier extensions based on equispaced data.  In particular, if $x_n = \frac{n}{M}$ we study the so-called \textit{equispaced} Fourier extension
\be{
\label{equiExt}
F_{N,M}(f) : = \underset{\phi \in
\cG_N}{\operatorname{argmin}} \sum_{|n| \leq M} | f(x_n) -
\phi(x_n) |^2,
}
and its finite-precision counterpart $G_{N,M}(f)$.  

Our primary interest shall lie with the case where $M = \gamma N$ for some $\gamma \geq 1$, i.e.\ where the number of points $M$ scales linearly with $N$.  In this case we refer to $\gamma$ as the \textit{oversampling} parameter.  Observe that \R{equiExt} results in an $(2M+1) \times (2N+1)$ least squares problem for the coefficients of $F_{N,M}(f)$.  We shall denote the corresponding matrix by $\bar{A}$.

Our main conclusions concerning the exact equispaced FE $F_{N,M}(f)$ are as follows (see \S \ref{ss:equiFEthy}):
\begin{enumerate}
\item[6.] The condition number of $\bar{A}$ is exponentially large as $N,M \rightarrow \infty$ with $M \geq N$.
\item[7.] The condition number of exact equispaced FE mapping $\kappa(F_{N,\gamma N})$ is exponentially large in $N$ whenever $M = \gamma N$ for $\gamma \geq 1$ fixed.  Moreover, the approximation $F_{N,\gamma N}(f)$ suffers from a Runge phenomenon for any fixed $\gamma \geq 1$.  In particular, the error $\| f - F_{N,\gamma N}(f)\|$ may diverge geometrically fast in $N$ for certain analytic functions $f$.
\item[8.] The scaling $M = \ord{N^2}$ is required to overcome the ill-conditioning and the Runge phenomenon in $F_{N,M}$.  In this case, $F_{N,M}(f)$ converges at the same rate as the exact continuous FE $F_N(f)$, i.e.\ geometrically fast in $N$.  Although the condition number of $\bar{A}$ remains exponentially large, the condition number of the mapping $\kappa(F_{N,M})$ is $\ord{1}$ for this scaling.
\end{enumerate}
These results lead to the following conclusion.  The exact (infinite-precision) equispaced FE $F_{N,M}$ with $M = \ord{N^2}$ attains the stability barrier of Platte, Trefethen \& Kuijlaars: namely, it is well-conditioned and converges root-exponentially fast in the parameter $M$.

However,  since the matrix $\bar{A}$ is always ill-conditioned,
one expects there to be differences between the exact
equispaced extension $F_{N,M}(f)$ and its numerical counterpart
$G_{N,M}(f)$.  In practice, one sees both differing stability and
convergence behaviour of $G_{N,M}(f)$, much like in the case of
continuous and discrete FEs.  Specifically, in \S
\ref{ss:equianalysis} we show the following:

\begin{enumerate}
\item[9.] The condition number $\kappa(G_{N,\gamma N})$ satisfies
\bes{ \kappa(G_{N,\gamma N}) \lesssim
\epsilon^{-a(\gamma;T)},\quad \forall N \in \bbN, } where
$\epsilon = \epsilon_{\mathrm{mach}}$ is the machine precision
used, and $0<a(\gamma;T) \leq 1$
is independent of $N$ and satisfies $a(\gamma;T) \rightarrow 0$
as $\gamma \rightarrow \infty$ for fixed $T$ (see \R{agammaT} for the definition of $a(\gamma;T)$).
\item[10.] The error $\| f - G_{N,\gamma N} ( f )\|$ behaves as follows:
\enum{
\item[(i)] If $N<N_2$, where $N_2$ is a function-independent breakpoint, $\| f - G_{N,\gamma N} ( f ) \|$
converges or diverges exponentially fast at the same rate as
$\nm{f-F_{N,\gamma N} ( f ) }$.
\item[(ii)] If $N_2 \leq N < N_1$, where $N_1$ is as introduced previously in \S\ref{ss:intro_summary}, then
$\nm{f-G_{N,\gamma N} ( f )}$ converges geometrically fast at the
same rate as $\nm{f-F_N ( f )}$, where $F_N ( f )$
is the exact continuous FE.
\item[(iii)] When $N = N_1$ the error
\bes{
\| f - G_{N_1,\gamma N_1}(f) \| \approx c_f \epsilon^{d_f-a(\gamma;T)},
}
where $c_f$ and $d_f$ are as in 5.\ of \S \ref{ss:intro_summary}.
\item[(iii)] If $N > N_1$ then $\nm{f-G_{N,\gamma N} ( f )}$
decays at least superalgebraically fast in $N$ down to a maximal achievable accuracy of order $\epsilon^{1-a(\gamma;T)}$.
}
\end{enumerate}
These results show that the condition number of the numerical equispaced FE is bounded whenever $M = \gamma N$, unlike for its exact analogue.  Moreover, after a (function-independent) regime of
possible divergence, we witness geometric convergence of $G_{N,\gamma N}(f)$ down to a certain accuracy.  As in the case of the continuous or discrete FEs, if the function $f$ is sufficiently analytic with small constant $c_f$ then the convergence effectively stops at this point.  If not, we witness a further regime of guaranteed superalgebraic convergence.  But in both cases, the maximal achievable accuracy is of order $\epsilon^{1-a(\gamma;T)}$, which, since $a(\gamma;T) \rightarrow 0$ as $\gamma \rightarrow \infty$, can be made arbitrarily close to $\epsilon$ by increasing $\gamma$.  Note that doing this both improves the condition number of the numerical equispaced FE and yields a
less severe rate of exponential divergence in the region $N<N_2$.  As we show via numerical computation of the relevant constants, double oversampling $\gamma = 2$ with $T=2$ gives perfectly adequate results in most cases.

The main conclusion of this analysis is that numerical equispaced FEs, unlike their exact counterparts, are able to circumvent the stability barrier of Platte, Trefethen \& Kuijlaars to an extent (see \S \ref{ss:PTKrelation} for a more detailed discussion).  Specifically, the numerical FE $F_{N,\gamma N}$ has a bounded condition number, and for all sufficiently analytic functions---namely, those analytic in the region $\cD(E(T))$---the convergence is geometric down to a finite accuracy of order $c_f \epsilon^{1-a(\gamma;T)}$.  This latter observation, namely the fact that the maximal accuracy is nonzero, is precisely the reason why the stability theorem, which requires geometric convergence for all $N$, does not apply.  On the other hand, for all other analytic functions (or those possessing large constants $c_f$) the convergence is at least superalgebraic for $N > N_1$ down to roughly $\epsilon^{1-a(\gamma;T)}$; again not in contradiction with the theorem.  Importantly, one never sees divergence of the numerical FE after 
the finite breakpoint $N_2$.

For this reason, we conclude that equispaced FEs are an attractive method for approximations from equispaced data.  To further support this conclusion we also remark that although the primary concern of this paper is analytic functions, equispaced FEs are also applicable to functions of finite regularity.  In this case, one witnesses algebraic convergence, with the precise order depending solely on the degree of smoothness (see Theorem \ref{t:specconv}).

\subsection{Relation to previous work}
One-dimensional FEs for overcoming the Gibbs and Runge phenomena
were studied in \cite{BoydFourCont} and
\cite{BoydRunge}, and applications to surface parametrizations
considered in \cite{brunoFEP}.  Analysis of the convergence of the
exact continuous and discrete FEs was presented by Huybrechs in \cite{huybrechs2010fourier} and Adcock \& Huybrechs in \cite{BADHFEResolution}.  The issue of resolution power was also addressed in the latter.  The content of the first half of this paper, namely analysis of exact/numerical FEs, follows on directly from this work.

A different approach to FEs, known as the FC--Gram method, was introduced in \cite{lyon2010high}.  This approach forms a central part of an extremely effective method for solving PDEs in complex geometries \cite{albin2011,bruno2010high}.  For previous work on using FEs for PDE problems (so-called Fourier \textit{embeddings}) see \cite{boyd2005fourier,pasquettiFourEmbed}.

Equispaced FEs of the form studied in this paper were
first  independently considered by Boyd \cite{BoydFourCont} and Bruno \cite{bruno2003review}, and later by Bruno et al.\ \cite{brunoFEP}.  In particular,
Boyd \cite{BoydFourCont} describes the use of truncated singular
value decompositions (SVDs) to compute equispaced FEs, and gives
extensive numerical experiments (see also \cite{BoydRunge}).  Bruno focuses on the use of Fourier extensions (also called \emph{Fourier continuations} in the above references) for the description of complicated smooth surfaces. He suggested in \cite{bruno2003review} a weighted least squares to obtain a smooth extension for this purpose, with numerical evidence supporting convergence results in \cite{brunoFEP}.  Most
recently Lyon has presented an analysis of equispaced FEs
computed using truncated SVDs \cite{LyonFESVD}.  In particular,
numerical stability and convergence (down to close to machine
precision) were shown.  In \S \ref{ss:equianalysis} we discuss
this work in more detail (see, in particular, Remark
\ref{r:Lyon}), and give further insight into some of the questions
raised in \cite{LyonFESVD}.

\subsection{Outline of the paper}
The outline of the remainder of this paper is as follows.  In \S \ref{s:FourExt} we recap properties of the continuous and discrete FEs from \cite{BADHFEResolution,huybrechs2010fourier}, including convergence and how to choose the extension parameter $T$.  Ill-conditioning of the coefficient map is proved in \S \ref{s:instability}, and in \S \ref{s:stability} we consider the stability of the numerical extensions and their convergence.  Finally, in \S \ref{s:equispacedI} we consider the case of equispaced FEs.

A comprehensive list of symbols is given at the end of the paper.

\section{Fourier extensions}\label{s:FourExt}
In this section we introduce FEs, and recap salient important aspects of \cite{BADHFEResolution,huybrechs2010fourier}.

\subsection{Two interpretations of Fourier extensions}\label{ss:FEinterp}
There are two important interpretations of FEs which inform their approximation properties and their stability, respectively.  These are described in the next two sections.

\subsubsection{Fourier extensions as polynomial approximations}\label{sss:FEpoly}
The space $\cG_{N}$ can be decomposed as $\cG_{N} = \cC_{N} \oplus \cS_{N}$, where
\bes{
\cC_{N} = \spn \left \{ \cos \tfrac{n \pi}{T} x : n=0,\ldots,N \right \},\quad \cS_{N} = \spn \left \{ \sin \tfrac{n \pi}{T} x : n=1,\ldots,N \right \},
}
consist of even and odd functions respectively.  Likewise, for $f$ we have
\bes{
f(x) = f_{e}(x) + f_{o}(x),\qquad f_e(x) = \tfrac{1}{2} \left [ f(x) + f(-x) \right ],\quad f_{o}(x)= \tfrac{1}{2} \left [ f(x) - f(-x) \right ],
}
and for any FE $f_N$ of $f$:
\be{
\label{FEdecomp}
f_{N} = f_{e,N} + f_{o,N},\quad f_{e,N} \in \cC_N,\  f_{o,N} \in \cS_{N}.
}
Throughout this paper we shall use the notation $f_N$ to denote an arbitrary FE of $f$ when not wishing to specify its particular construction.  From \R{FEdecomp}, it follows that the problem of approximating $f$ via a FE $f_N$ decouples into two problems $f_{e,N}\approx f_e $ and $f_{o,N} \approx f_o$ in the subspaces $\cC_N$ and $\cS_N$ respectively on the half-interval $[0,1]$.

Let us define the mapping $y = y(x) : [0,1] \rightarrow [c(T),1]$ by $ y = \cos
\tfrac{\pi}{T} x$,
where $c(T) = \cos \frac{\pi}{T}$. The functions $\cos \frac{n \pi}{T} x$
and $\sin \frac{(n+1) \pi}{T} x / \sin \frac{\pi}{T} x$ are
algebraic polynomials  of degree $n$ in $y$.  Therefore $\cC_{N}$
and $\cS_{N}$ are (up to multiplication by $\sin \frac{\pi}{T} x$
for the latter) the subspaces $\bbP_N$ and $\bbP_{N-1}$
of polynomials of degree $N$ and $N-1$ respectively in the
transformed variable $y$. Letting \bes{ g_1(y) = f_e(x),\quad
g_2(y) = \frac{f_o(x)}{\sin \frac{\pi}{T} x},\qquad
 g_{1,N}(y) = f_{e,N}(x),\quad g_{2,N}(y) = \frac{f_{o,N}(x)}{\sin \tfrac{\pi}{T} x} ,
}
with $g_{1,N}(y) \in \bbP_{N}$ and $g_{2,N}(y) \in \bbP_{N-1}$, we conclude that the FE approximation $f_N$ in the variable $x$ is completely equivalent to two polynomial approximations in the transformed variable $y \in [c(T),1]$.  

This fact is central to the analysis of FEs.  It allows one to use the rich literature on polynomial approximations to determine the theoretical behaviour of the continuous and discrete FEs (see \S \ref{s:FEconv}).

\rem{
The interpretation of $f_N$ in terms of polynomials is solely for the purposes of analysis.  We always perform computations in the $x$-domain using the standard trigonometric basis for $\cG_N$ (see \S \ref{s:FEtypes}).
}

The interval $[c(T),1] \subseteq (-1,1]$ is not standard.
It is thus convenient to map it affinely to $[-1,1]$.  Let
\bes{  z:  = z(y)= 2 \frac{ y-c(T) }  { 1-c(T) } - 1
\in [-1,1].
}
Observe that $y = y(z) =  c(T) +
\frac{1-c(T)}{2}(z+1)$. Let $m : [0,1]
\rightarrow [-1,1]$ be the mapping $x \mapsto z$, i.e.\ \be{
\label{mdef} z = m(x) = 2 \frac{\cos \frac{\pi}{T} x -
c(T)}{1-c(T)} - 1. } Note that $x = m^{-1}(z) =
\frac{T}{\pi} \arccos \left [ c(T) + \frac{1-c(T)}{2}(z+1) \right
]$.  If we now define \be{ \label{hidef} h_{i}(z) = g_{i}(y(z)),
\quad i=1,2, } then the FE $f_N$ is equivalent to
the two polynomial approximations \be{ \label{E:hNdef} h_{1,N}(z)
= g_{1,N}(y(z)) = f_{e,N}(m^{-1}(z)),\quad h_{2,N}(z) =
g_{2,N}(y(z)) = \frac{f_{o,N}(m^{-1}(z))}{\sin \left (
\frac{\pi}{T} m^{-1}(z) \right)}, } of degree $N$ and $N-1$
respectively in the new variable $z \in [-1,1]$.

\subsubsection{Fourier extensions as frame approximations}\label{sss:FEframe}

\defn{
Let $\rH$ be a Hilbert space with inner product $\ip{\cdot}{\cdot}$ and norm $\nm{\cdot}$.  A set $\{ \phi_n \}^{\infty}_{n =1} \subseteq \rH$ is a frame for $\rH$ if (i) $\spn \{ \phi_n \}^{\infty}_{n=1}$ is dense in $\rH$ and (ii) there exist $c_1,c_2 > 0$ such that
\be{
\label{frameprop}
c_1 \nm{f}^2 \leq \sum^{\infty}_{n=1} | \ip{f}{\phi_n} |^2 \leq c_2 \nm{f}^2,\quad \forall f \in \rH.
}
If $c_1 = c_2$ then $\{ \phi_n \}^{\infty}_{n=1}$ is referred to
as a tight frame. } Introduced by Duffin \& Schaeffer
\cite{duffinshaeffer}, frames are
vitally important in signal processing
\cite{christensen2003introduction}.  Note that all orthonormal,
indeed Riesz, bases are frames, but a frame need not be a basis.
In fact, frames are typically \textit{redundant}: any element $f
\in \rH$ may well have infinitely many representations of the form $f = \sum^{\infty}_{n=1} \alpha_n \phi_n$ with coefficients $\{ \alpha_n
\}^{\infty}_{n=1} \in l^2(\bbN)$.

The relevance of frames to Fourier extensions is due to the following observation:

\lem{
[\cite{BADHFEResolution}] The set $\{ \frac{1}{\sqrt{2T}} \E^{\I
\frac{n \pi}{T} x} \}_{n \in \bbZ}$ is a tight frame for
$\rL^2(-1,1)$ with $c_1=c_2=1$. }
Note that $\{ \frac{1}{\sqrt{2T}} \E^{\I \frac{n \pi}{T}
x} \}_{n \in \bbZ}$ is an orthonormal basis for $\rL^2(-T,T)$: it is precisely the standard Fourier basis on $[-T,T]$.
However, it forms only a
frame when considered as a subset of $\rL^2(-1,1)$.  This fact means that ill-conditioning may well be an issue in numerical algorithms for computing FEs, due to the possibility of redundancies.  As it happens, it is trivial to see that the set $\{ \frac{1}{\sqrt{2T}} \E^{\I \frac{n \pi}{T} x} \}_{n \in \bbZ}$ is redundant: 
\lem{
Let $f \in \rL^2(-1,1)$ be arbitrary, and suppose that $\tilde{f} \in \rL^2(-T,T)$ is such that $f = \tilde{f}$ a.e.\ on $[-1,1]$.  If $\phi_n(x) = \frac{1}{\sqrt{2 T}} \E^{\I \frac{n \pi}{T} x}$ and $\alpha_n = \ip{\tilde{f}}{\phi_n}_{[-T,T]}$, then
\be{
\label{f_rep}
f = \sum_{n \in \bbZ} \alpha_n \phi_n\ \mbox{a.e.}
}
In particular, there are infinitely many sequences $\{ \alpha_n \}_{n \in \bbZ} \in l^2(\bbZ)$ for which $f = \sum_{n \in \bbZ} \alpha_n \phi_n$.
}
\prf{
The sum $\sum_{n \in \bbZ} \alpha_n \phi_n$ is the Fourier series of $\tilde{f}$ on $[-T,T]$.  Thus it coincides with $\tilde{f}$ a.e.\ on $[-T,T]$, and hence $f$ when restricted to $[-1,1]$.  Since there are infinitely many  possible $\tilde{f}$, each giving rise to a different sequence $ \{ \alpha_n \}_{n \in \bbZ}$, the result now follows.
}
This lemma is valid for arbitrary $f \in \rL^2(-1,1)$.  When $f$ has higher regularity---say $f \in \rH^k(-1,1)$, where $\rH^k(-1,1)$ is the $k^{\rth}$ standard Sobolev space on $(-1,1)$---it is useful to note that there exist extensions $\tilde f$ with the same regularity on the torus $\bbT = [-T,T)$.  This is the content of the next result.  For convenience, given a domain $I$, we now write $\nm{\cdot}_{\rH^k(I)}$ for the standard norm on $\rH^k(I)$:

\lem{
\label{l:smoothextension}
Let $f \in \rH^k(-1,1)$ for some $k \in \bbN$.  Then there exists an extension $\tilde f \in \rH^k(\bbT)$ of $f$ satisfying $\| \tilde f\|_{\rH^k(\bbT)} \leq c_k(T) \| f \|_{\rH^k(-1,1)}$, where $c_k(T) > 0$ is independent of $f$.  Moreover, $f = \sum_{n \in \bbZ} \alpha_n \phi_n$, where $\alpha_n = \ip{\tilde{f}}{\phi_n}_{[-T,T]}$ satisfies $\alpha_n  = \ord{n^{-k}}$ as $|n| \rightarrow \infty$.
}
\prf{
The first part of the lemma follows directly from the proof of Theorem 2.1 in \cite{BADHFEResolution}.  The second follows from integrating by parts $k$ times and using the fact that $\tilde f$ is periodic.
}

This lemma, which shall be important later when studying numerical FEs, states that there exist representations of $f$ in the frame $\{ \frac{1}{\sqrt{2T}} \E^{\I \frac{n \pi}{T} x} \}_{n \in \bbZ}$ that have nice (i.e.\ rapidly decaying) coefficients and which cannot grow large on the extended region $[-T,T] $.

\subsection{The continuous and discrete Fourier extensions}\label{s:FEtypes}
We now describe the two types of FEs we consider in the first part of this paper.

\subsubsection{The continuous Fourier extension}\label{sss:continuous}
The continuous FE of $f \in \rL^2(-1,1)$, defined by \R{L2}, is the orthogonal projection onto $\cG_N$.  Computation of this extension involves solving a linear system.  Let us write $F_{N}(f) = \sum^{N}_{n=-N} a_n \phi_{n}$ with unknowns $\{ a_n \}^{N}_{n=-N}$.  If $a = (a_{-N},\ldots,a_{N})^{\top}$ and $b=(b_{-N},\ldots,b_{N})^{\top}$, where
\be{
\label{E:fcoef}
b_{n} =\ip{f}{\phi_n} = \int^{1}_{-1} f(x) \overline{\phi_{n}(x)} \D x,\quad n=-N,\ldots,N,
}
and $A \in \bbC^{(2N+1) \times (2N+1)}$ is the matrix with $(n,m)^{\rth}$ entry
\be{
 \label{E:A}
A_{n,m} = \ip{\phi_m}{\phi_n} = \int^{1}_{-1} \phi_{m}(x)
\overline{\phi_n(x)} \D x,\quad n,m = -N,\ldots,N, } then $a$ is
the solution of the linear system $A a = b$.  We refer to the values $\{ a_n \}^{N}_{n=-N}$ as the \textit{coefficients} of the FE $F_N(f)$.  Note that the matrix $A$ is a
Hermitian positive-definite, Toeplitz matrix with $A_{n,m} =
A_{n-m}$, where $A_0 = \frac{1}{T}$ and $A_n = \frac{\sin
\frac{n \pi }{T}}{n \pi}$ otherwise. In fact, $A$ coincides with the so-called \textit{prolate} matrix \cite{SlepianV,varah}. We
shall discuss this connection further in \S \ref{ss:ASVD}.

For later use, we also note the following characterization of $F_N(f)$:

\prop{ [\cite{BADHFEResolution,huybrechs2010fourier}]
\label{p:exactinterp} Let $F_N(f)$ be the continuous FE \R{L2} of
a function $f$, and let $h_{i}(z)$ and $h_{i,N}(z)$ be given by
\R{hidef} and \R{E:hNdef} respectively (i.e.\ the symmetric
and anti-symmetric parts of $f$ and $f_N$ with the coordinate transformed from the trigonometric argument $x$ to the polynomial argument $z$).  Then $h_{1,N}(z)$ and $h_{2,N}(z)$ are the truncated expansions of $h_1(z)$
and $h_2(z)$ respectively in polynomials orthogonal with respect to the weight functions \be{ \label{weightfns} w_{1}(z) = \left [
(1-z)(z-m(T)) \right ]^{-\frac12},\quad w_{2}(z) = \left [
(1-z)(z-m(T)) \right ]^{\frac12},\quad z \in [-1,1], } where $m(T)
= 1 - 2 \mathrm{cosec}^2 \left ( \tfrac{\pi}{2T} \right ) < -1$.
In other words, $h_{i,N}(z)$, $i=1,2$, is the orthogonal
projection of $h_i(z)$ onto $\bbP_{N+1-i}$ with respect to the
weighted inner product $\ip{\cdot}{\cdot}_{w_i}$ with weight
function $w_i$. }

\subsubsection{The discrete Fourier extension}\label{sss:discrete}
The discrete FE $\tilde{F}_N(f)$ is defined by \R{discL2}.  To use this extension it is first necessary to choose nodes $\{ x_n \}^{N}_{n=-N}$.  This question was considered in \cite{BADHFEResolution}, and a solution was obtained by exploiting the characterization of FEs as polynomial approximations in the transformed variable $z$.

A good system of nodes for polynomial interpolation is given by the
Chebyshev nodes \be{ \label{Cheb} z_n = \cos \left (
\frac{(2n+1) \pi}{2N+2} \right ),\quad n=0,\ldots,N. } Mapping
these back to the $x$-variable and symmetrizing about $x=0$ leads
to the so-called \textit{mapped symmetric Chebyshev} nodes \be{
\label{mapCheb} x_n = - x_{-n-1} = \frac{T}{\pi} \arccos \left [
\frac{1}{2} (1-c(T)) \cos \left ( \frac{(2n+1) \pi}{2N+2} \right
)+ \frac{1}{2}(1+c(T)) \right ],\quad n=0,\ldots,N. } This gives a
set of $2N+2$ nodes.  Therefore, rather than \R{discL2}, we define
the discrete FE by \be{ \label{discL2alt} \tilde{F}_N(f) :  =
\underset{\phi \in \cG'_N}{\operatorname{argmin}}
\sum^{N}_{n=-N-1} | f(x_n) - \phi(x_n) |^2, } from now on, where $
\cG'_N = \cC_N \oplus \cS_{N+1}$.  Exploiting the relation between
FEs and polynomial approximations once more, we now obtain the
following:

\prop{
\label{p:discinterp}
Let $f_N = \tilde{F}_N(f)  \in \cG'_N$ be the discrete FE \R{discL2alt} based on the nodes \R{mapCheb}, and let $h_i(z)$ and $h_{i,N}(z) \in \bbP_N$ be given by \R{hidef} and \R{E:hNdef} respectively.  Then $h_{i,N}(z)$, $i=1,2$ is the $N^{\rth}$ degree polynomial interpolant of $h_i(z)$ at the Chebyshev nodes \R{Cheb}.
}

Write $\phi_{n}(x) = \cos \frac{n \pi}{T} x$, $\phi_{-(n+1)}(x) = \sin \frac{n+1}{T} \pi x$, $n \in \bbN$, and let $\tilde{F}_{N}(f)(x) = \sum^{N}_{n=-N-1} a_n \phi_n(x)$.  If $a = (a_{-N-1},\ldots,a_{N})^{-T}$ and $\tilde A \in \bbR^{(2N+2) \times (2N+2)}$ has $(n,m)^{\rth}$ entry
\be{
 \label{E:A_tilde}
\tilde{A}_{n,m} = \sqrt{\frac{\pi}{N+1}}\phi_{m}(x_n), \qquad n,m=-N-1,\ldots,N,
}
then we have $\tilde A a =\tilde b$, where $\tilde b = (\tilde b_{-N-1},\ldots,\tilde b_N)^{\top}$ and $\tilde b_n = \sqrt{\frac{\pi}{N+1}} f(x_n)$.  

The following lemma concerning the matrix $\tilde A$ will prove useful in what follows:

\lem{[\cite{BADHFEResolution}] \label{l:weightedOpt} The matrix $A_W = (\tilde A)^* \tilde
A$ has entries \bes{ \ip{\phi_n}{\phi_m}_W := \int^{1}_{-1}
\phi_n(x) { \phi_m(x) } W(x) \D x,\quad
n,m=-N-1,\ldots,N, } where $W$ is the positive, integrable weight
function given by $W(x) = \frac{\sqrt{2} \pi}{T}
\frac{\cos \frac{\pi}{2T} x}{\sqrt{ \cos \frac{\pi}{T} x-\cos
\frac{\pi}{T}}}$.}
This lemma implies that the left-hand side of the normal
equations of the discrete FE are the equations of a continuous FE based on the weighted least-squares minimization with weight function $W$.

\subsection{Convergence of exact Fourier extensions}\label{s:FEconv}
A detailed analysis of the convergence of the exact continuous FE, which we now recap, was carried out in \cite{BADHFEResolution,huybrechs2010fourier}.  We commence with the following theorem:

\thm{
[\cite{BADHFEResolution}] \label{t:specconv}
Suppose that $f \in \rH^{k}(-1,1)$ for some $k \in \bbN$ and that $T>1$.  If $F_N(f)$ is the continuous FE of $f$ defined by \R{L2}, then
\begin{equation}
\label{specineq}
\| f - F_N(f) \| \leq c_{k}(T) N^{-k} \| f \|_{\rH^{k}(-1,1)},\quad \forall N \in \bbN,
\end{equation}
where $c_{k}(T) > 0 $ is independent of $f$ and $N$.
}

This theorem confirms \textit{algebraic} convergence of $F_N(f)$ whenever the approximated function $f$ has finite degrees of smoothness, and \textit{superalgebraic} convergence, i.e.\ faster than any fixed algebraic power of $N^{-1}$, whenever $f \in \rC^\infty[-1,1]$.

Suppose now that $f$ is analytic.  Although superalgebraic convergence is guaranteed by Theorem \ref{t:specconv}, it transpires that the convergence is actually geometric.  This is a direct consequence of the interpretation of the $F_N(f)$ as the sum of two polynomial expansions in the transformed variable $z$ (Proposition \ref{p:exactinterp}).  To state the corresponding theorem, we first require the following definition:

\defn{
\label{d:Bernstein}
The Bernstein ellipse $\cB(\rho) \subseteq \bbC$ of index
$\rho \geq 1$ is given by
\bes{ \cB(\rho) = \left \{
\tfrac{1}{2} \left ( \rho^{-1} \E^{\I \theta} + \rho \E^{-\I
\theta} \right ) : \theta \in [-\pi,\pi ] \right \}. } }
Given a compact region bounded by the Bernstein ellipse $\cB(\rho)$, we shall write
\be{
\label{Mapped_Bernstein}
\cD(\rho) \subseteq
\bbC
}
for its image in the complex $x$-plane under the mapping $x
= m^{-1}(z)$, where $m$ is as in \R{mdef}.

\thm{ [\cite{BADHFEResolution}, \cite{huybrechs2010fourier}]
\label{t:expconv} Suppose that $f$ is analytic in
$\cD(\rho^*)$ and continuous on its boundary.  Then $\| f -
F_N(f) \|_{\infty} \leq c_f \rho^{-N}$, where $\rho = \min \left \{ \rho^*,
E(T) \right \}$, $c_f > 0$ is proportional
to $\max_{x \in \cD(\rho)}| f(x) |$, and $E(T)$ is as in \R{cteET}.
}
\prf{ A full proof was given in
\cite[Thm 2.3]{BADHFEResolution}. The expansion $g_N$ of an
analytic function $g$ in a system of orthogonal polynomials with
respect to some integrable weight function satisfies $\| g - g_N
\|_{\infty} \leq c_g \rho^{-N}$, where $c_g$ is proportional to
$\max_{z \in \cB(\rho)} | g(z)|$ \cite{rivlin1990chebyshev}.
In view of Proposition \ref{p:exactinterp}, it remains only to
determine the maximal parameter $\rho$ of Bernstein ellipse
$\cB(\rho)$ within which $h_1(z)$ and $h_2(z)$ are analytic.

The mapping $z = m(x)$ introduces a square-root type singularity into the functions $h_{i}(z)$ at the point $z = m(T) < -1$.  Hence the maximal possible value of the parameter $\rho$ satisfies
\be{
\label{Bernequiv}
\tfrac{1}{2} (\rho + \rho^{-1}) = - m(T).
}
Observe that if $\psi(t) = t + \sqrt{t^2-1} $ then
\be{
\label{Emequiv}
\psi(m(T)) = E(T).
}
Thus, since $\rho > 1$, the solution to \R{Bernequiv} is precisely $\rho = E(T)$.  Conversely, any singularity of $f$ introduces a singularity of $h_i(z)$, which also limits this value.  Hence we obtain the stated minimum.
}
Theorem \ref{t:expconv} shows that if $f$ is analytic in a sufficiently large region (for example, if $f$ is entire) then the rate of geometric convergence is precisely $E(T)$.  Recall that the parameter $T$ can be chosen by the user.  In the next section we consider the effect of different choices of $T$.

\rem{
Although Theorems~\ref{t:specconv} and \ref{t:expconv} are stated for $F_N(f)$, they also hold for the discrete FE $\tilde{F}_N(f)$, since the latter is equivalent to a sum of Chebyshev interpolants (Proposition \ref{p:discinterp}).
}

\subsection{The choice of $T$} \label{ss:Tchoice}
Note that $E(T) \sim 1 + \pi(T-1)$ as $T \rightarrow 1^+$ and $E(T)  \sim \frac{16}{\pi^2} T^2$ when $T \rightarrow \infty$.  Thus, small $T$ leads to a slower rate of geometric convergence, whereas large $T$ gives a faster rate.    As discussed in \cite{BADHFEResolution}, however, a larger value of $T$ leads to a worse resolution power, meaning that more degrees of freedom are required to resolve oscillatory behaviour.  On the other hand, setting $T$ sufficiently close to $1$ yields a resolution power that is arbitrarily close to optimal.

In \cite{BADHFEResolution} a number of fixed values of $T$ were
used in numerical experiments.  These typically give good results, with small values of $T$
being particularly well suited to oscillatory functions.  Another approach for choosing $T$ was also discussed. This involves letting \be{ \label{Tkte}
T =T(N;\epsilon_{\mathrm{tol}})= \frac{\pi}{4} \left ( \arctan
\left ( (\epsilon_{\mathrm{tol}})^{\frac{1}{2N}} \right ) \right )^{-1}, } where
$\epsilon_{\mathrm{tol}} \ll 1$ is some fixed tolerance (note that this
is very much related to the Kosloff Tal--Ezer map in spectral
methods for PDEs \cite{boyd,KTEmapped}---see
\cite{BADHFEResolution} for a discussion).  This choice of $T$,
which now depends on $N$, is such that $E(T)^{-N} = \epsilon_{\mathrm{tol}}$.  Although
this limits the best achievable accuracy of the FE
with this approach to $\ord{\epsilon_{\mathrm{tol}}}$, setting
$\epsilon_{\mathrm{tol}} = 10^{-14}$ is normally sufficient in
practice.  Numerical experiments
in \cite{BADHFEResolution} indicate that this works well,
especially for oscillatory functions.  In fact, since
\be{
\label{Tkteasymp} T(N;\epsilon_{\mathrm{tol}}) \sim 1 - \frac{\log (\epsilon_{\mathrm{tol}})}{\pi N} +
\ord{N^{-2}},\quad N \rightarrow \infty,}
this approach has formally optimal resolution power.

\rem{
The strategy \R{Tkte} is particularly good for oscillatory problems.  However, if this is not a concern, a practical choice appears to be $T=2$.  In this case, the FE has a particular symmetry that can be exploited to allow for its efficient computation in only $\ord{N (\log N)^2}$ operations \cite{LyonFast}.
}

\section{Condition numbers of exact Fourier extensions}\label{s:instability}
The redundancy of the frame $\{ \frac{1}{\sqrt{2 T}} \E^{\I \frac{n \pi}{T} \cdot} \}_{n \in \bbZ}$ means that the matrices associated with the continuous and discrete FEs are ill-conditioned.  We next derive bounds for the condition number of these matrices. The spectrum of $A$ is considered further in \S \ref{ss:ASVD}, and the condition numbers of the FE mappings $f \mapsto F_{N}(f)$ and $f \mapsto \tilde{F}_N(f)$ are discussed in \S \ref{ss:cond_numb_exact}.

\subsection{The condition numbers of the continuous and discrete FE matrices}
\label{ss:condnumb}

\thm{
\label{t:exactCondNumb}
Let $A$ be the matrix \R{E:A} of the continuous FE.  Then the condition number of $A$ is $\ord{E(T)^{2N}}$ for large $N$.  Specifically, the maximal and minimal eigenvalues satisfy
\ea{
\label{exactCondNumb}
T^{-1} \leq \lambda_{\max}(A) \leq 1,\qquad c_1(T) N^{-3} E(T)^{-2N} \leq \lambda_{\min}(A) \leq c_2(T) N^{2} E(T)^{-2N},
}
where $c_1(T)$ and $c_2(T)$ are positive constants with $c_1(T),c_2(T) = \ord{1}$ as $T \rightarrow 1^+$.
}

\prf{
It is a straightforward exercise to verify that
\be{
\label{exactmateval}
\lambda_{\min}(A) =  \min_{\phi \in \cG_N} \left \{\| \phi \|^2 :  \| \phi \|_{[-T,T]} = 1 \right \},\quad \lambda_{\max}(A) =  \max_{\phi \in \cG_N} \left \{\| \phi \|^2 :  \| \phi \|_{[-T,T]} = 1 \right \}.
}
Using the fact that $\| \phi \| \leq \| \phi \|_{[-T,T]}$, we first notice that $\lambda_{\max}(A) \leq 1 $.  On the other hand, setting $\phi =\frac{1}{\sqrt{2 T}}$, we find that $\lambda_{\max}(A) \geq T^{-1}$, which completes the result for $\lambda_{\max}(A)$.

We now consider $\lambda_{\min}(A)$.  Recall that any $\phi \in \cG_{N}$ can be decomposed into even and
odd parts $\phi_{e}$ and $\phi_{o}$, with each function
corresponding to a polynomial in the transformed variable $z$.
Hence,
 \be{ \label{char} \lambda_{\min}(A) =
\min_{\substack{\phi \in \cG_N \\ \phi \neq 0}}  \left \{\frac{\|
\phi \|^2}{\| \phi \|^2_{[-T,T]} } \right \} = \min_{\substack{p_1 \in \bbP_N, p_2 \in \bbP_{N-1} \\ \|p_1\| +  \|p_2\| \neq 0}}
\left \{ \frac{\| p_1 \|^2_{w_1} + \| p_2
\|^2_{w_2}}{\| p_1 \|^2_{w_1,[m(T),1]} + \| p_2
\|^2_{w_2,[m(T),1]}} \right \}, } where $w_i$, $i=1,2$, is given
by \R{weightfns}.  Since the weight function $w_{i}$ is
integrable, we have \be{ \label{intweight} \| p_{i}
\|_{w_i,[m(T),1]} \leq \sqrt{C_{i}(T)} \| p_i
\|_{\infty,[m(T),1]},\quad i=1,2, } where $C_{i}(T) =
\int^{1}_{m(T)} \D w_{i}$, $i=1,2$.  Moreover, by Remez's
inequality, \bes{ \| p \|_{\infty,[m(T),1]} \leq \| T_{N}
\|_{\infty,[m(T),1]} \| p \|_{\infty},\quad \forall p \in \bbP_N,
} where $T_{N} \in \bbP_N$ is the $N^{\rth}$ Chebyshev polynomial.
Since $T_{N}$ is monotonic outside $[-1,1]$, we have $\| T_{N}
\|_{\infty,[m(T),1]} = | T_{N}(m(T)) |$.  Moreover, due to the formula \bes{
T_{N}(x) = \frac{1}{2} \left [ \left ( x - \sqrt{x^2-1} \right )^n
+  \left ( x + \sqrt{x^2-1} \right )^n \right ], } an application
of \R{Emequiv} gives \be{ \label{remez} \| T_{N}
\|_{\infty,[m(T),1]} = \frac{1}{2} \left [ E(T)^N + E(T)^{-N} \right ] < E(T)^N,\quad \forall N \in \bbN, T>1. }
Next we note that $w_{1}(z) \geq D_1(T)
$ and $w_{2}(z) \geq D_2(T) \sqrt{1-z^2} $, $ \forall z \in [-1,1]$, for positive constants $D_1(T)$ and $D_2(T)$.
Moreover, there exist constants $d_1,d_2 > 0$
independent of $T$ such that \bes{ \| p \|_{\infty} \leq d_1 N \|
p \|,\quad \| p \|_{\infty} \leq d_2 N^{\frac{3}{2}}  \| p
\|_v,\quad p \in \bbP_N, } where $v(z) = \sqrt{1-z^2}$ (this follows from expanding $p$ in orthonormal polynomials $\{ p_n \}_{n \in \bbN}$ on $[-1,1]$ corresponding to the weight function $w(z) = 1$, i.e.\ Legendre polynomials, or $w(z) = v(z)$, i.e.\ Chebyshev polynomials of the second kind, and using the known estimate $\| p_n \|_{\infty} = \ordu{n^{\frac12}}$ for the former and $\| p_n \|_{\infty} = \ordu{n^{\frac32}}$ for the latter \cite[chpt.\ X]{bateman}).  Therefore \be{ \label{claim} \| p \|_{\infty}
\leq \frac{d_i}{\sqrt{D_i(T)}} N^{\frac{1+i}{2}} \| p \|_{w_i},\quad
\forall p \in \bbP_N,\quad i=1,2. } Substituting \R{intweight},
\R{remez} and \R{claim} into \R{char} now gives \bes{
\lambda_{\min}(A)\geq \frac{1}{\max \{ C_1(T)/D_1(T),C_2(T)/D_2(T)
\} } N^{-3} E(T)^{-2N}, } which gives the lower bound in
\R{exactCondNumb}.

For the upper bound, we set $p_2 = 0$ and $p_1 = T_N$ in \R{char} to give
\be{
\label{lminupper}
\lambda_{\min}(A) \leq \frac{ \| T_{N} \|^2_{w_1}}{\| T_{N} \|^2_{w_1,[m(T),1]}} \leq \frac{C_1(T)}{\| T_{N} \|^2_{w_1,[m(T),1]}}.
}
Using \R{remez} we note that $\| T_N \|_{\infty,[m(T),1]}  \geq \frac{1}{2} E(T)^N$.  Recall also that $\| p \|_{\infty} \leq d_1 N \| p \|$,$\forall p \in \bbP_N$.  Scaling this inequality to the interval $[m(T),1]$ now gives
\bes{
\| p \|_{\infty,[m(T),1]} \leq d_1 \sqrt{\frac{2}{1-m(T)}}  N \| p \|_{[m(T),1]} = \sqrt{C_3(T)} N\| p \|_{[m(T),1]}.
}
Note also that  $w_1(z) \geq D_3(T)$, $\forall z \in [m(T),1]$.  Therefore,
\eas{
\| T_N  \|_{w_1,[m(T),1]} \geq \sqrt{D_3(T)}  \| T_N  \|_{[m(T),1]} \geq  \frac{\sqrt{D_3(T)}}{\sqrt{C_3(T)} N} \| T_N \|_{\infty,[m(T),1]} \geq \frac{\sqrt{D_3(T)} }{2 \sqrt{C_3(T)} N} E(T)^N.
}
Substituting this into \R{lminupper} now gives the result.
}

We now consider the case of the discrete FE:

\thm{ \label{t:discCondNumb} Let $\tilde A$ be the matrix
\R{E:A_tilde} of the discrete FE.  Then the
condition number of $\tilde A$ is $\ord{E(T)^{N}}$
for large $N$. Specifically, the maximal and minimal singular
values of $\tilde{A}$ satisfy
\ea{ \label{discCondNumb} c_1(T) \leq
\sigma_{\max}(\tilde A) \leq c_2(T) N^{\frac{3}{2}},\qquad d_1(T)
N^{-\frac32} E(T)^{-N} \leq \sigma_{\min}(\tilde A) \leq d_2(T)
N^{\frac52} E(T)^{-N}, }
where $c_1(T),c_2(T),d_1(T),d_2(T)$ are positive constants that
are $\ord{1}$ as $T \rightarrow 1^+$. } \prf{ Using Lemma
\ref{l:weightedOpt}, the values $\sigma^2_{\min}(\tilde A)$ and
$\sigma^2_{\max}(\tilde A)$ may be expressed as in
\R{exactmateval} (with $\nm{\cdot}$ replaced by $\nm{\cdot}_W$).
Note that $W(0) \| \phi \|^2 \leq  \| \phi \|^2_W \leq \|
\phi \|^2_{\infty} \int^{1}_{-1} \D W.$ It is a
straightforward exercise (using the bound \R{claim} and the fact
that $\phi$ can be expressed as the sum of two polynomials) to
show that $\| \phi \|_{\infty} \leq C_1(T) N^{\frac32} \| \phi \|$,
where $C_1(T) = \ord{1}$ as $T \rightarrow 1^+$.  Thus we obtain
\bes{ W(0) \frac{\| \phi \|^2}{\| \phi \|^2_{[-T,T]}}  \leq
\frac{\| \phi \|^2_W}{\| \phi \|^2_{[-T,T]}} \leq \left ( C_1(T)^2
\int^{1}_{-1} \D W \right ) N^3 \frac{\| \phi \|^2}{\| \phi
\|^2_{[-T,T]}} . } The result now follows immediately from the
bounds \R{exactCondNumb}. }

Theorems \ref{t:exactCondNumb} and \ref{t:discCondNumb} demonstrate that the condition numbers of the continuous and discrete FE matrices grow exponentially in $N$.  This establishes conclusion 1.\ of \S \ref{s:introduction}.

\rem{
\label{r:factorization}
Although exponentially large, the matrix of the discrete FE is substantially less poorly conditioned than that of the continuous FE.  In particular, the condition number is of order $E(T)^{N}$ as opposed to $E(T)^{2N}$.  This can be understood using Lemma \ref{l:weightedOpt}.  The normal form $A_W = (\tilde A)^* \tilde A$ of the discrete FE matrix is a continuous FE matrix with respect to the weight function $A_W$.  Hence $\kappa(\tilde A) =\sqrt{\kappa(A_W)} \approx \sqrt{\kappa(A)} \approx E(T)^N$.  As we shall see later, this property also translates into superior performance of the numerical discrete FE over its continuous counterpart (see \S \ref{ss:numsolnanalysis}).
}

Since the constants in Theorems \ref{t:exactCondNumb} and \ref{t:discCondNumb} are bounded as $T \rightarrow 1^+$, this allows one also to determine the condition number in the case that  $T \rightarrow 1^+$ as $N \rightarrow \infty$ (see \S \ref{ss:Tchoice}).  In particular, if $T$ is given by \R{Tkte}, then $\kappa(A)$ and $\kappa(\tilde A)$ are (up to possible small algebraic factors in $N$) of order $(\epsilon_{\mathrm{tol}})^{-2}$ and $(\epsilon_{\mathrm{tol}})^{-1}$.

\subsection{The singular value decomposition of $A$}\label{ss:ASVD}
Although we have now determined the condition number of $A$, it is possible to give a rather detailed analysis of its spectrum.  This follows from the identification of $A$ with the well-known prolate matrix, which was analyzed in detail by Slepian \cite{SlepianV,varah}.  We now review some of this work.

\begin{figure}
\begin{center}
$\begin{array}{ccc}
\includegraphics[width=6.25cm]{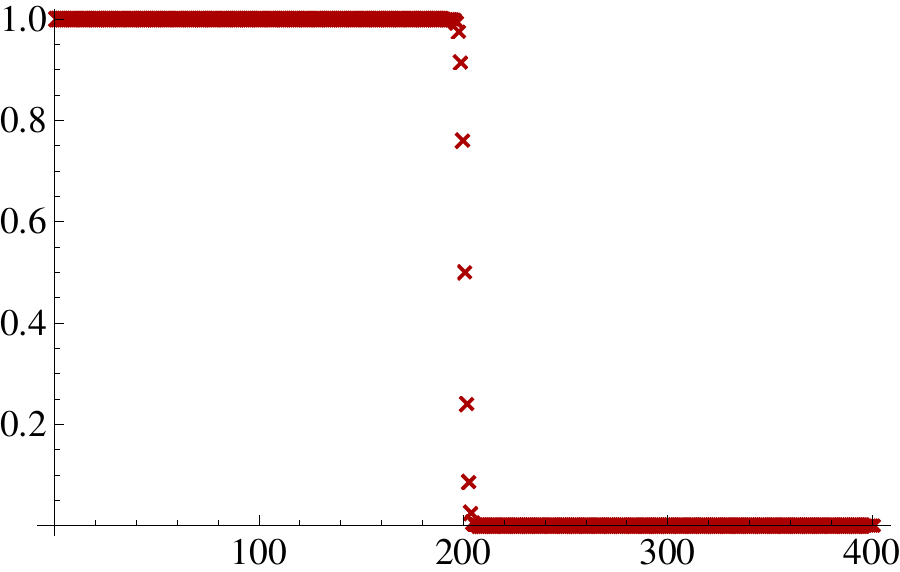}  & \hspace{1pc} &
\includegraphics[width=6.25cm]{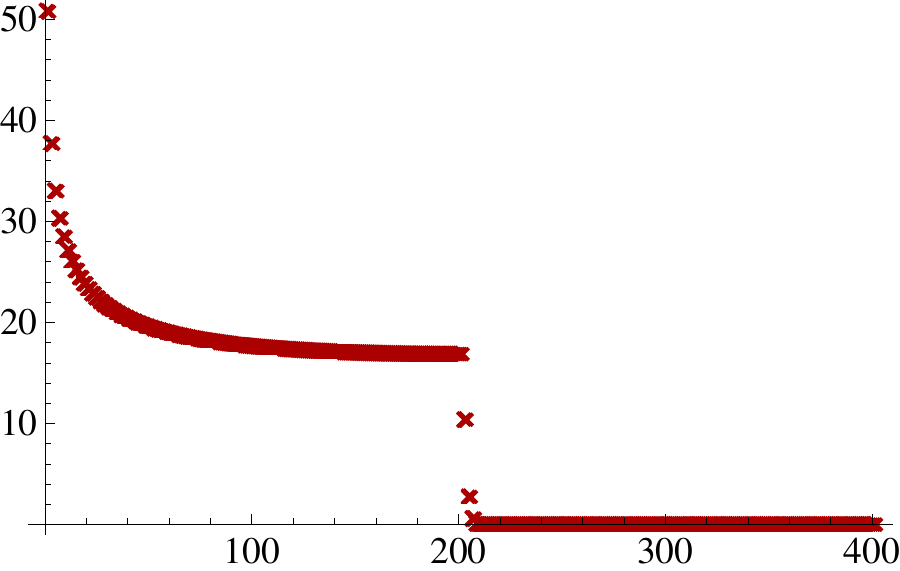}
\end{array}$
\caption{\small  Eigenvalues of the matrices \R{E:A} (left) and
\R{E:A_tilde} (right) for $N=200$ and $T=2$.} \label{f:ExactEvals}
\end{center}
\end{figure}

Following Slepian \cite{SlepianV}, let
$P(N,W) \in \bbC^{N \times N}$ be the prolate matrix with entries \bes{
P(N,W)_{m,n} = \left \{ \begin{array}{cl} \frac{\sin
2 \pi W (m-n)}{\pi (m-n)} &  m \neq n \\ 2 W & m = n, \end{array} \right . \quad m,n=0,\ldots,N-1,
}
where $0<W <\frac{1}{2}$ is fixed, and write $1> \lambda_0(N,W)> \ldots> \lambda_{N-1}(N,W)> 0$ for its eigenvalues.  Note that
 \be{
  \label{lambdasymmetry}
\lambda_{k}(N,\tfrac{1}{2}-W) =  1 - \lambda_{N-1-k}(N,W).
}
The following asymptotic results are found in \cite{SlepianV}:
\begin{itemize}
\item[(i)] For fixed and small $k$, 
\be{ \label{smalleval} 1- \lambda_{k}(N,W) \sim \sqrt{\pi}
(k!)^{-1} 2^{(14k+9)/4} \alpha^{(2k+1)/4} (2-\alpha)^{-(k+1/2)}
N^{k+1/2} \beta ^{-N}, } where $\alpha = 1 - \cos 2 \pi W$ and
$\beta = 
\frac{\sqrt{2}+\sqrt{\alpha}}{\sqrt{2}-\sqrt{\alpha}} $.

\item[(ii)] For large $N$ and $k$ with $k = \lfloor 2 W N (1-\epsilon) \rfloor$ and $0 < \epsilon < 1$, $1 - \lambda_{k}(N,W) \sim \E^{-c_1 - c_2 N}$ for explicitly known constants $c_1,c_2$ depending only on $W$ and $\epsilon$.
\item[(iii)] For large $N$ and $k$ with $ k = \lfloor 2 W N + (b/\pi) \log N \rfloor$, $\lambda_{k}(N,W) \sim \frac{1}{1+\E^{\pi b}}$.
\end{itemize}
(Slepian also derives similar asymptotic results for the
eigenvectors of $P(N,W)$ \cite{SlepianV}).  From these
results we conclude that the eigenvalues of the prolate matrix
cluster exponentially near $0$ and $1$ and have a transition
region of width $\ord{\log N}$ around $k = 2 W N$.  This is shown
in Figure \ref{f:ExactEvals}.

The matrix $A$ of the continuous FE is precisely the prolate matrix
$P(2N+1,\frac{1}{2T})$.  In this case, the parameter $\beta$ in \R{smalleval} is given by
\bes{
\beta = \frac{\sqrt{2}+\sqrt{\alpha}}{\sqrt{2}-\sqrt{\alpha}} = \cot^2
\left (\frac{\pi}{4 T} \right ) =E(T).}
Applying Slepian's analysis, we now see that the eigenvalues of $A$ cluster exponentially at rate $E(T)^2$ near zero and one (note that $A$
corresponds to a prolate matrix of size $2N$), and in particular, that
the condition number is $\ord{E(T)^{2N}}$.  The latter estimate agrees with that given in Theorem \ref{t:exactCondNumb}.  We remark, however, that Theorem
\ref{t:exactCondNumb} gives bounds for the minimal
eigenvalue of $A$ that hold for all $N$ and $T$, unlike
\R{smalleval}, which holds only for fixed $T$ and sufficiently
large $N$.  Hence Theorem \ref{t:exactCondNumb} remains
valid when $T$ is varied with $N$, an option which, as discussed
in \S \ref{ss:Tchoice}, can be advantageous in practice.

Since the matrix $\tilde A$ of the discrete FE is related to $A$ (see Lemma \ref{l:weightedOpt}), we
expect a similar structure for its singular values.  This is
illustrated in Figure \ref{f:ExactEvals}.  As is evident, the only qualitative difference between $\tilde A$ and $A$ is found in the large
singular values. The other features---the narrow transition
region and the exponential clustering of singular values near $0$---are much the same.

\rem{
\label{sss:T2}
The choice $T=2$ ($W = \frac{1}{4}$) is special.  As shown by \R{lambdasymmetry}, the eigenvalues $\lambda_{k}(N,W)$ are symmetric in this case, and the transition region occurs at $k = \frac{1}{2} N$.  This is unsurprising.  When $T=2$, the frame $\{ \E^{\I \frac{n \pi}{2} x} \}_{n \in \bbZ}$ decomposes into two orthogonal bases, related to the sine and cosine transforms.  Using this decomposition and the associated discrete transforms for each basis, M. Lyon has introduced an $\ord{N (\log N)^2 }$ complexity algorithm for computing FEs \cite{LyonFast}.
}

\subsection{Numerical examples}\label{ss:S5numexamp}

\begin{figure}
\begin{center}
$\begin{array}{ccc}
 \includegraphics[width=6.25cm]{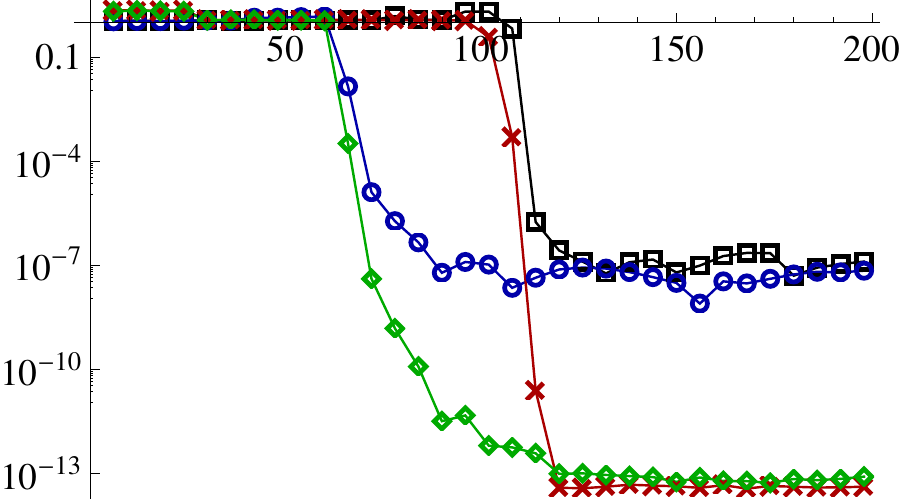}    & \hspace{1pc} &\includegraphics[width=6.25cm]{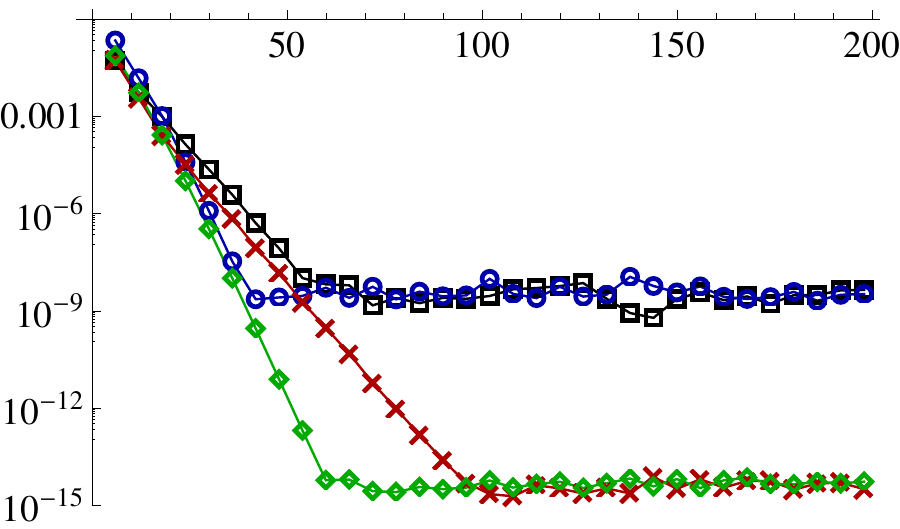}      \\
 f(x) = \E^{25 \sqrt{5} \pi \I x} &&  f(x) = \frac{1}{1+25
 x^2}  \\ \\
  \includegraphics[width=6.25cm]{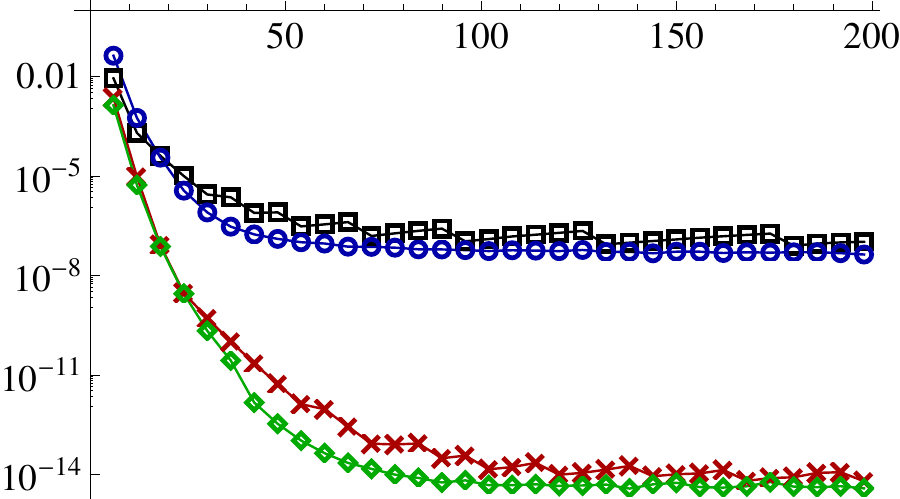}    & &\includegraphics[width=6.25cm]{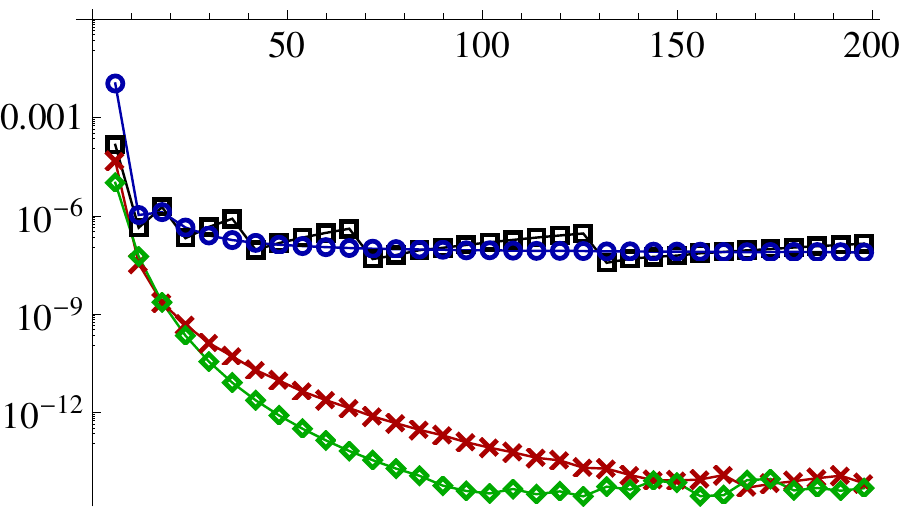}      \\
 f(x) = \frac{1}{8-7 x} & &  f(x) = | x|^7
  \end{array}$
 \caption{\small The error $\| f - f_N \|_{\infty}$,
where $f_N = G_N(f)$ (squares and circles) or $f_N =
\tilde{G}_N(f)$ (crosses and diamonds) and $T=2$ (squares/crosses)
or $T = T(N;\epsilon_{\mathrm{tol}})$ (circles/diamonds) with
$\epsilon_{\mathrm{tol}} = 10^{-14}$.}  \label{f:FnApp}
\end{center}
\end{figure}

We now consider several numerical examples of the continuous and discrete FEs.  In Figure \ref{f:FnApp} we plot the error $\| f - f_N \|_{\infty}$ against $N$ for various choices of $f$.  Here the extension $f_N$ is the numerically computed continuous or discrete FE---i.e.\ the result of solving the corresponding linear system in standard precision (recall Remark \ref{r:numsolver}).  Henceforth, we use the notation $G_N(f)$ and $\tilde{G}_N(f)$ for these \textit{numerical} extensions, so as to distinguish them from their \textit{exact} counterparts $F_N(f)$ and $\tilde{F}_N(f)$.  Note that the word `exact' in this context refers to exact arithmetic.  We do not mean exact in the sense that $F_N(f) = f$ for $f \in \cG_N$.

At first sight, Figure \ref{f:FnApp} appears somewhat surprising: for all three functions we obtain good accuracy, and there is no drift or growth in the error, even in the case where $f$ is nonsmooth or has a complex singularity near $x=0$.  Evidently the ill-conditioning of the FE matrices established in Theorems \ref{t:exactCondNumb} and \ref{t:discCondNumb} appears to have little effect on the numerical extensions $G_N(f)$ and $\tilde{G}_{N}(f)$.  The purpose of \S \ref{s:stability} is to offer an explanation of this phenomenon.

In Figure \ref{f:FnApp} we also compare two choices of $T$: fixed $T=2$ and the $N$-dependent value \R{Tkte} with $\epsilon_{\mathrm{tol}} = 10^{-14}$.  Note that the latter typically outperforms the fixed value $T=2$, especially for oscillatory functions.  This is unsurprising in view of the discussion in \S \ref{ss:Tchoice}.  

Figure \ref{f:FnApp} also illustrates an important disadvantage of the continuous FE: namely, the approximation error levels off at around $\sqrt{\epsilon_{\mathrm{mach}}}$, where $\epsilon_{\mathrm{mach}} \approx 10^{-16}$ is the machine precision used, as opposed to around $\epsilon_{\mathrm{mach}} $ for the discrete extension.  Our analysis in \S \ref{s:stability} will confirm this phenomenon.  Note that the differing behaviour between the continuous and discrete extensions in this respect can be traced back to the observation made in Remark \ref{r:factorization}.

\subsection{Condition numbers of the exact continuous and discrete FE mappings}\label{ss:cond_numb_exact}
The exponential growth in the condition numbers of the continuous and discrete FE matrices imply extreme sensitivity in the FE coefficients to perturbations.  However, the numerical results of Figure \ref{f:FnApp} indicate that the FE approximations themselves are far more robust.  Although we shall defer a full explanation of this difference to \S \ref{s:stability}, it is possible to give a first insight by determining the condition numbers of the mappings $F_N$ and $\tilde{F}_N$.

For vectors $b \in \bbC^{2N+1}$ and $\tilde{b} \in \bbC^{2N+2}$ let us write, with slight abuse of notation, $F_N(b)$ and $\tilde{F}_N(\tilde b)$ for the corresponding continuous and discrete Fourier extensions whose coefficient vectors are the solutions of the linear systems $A a = b$ and $\tilde{A} a = \tilde{b}$ respectively.   We now define the condition numbers
\be{
\label{conditionFE}
\kappa(F_N) =  \sup \left \{ \| F_{N}(b) \| : b \in \bbC^{2N+1},\ \nm{b} = 1 \right \},\ \kappa(\tilde{F}_N) = \sup \left \{ \| F_{N}(b) \|_W : b \in \bbC^{2N+2},\ \nm{b} = 1 \right \}.
}
Here $\nm{\cdot}$ denotes the usual $l^2$ vector norm, and $W$ is the weight function of Lemma \ref{l:weightedOpt}.  Note that \R{conditionFE} gives the \textit{absolute} condition numbers of $F_N$ and $\tilde{F}_N$, as opposed to the more standard \textit{relative} condition number \cite{TrefethenBau}.  The key results of this paper can easily be reformulated for the latter.  However, we shall use \R{conditionFE} throughout, since it coincides with the definition given in \cite{TrefPlatteIllCond} for linear mappings such as FEs.  The work of \cite{TrefPlatteIllCond} will be particularly relevant when considering equispaced FEs in \S \ref{s:equispacedI}.

We now have the following result:
\lem{
\label{l:ExactCondNumb}
The condition numbers of the exact continuous and discrete FEs satisfy
\bes{
\kappa(F_N)= 1/\sqrt{\lambda_{\min}(A)},\qquad \kappa(\tilde{F}_N) = 1.
}
}
\prf{
Write $F_N(b) = \sum^{N}_{n=-N} a_n \phi_n$, where $A a = b$.  We have $\| F_N(b) \|^2 = a^* A a = b^* A^{-1} b$, and therefore $\kappa(F_N) = 1/\sqrt{\lambda_{\min}(A)}$, as required.  For the second result, we note that $\| \tilde{F}_N(\tilde{b}) \|^2 = a^* A_W a$, where $A_W = (\tilde{A})^* \tilde{A}$ is the matrix of Lemma \ref{l:weightedOpt}.  Since $\tilde A a = \tilde b$ the second result now follows.
}
As with the FE matrices, this lemma shows that condition number of the discrete mapping $\tilde F_N$, which is identically one, is much better than that of the continuous mapping $F_N$.  Similarly, the reason can be traced back to Remark \ref{r:factorization}.  Note that this lemma establishes 2.\ of \S \ref{s:introduction}.

At first, it may seem that the fact that $\kappa(\tilde{F}_N) =1 $ explains the observed numerical stability in Figure \ref{f:FnApp}.  However, since $\lambda_{\min}(A)$ is exponentially small (Theorem \ref{t:exactCondNumb}), the above lemma clearly does not explain the lack of drift in the numerical error in the case of the continuous FE.  This is symptomatic of a larger issue: in general, the exact FEs $F_N(f)$ and $\tilde{F}_N(f)$ differ substantially from their numerical counterparts $G_N(f)$ and $\tilde{G}_N(f)$.  As we show in the next section, there are important differences in both their stability \textit{and} their convergence.  In particular, any analysis based solely on $F_N$ and $\tilde{F}_N$ is insufficient to describe the behaviour of the numerical extensions $G_N$ and $\tilde{G}_N$.

\section{The numerical continuous and discrete Fourier extensions}\label{s:stability}
We now analyze the numerical FEs $G_N$ and $\tilde{G}_N$, and describe both when and how they differ from the exact extensions $F_N$ and $\tilde{F}_N$.

\subsection{The norm of the exact FE coefficients}\label{ss:magnitude}
In short, the reason for this difference is as follows.  Since the FE matrices $A$ and $\tilde A$ are so ill-conditioned, the coefficients of the exact FEs $F_N$ and $\tilde{F}_N$ will not usually be obtained in finite precision computations.    To explain exactly how this affects stability and convergence, we first need to determine when this will occur.  We require the following theorem:

\thm{ \label{t:coefficients} Suppose that $f$ is analytic in $\cD(\rho^*)$ and continuous on its boundary.  If $a \in \bbC^{2N+1}$ is the vector of coefficients
of the continuous FE $F_N(f)$ then \be{ \label{coefbd1} \| a \| 
\leq c_f \left \{ \begin{array}{ll}  \left ( \frac{E(T)}{\rho^*}
\right )^N & \rho^* < E(T ), \\  N & \rho^* \geq E(T), \end{array}
\right. } where $c_f$ is proportional to $\max_{x \in \cD(\rho)} | f(x) |$.  If $f
\in \rL^2(-1,1)$, then \be{ \label{coefbd2} \| a
\| \leq c \| f \| E(T)^N, }
for some $c>0$ independent of $f$ and $N$.
 }
\prf{ Write $F_N(f) = f_N = f_{e,N}
+ f_{o,N}$, where $f_{e,N}$ and $f_{o,N}$ are the even and odd
parts of $f_N$ respectively.  Since the set $\{ \phi_n \}_{n \in
\bbZ}$ is orthonormal over $[-T,T]$ we find that \eas{ \| a \| =
\| f_N \|_{[-T,T]}  \leq 2 \left ( \| f_{e,N} \|_{[0,T]} + \|
f_{o,N} \|_{[0,T]} \right ) \leq 2 \sqrt{T} \left ( \| f_{e,N}
\|_{\infty,[0,T]} + \| f_{o,N} \|_{\infty,[0,T]} \right ). }
Recall from \S \ref{sss:FEpoly} that $f_{e,N}(x) = h_{1,N}(z)$ and
$f_{o,N}(x) =\sin  \left ( \tfrac{\pi}{T} m^{-1}(z) \right )
h_{2,N}(z)$, where $h_{i,N} \in \bbP_{N+1-i}$, $i=1,2$, is defined
by \R{E:hNdef}.  Thus, $\| a \| \leq c \left ( \| h_{1,N}
\|_{\infty,[m(T),1]} + \| h_{2,N} \|_{\infty,[m(T),1]} \right )$
for some $c >0$.  Consider
$h_{1,N}(z)$.  This is precisely the expansion of the function $h_1(z) = f_1(m^{-1}(z))$ in
polynomials $\{ p_n \}^{\infty}_{n=0}$ orthogonal with respect to
the weight function $w_{1}$: i.e.\ $h_{1,N} = \sum^{N}_{n=0}
\ip{h_1}{p_n}_{w_1} p_n$.
Therefore \bes{ \| h_{1,N} \|_{\infty,[m(T),1]} \leq
\sum^{N}_{n=0} | \ip{h_1}{p_n}_{w_1} | \| p_n
\|_{\infty,[m(T),1]}. } It is known that $\| p_n
\|_{\infty,[m(T),1]}  \leq c E(T)^n$ \cite{huybrechs2010fourier}.
Also, since $h_1$ is analytic in $\cB(\rho^*)$ we have $| \ip{h_1}{p_n}_{w_1} |
\leq c_{f} (\rho^*)^{-n}$.  Hence \bes{ \| h_{1,N}
\|_{\infty,[m(T),1]} \leq c_f \sum^{N}_{n=0} \left (
E(T) / \rho^* \right )^{n}, } which gives \R{coefbd1}.  For
\R{coefbd2} we use the bound $| \ip{h_1}{p_n}_{w_1} | \leq \| h_1
\|_{w_1} \leq c \| f \|$ instead. }

\cor{ Let $f$ be as in Theorem \ref{t:coefficients}.  Then the
vector of coefficients $a \in \bbC^{2N+2}$ of the discrete Fourier
extension $\tilde{F}_N(f)$ of $f$ satisfies the same bounds as
those given in Theorem \ref{t:coefficients}. } \prf{The functions
$h_{i,N}$, $i=1,2$ are the polynomial interpolants of $h_i$ at the
nodes \R{Cheb} (Proposition \ref{p:discinterp}). Write
$h_{i,N}(z) = \sum^{N}_{n=0} \tilde d_n T_n(z)$, where
$T_n(z)$ is the $n^{\rth}$ Chebyshev polynomial, and let
$\hat{d}_n = \ip{h_i}{T_n}_{w}$ be the Chebyshev polynomial
coefficient of $h_i$.  Note that $| \hat{d}_n | \leq c_f
(\rho^*)^{-n}$.  Due to aliasing formula $\tilde d_n
= \hat{d}_n + \sum_{k \neq 0} (\hat{d}_{2kN+n} + \hat{d}_{2kN-n})$
(see \cite[Eqn. (2.4.20)]{SMSD}) we obtain \bes{ | \tilde d_n |
\leq c_f \left ( (\rho^*)^{-n} + \sum^{\infty}_{k=1} (\rho^*)^{-2 k N - n}
+\sum^{\infty}_{k=1}(\rho^*)^{-2kN + n} \right ) \leq c_f \left (
(\rho^*)^{-n} + (\rho^*)^{n-2N} \right ) \leq c_f (\rho^*)^{-n}.} The
result now follows along the same lines as the proof of Theorem
\ref{t:coefficients}. }

To compute the continuous or discrete FE we need to
solve the linear system $A a = b$ (respectively $\tilde A a
=\tilde b$).  When $N$ is large, the columns of $A$ ($\tilde A$)
become near-linearly dependent, and, as shown in \S
\ref{ss:ASVD}, the numerical rank of $A$ is roughly $1/T$ times
its dimension.  Now suppose we solve $A a = b$ with a standard numerical solver.
Loosely speaking, the solver will use the extra degrees of freedom
to construct approximate solutions $a'$ with small norms. The previous theorem and corollary therefore suggest the following.  In general, only in those cases where $f$ is analytic with $\rho^* \geq E(T)$ can we expect the theoretical coefficient vector $a$ to be produced by the numerical solver for all $N$.  Outside of this case, we may well have that $a' \neq a$
for sufficiently large $N$, due to the potential for exponential growth of the latter.  Hence, in this case, the numerical extension $G_N(f)$ will not
coincide with the exact extension $F_N(f)$.

This raises the following question: if the numerical solver does
not give the exact coefficients vector, then what does it yield?
The following proposition confirms the existence of infinitely
many approximate solutions of the equations $A a = b$ with small
norm coefficient vectors:

\prop{
\label{p:appcoeff}
Suppose that $f \in \rH^k(-1,1)$.  Then there exist $a^{[N]} \in \bbC^{2N+1}$, $N \in \bbN$, satisfying
\be{
\label{coeffbd}
\| a^{[N]} \| \leq c_k(T) \| f \|_{\rH^k(-1,1)},
}
and
\be{
\label{coefferr}
\| A a^{[N]} - b \| \leq c_k(T) N^{-k}  \| f \|_{\rH^k(-1,1)},
}
where $c_k(T)$ is the constant of Lemma \ref{l:smoothextension}.  Moreover, if $g_N = \sum_{|n| \leq N} a^{[N]}_{n} \phi_n$ then
\be{
\label{apperr}
\| f - g_N \| \leq c_k(T) N^{-k} \| f \|_{\rH^k(-1,1)}.
}
}
\prf{
Let $\tilde f \in \rH^k(\bbT)$ be the extension guaranteed by Lemma \ref{l:smoothextension}, and write $a^{[N]}$ for the vector of its first $2N+1$ Fourier coefficients on $\bbT = [-T,T)$.  By Bessel's inequality, $\| a^{[N]} \| \leq \| \tilde f \|_{[-T,T]} \leq c_k(T) \| f \|_{\rH^k(-1,1)}$ which gives \R{coeffbd}.  For \R{coefferr}, we merely note that $(A a^{[N]} - b)_n = \ip{f - g_N}{\phi_n}$.  Using the frame property \R{frameprop} we obtain $\| A a^{[N]} - b \| \leq \| f - g_N \|$.  Thus, \R{coefferr} follows directly from \R{apperr},
 and the latter is a standard result of Fourier analysis (see \cite[eqn. (5.1.10)]{SMSD}, for example).
}

This proposition states that there exist vectors with norms bounded
 independently of $N$ that solve the equations $A a =b$
 up to an error of order $N^{-k}$.  Moreover, these vectors
   yield extensions which converge algebraically fast to $f$ at rate $k$.  Whilst it does not imply that these are the vectors produced by the numerical solver, it does indicate that, in the case where the exact extension $F_N(f)$ or $\tilde F_N(f)$ has a large coefficient norm, geometric convergence of the numerical extension $G_N(f)$ or $\tilde G_N(f)$ may be sacrificed for superalgebraic convergence so as to retain boundedness of the computed coefficients.

\begin{figure}[t]
\begin{center}
$\begin{array}{ccc}
\includegraphics[width=6.25cm]{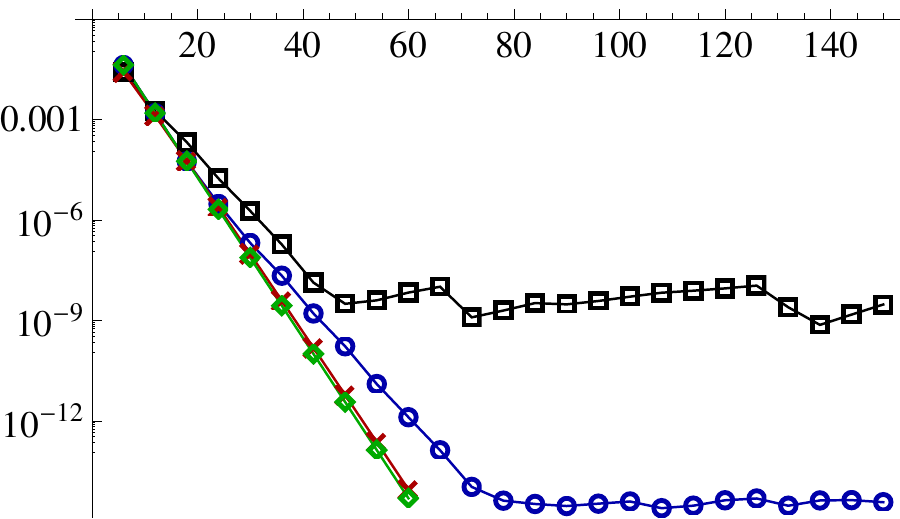}  & \hspace{1pc} &
\includegraphics[width=6.25cm]{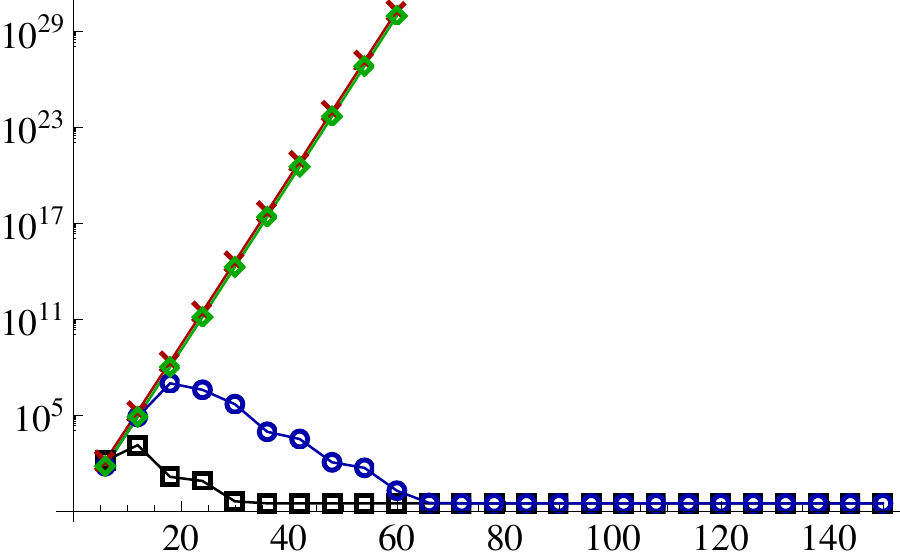}
\\
\includegraphics[width=6.25cm]{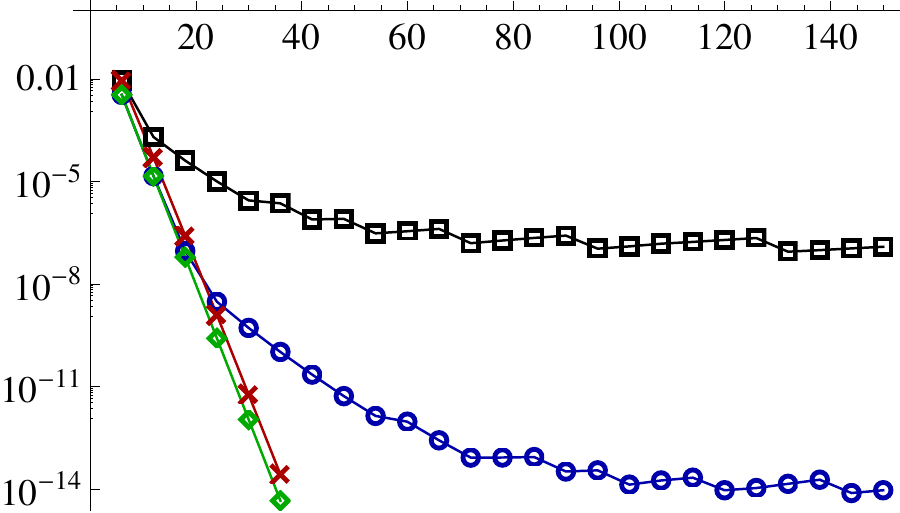}  & \hspace{1pc} &
\includegraphics[width=6.25cm]{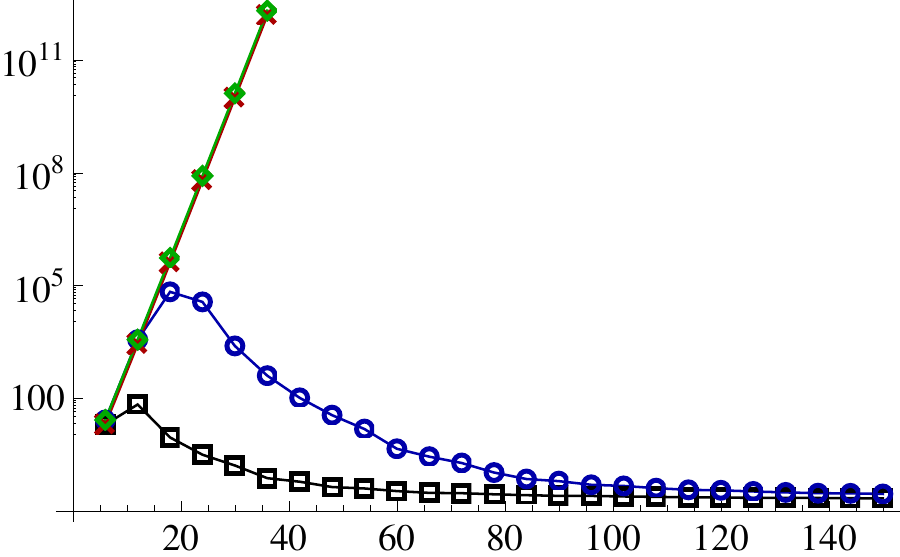}
\\
\includegraphics[width=6.25cm]{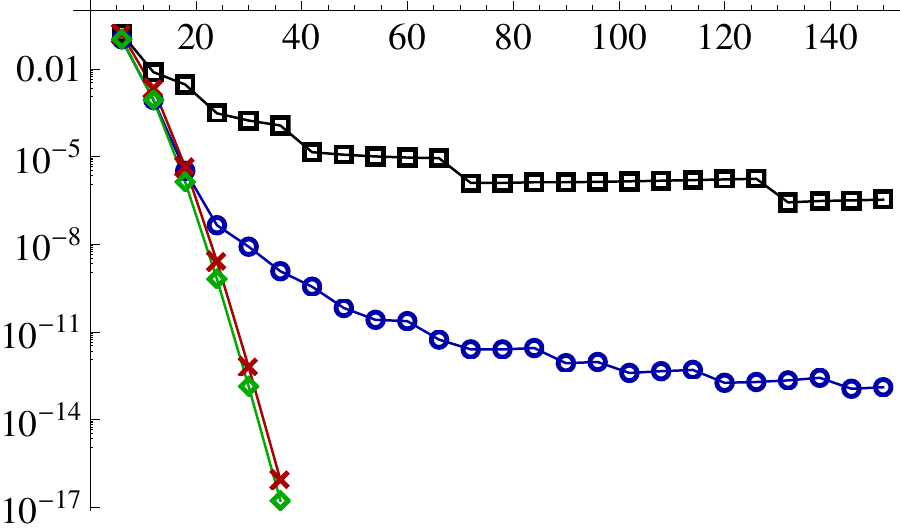}  & \hspace{1pc} &
\includegraphics[width=6.25cm]{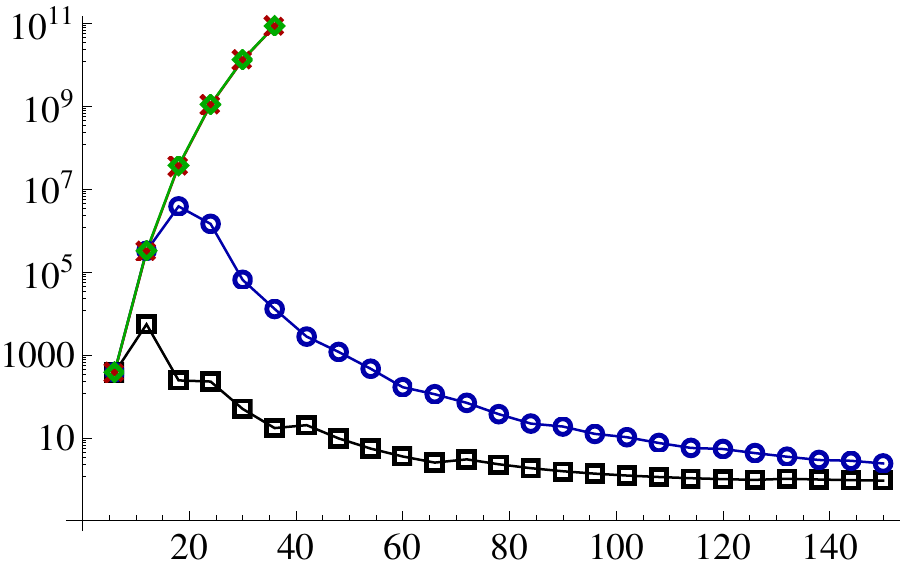}
\end{array}$
\caption{\small Comparison of the numerical continuous and
discrete FEs $G_N(f)$ and $\tilde{G}_{N}(f)$ (squares and
circles) and their exact counterparts $F_N(f)$ and
$\tilde{F}_N(f)$ (crosses and diamonds) for $T=2$. Left: the
uniform error $\| f - f_N \|_{\infty}$ against $N$.  Right: the
norm $\| a \|$ of the coefficient vector.  Top row: $f(x) =  \frac{1}{1+16 x^2}$.  Middle row: $f(x) =
\frac{1}{8-7 x}$.
Bottom row: $f(x) = 1+\frac{\cosh 40 x}{\cosh 40}$. }
\label{f:FnAppCoeff}
\end{center}
\end{figure}

This hypothesis is verified numerically in Figure \ref{f:FnAppCoeff} (all computations were carried out in \textit{Mathematica}, with additional precision used to compute the exact FEs and standard precision used otherwise).  Geometric convergence of the exact extension is replaced by slower, but still high-order convergence for sufficiently large $N$.  Note that the `breakpoint' occurs at roughly the same
value of $N$ regardless of the function being approximated. Moreover, the breakpoint
occurs at a larger value of $N$ for the discrete extension than
for the continuous extension.

These observations will
be established rigorously in the next section.  However, we now
make several further comments on Figure
\ref{f:FnAppCoeff}.  First, note that the breakdown of geometric
convergence is far less severe for the classical Runge function
$f(x) = \frac{1}{1+16 x^2}$ than for the functions $f(x) = \frac{1}{8-7
x}$ and $f(x) = 1 + \frac{\cosh 40 x}{\cosh
40}$. This can be explained by the behaviour of these functions near $x=\pm 1$.  The Runge function $f(x) = \frac{1}{1+16 x^2}$ is reasonably flat near $x = \pm 1$.  Hence it possesses extensions with high degrees of smoothness which do not grow large on the extended domain $[-T,T]$.  Conversely, the other two functions have boundary layers near $x=1$ (also $x=-1$ for the latter).  Therefore any smooth extension will be large on $[-T,T]$, and by Parseval's relation, the coefficient vectors corresponding to the Fourier series of this extension will also have large norm.

\begin{figure}[t]
\begin{center}
$\begin{array}{ccc}
\includegraphics[width=6.25cm]{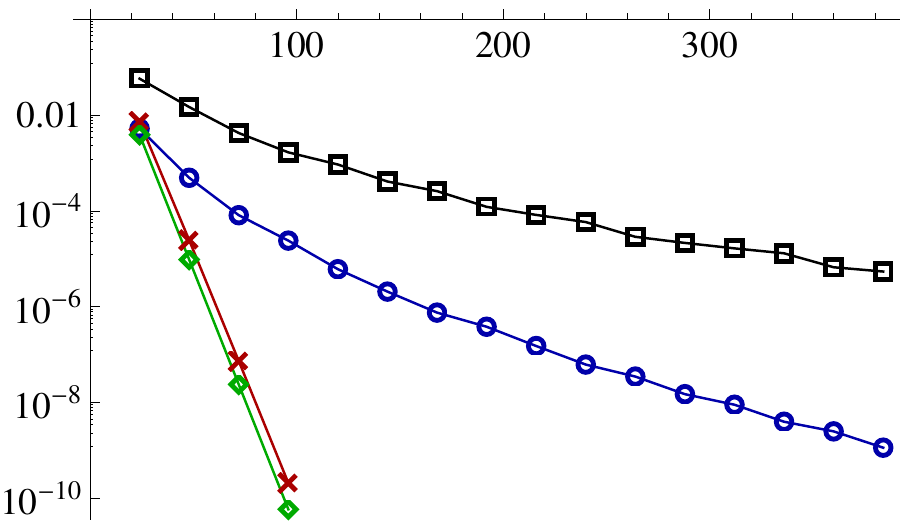}   & \hspace{1pc} &\includegraphics[width=6.25cm]{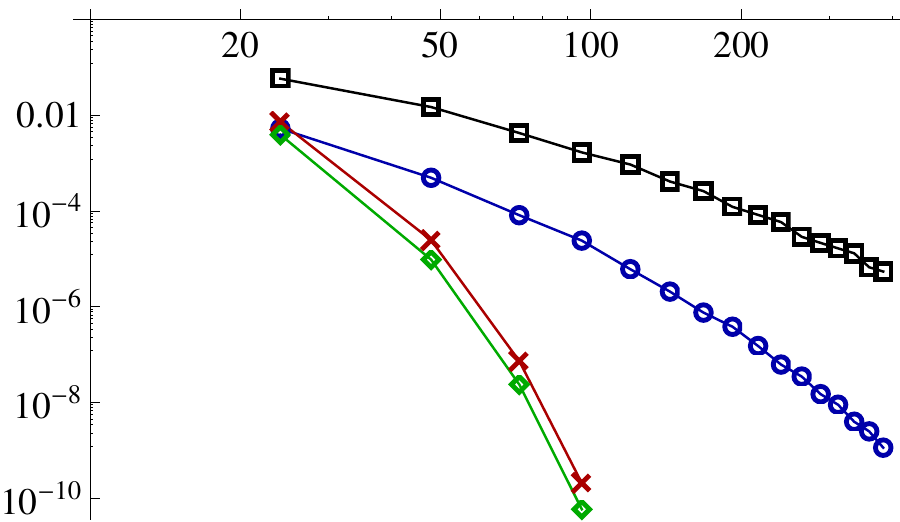}  \\
\end{array}$
\caption{\small Comparison of the numerical continuous and discrete FEs $G_N(f)$ and $\tilde{G}_{N}(f)$ (squares and circles) and their exact counterparts $F_N(f)$ and $\tilde{F}_N(f)$ (crosses and diamonds) for $T=2$ and $f(x) = \frac{1}{101-100 x}$.  Left: the uniform error in log scale.  Right: the uniform error in log-log scale. }  \label{f:FnAppSpectral}
\end{center}
\end{figure}

Second, although it is not apparent from Figure \ref{f:FnAppCoeff}
that the convergence rate beyond the breakpoint is truly superalgebraic, this is in fact the case.  This is confirmed by Figure \ref{f:FnAppSpectral}.  In the right-hand diagram we plot the error against $N$ in log-log scale.  The slight downward curve in the error indicates superalgebraic convergence.  Had the convergence rate been algebraic of fixed order then the error would have followed a straight line.

\subsection{Analysis of the numerical continuous and discrete FEs}\label{ss:numsolnanalysis}
We now wish to analyze the numerical extensions $G_N(f)$ and
$\tilde{G}_N(f)$.  Since the numerical solvers
used in environments such as \textit{Matlab} or
\textit{Mathematica} are difficult to analyze directly, we shall look at the result of solving $A a = b$ (or $\tilde A a
=\tilde b$) with a truncated singular value decomposition (SVD).
This represents an idealization of the numerical solver.  Indeed,
neither \textit{Matlab}'s $\backslash$ or \textit{Mathematica}'s
\texttt{LeastSquares} actually performs a truncated SVD.  However,
in practice, this simplification appears reasonable: numerical
experiments indicate that these
standard solvers give roughly the same approximation errors as the truncated
SVD with suitably small truncation parameter (typically $\epsilon
= 10^{-14}$).  We shall also assume throughout that the truncated
SVD is computed without error.  However, this also seems
fair: in experiments, we observe that the finite-precision SVD gives similar results to the
numerical solver whenever the tolerance is sufficiently small.

Suppose that $A$ (respectively $\tilde A$) has SVD $U S V^*$ with $S$ being the diagonal matrix of singular values.  Given a truncation parameter $\epsilon > 0$, we now consider the solution
\be{
\label{aeps}
a_{\epsilon} = V S^{\dag} U^* b,
}
where $S^\dag$ is the diagonal matrix with $n^{\rth}$ entry $1/\sigma_n$ if $\sigma_n > \epsilon$ and $0$ otherwise.  We write
\bes{
H_{N,\epsilon}(f) = \sum_{|n| \leq N} (a_{\epsilon})_n \phi_n,
}
for the corresponding FE.  Suppose that $ v_n \in \bbC^{2N+1}$ is the right singular vector of $A$ with singular value $\sigma_n$, and let
\bes{
\Phi_n = \sum_{|m| \leq N} (v_n)_m \phi_m \in \cG_N,
}
be the Fourier series corresponding to $v_n$.  Note that the functions $\Phi_n$ are orthonormal with respect to $\ip{\cdot}{\cdot}_{[-T,T]}$ and span $\cG_N$.  Also, if we define $\cG_{N,\epsilon} = \spn \{ \Phi_n : \sigma_n > \epsilon \} \subseteq \cG_N$, then we have $H_{N,\epsilon}(f) \in \cG_{N,\epsilon}$.

We now consider the cases of the continuous and discrete FEs separately.

\subsubsection{The continuous Fourier extension}\label{sss:SVDexact}

In this case, since $A$ is Hermitian and positive definite, the singular vectors $v_n$ are actually eigenvectors of $A$ with $A v_n = \sigma_n v_n$.  By definition, we have $\ip{\Phi_n}{\Phi_m} = (v_n)^* A v_m = \sigma_n \delta_{n,m}$, and therefore
\be{
\label{exactTSVD}
H_{N,\epsilon}(f) = \sum_{n : \sigma_n > \epsilon} \frac{1}{\sigma_n} \ip{f}{\Phi_n} \Phi_n.
}
Our main result is as follows:

\thm{ \label{t:exactSVD} Let $f \in \rL^2(-1,1)$ and suppose that
$H_{N,\epsilon}(f)$ is given by \R{exactTSVD}.  Then \be{
\label{exactTSVDerr} \| f - H_{N,\epsilon}(f) \| \leq \| f - \phi
\| + \sqrt{\epsilon} \| \phi \|_{[-T,T]},\quad \forall \phi \in
\cG_N, } and \be{ \label{exactTSVDcoeff} \| a_{\epsilon} \| = \|
H_{N,\epsilon}(f) \|_{[-T,T]} \leq \frac{1}{\sqrt{\epsilon}} \| f
- \phi \| + \| \phi \|_{[-T,T]},\quad \forall \phi \in \cG_N. } }
\prf{ The function $H_{N,\epsilon}(f)$ is the orthogonal
projection of $f$ onto $\cG_{N,\epsilon}$ with respect to
$\ip{\cdot}{\cdot}$.  Hence for any $\phi \in \cG_{N}$ we have $\|
f - H_{N,\epsilon}(f) \| \leq \| f - H_{N,\epsilon}(\phi)\| \leq
\| f - \phi \| + \| \phi - H_{N,\epsilon}(\phi) \|$.  Consider the
latter term.  Since $\phi \in \cG_N$, the observation that the functions $
\Phi_n$ are also orthonormal on $[-T,T]$ gives \bes{ \| \phi -
H_{N,\epsilon}(\phi) \|^2 = \left \| \sum_{n : \sigma_n <
\epsilon} \ip{\phi}{\Phi_n}_{[-T,T]} \Phi_n \right \|^2=  \sum_{n:
\sigma_n < \epsilon} \sigma_n |\ip{\phi}{\Phi_n}_{[-T,T]} |^2 \leq
\epsilon \| \phi \|^2_{[-T,T]}. } This yields \R{exactTSVDerr}.
For \R{exactTSVDcoeff} we first write $\| H_{N,\epsilon}(f)
\|_{[-T,T]} \leq \| H_{N,\epsilon}(f-\phi)  \|_{[-T,T]} + \|
H_{N,\epsilon}(\phi) \|_{[-T,T]}$.  By orthogonality, \bes{ \|
H_{N,\epsilon}(f-\phi)  \|^2_{[-T,T]} = \sum_{n: \sigma_n >
\epsilon} \frac{1}{\sigma^2_n} | \ip{f-\phi}{\Phi_n} |^2 \leq
\frac{1}{\epsilon} \sum_{n: \sigma_n > \epsilon}
\frac{1}{\sigma_n} | \ip{f-\phi}{\Phi_n} |^2 =  \frac{1}{\epsilon}
\| H_{N,\epsilon}(f-\phi) \|^2. } Since $H_{N,\epsilon}$ is an
orthogonal projection, we conclude that $\| H_{N,\epsilon}(f-\phi)
\|^2_{[-T,T]} \leq \frac{1}{\epsilon} \| f - \phi \|^2$, which
gives the first term in \R{exactTSVDcoeff}.  For the second, we
notice that \bes{ \| H_{N,\epsilon}(\phi) \|^2_{[-T,T]}=\sum_{n : \sigma_n > \epsilon} | \ip{\phi}{\Phi_n}_{[-T,T]} |^2 \leq \|
\phi \|^2_{[-T,T]}, } since $\phi \in \cG_N$. }
This theorem allows us to explain the behaviour of the numerical FE $G_N(f)$.
Suppose that $f$ is analytic in $\cD(\rho)$ and continuous on its boundary, where $\rho <
E(T)$ and $\cD(\rho)$ is as in Theorem \ref{t:expconv}. Set
$\phi = F_N(f)$ in \R{exactTSVDerr}, where $F_N(f)$ is the exact
continuous FE. Then Theorems \ref{t:expconv} and
\ref{t:coefficients} give \be{\label{one} \| f -
H_{N,\epsilon}(f) \| \leq c_f \left ( 1 + \sqrt{\epsilon} E(T)^N
\right ) \rho^{-N}. } For small $N$, the first term in the
brackets dominates, and we see geometric convergence of
$H_{N,\epsilon}(f)$, and therefore also
$G_N(f)$, at rate $\rho$.  Convergence continues as such until the breakpoint
\be{ \label{exactbreakpt}
N_0 = N_0(\epsilon,T) : = - \frac{\log \epsilon}{2 \log E(T)},
}
at which point the second term dominates and the bound begins to increase.  On the other hand, Proposition \ref{p:appcoeff} establishes the existence
of functions $\phi \in \cG_N$ with bounded coefficients which
approximate $f$ to any given algebraic order.  Substituting
such a function $\phi$ into \R{exactTSVDerr} gives \be{ \| f -
H_{N,\epsilon}(f) \| \leq c_k(T) \left ( N^{-k} + \sqrt{\epsilon}
\right ) \| f \|_{\rH^k(-1,1)},\quad \forall N, k \in \bbN. }
Therefore, once $N
> N_0(\epsilon,T)$ we expect at least superalgebraic convergence of
$H_{N,\epsilon}(f)$  down to a maximal
achievable accuracy of order $\sqrt{\epsilon} \| f \|$.  Note that at the breakpoint $N = N_0$, the error satisfies
\be{
\label{four}
\| f - H_{N_0,\epsilon} (f) \| \leq 2 c_f (\sqrt{\epsilon})^{d_f},\qquad d_f = \frac{\log \rho}{\log E(T)} \in (0,1].
}
If $f$ is analytic in $\cD(E(T))$, and if $c_f = \max_{x \in \cD(\rho)} | f(x) |$ is not too large, then $f$ is already approximated to order $\sqrt{\epsilon}$ accuracy at this point.  It is only in those cases where either $\rho < E(T)$ or where $c_f$ is large (or both) that one sees the second phase of superalgebraic convergence.

Theorem \ref{t:coefficients} also explains the behaviour of the coefficient norm $\| a_{\epsilon} \|$.  Observe that breakpoint $N_0(\epsilon,T)$ is (up to a small constant) the largest $N$ for which all singular values of $A$ are included in its truncated SVD (see Theorem \ref{t:exactCondNumb}).  Thus, when $N < N_0(\epsilon,T)$, we have $H_{N,\epsilon}(f) = F_{N}(f)$, and Theorem \ref{t:coefficients} indicates exponential growth of $\| a_{\epsilon} \|$.  On the other hand, once $N > N_0(\epsilon,T)$, we use \R{exactTSVDcoeff} to obtain
\bes{
\| a_{\epsilon} \| \leq c_k(T) \left (  N^{-k} / \sqrt{\epsilon} + 1 \right ) \| f \|_{\rH^k(-1,1)},\quad \forall N, k \in \bbN.
}
In particular, for $N > N_{0}(\epsilon,T)$, we expect decay of $\| a_{\epsilon} \|$ down from its maximal value at $N = N_{0}(\epsilon,T)$.

\begin{figure}[t]
\begin{center}
$\begin{array}{ccc}
\includegraphics[width=6.25cm]{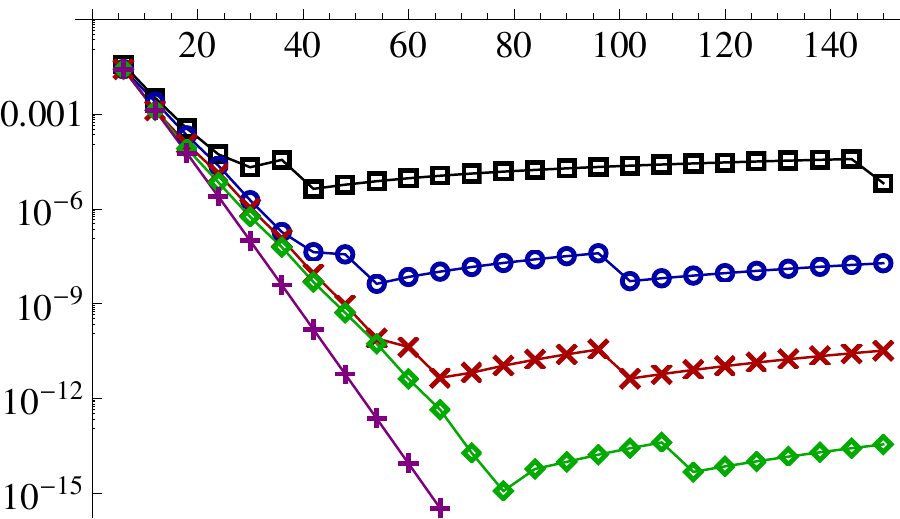} & \hspace{1pc} & \includegraphics[width=6.25cm]{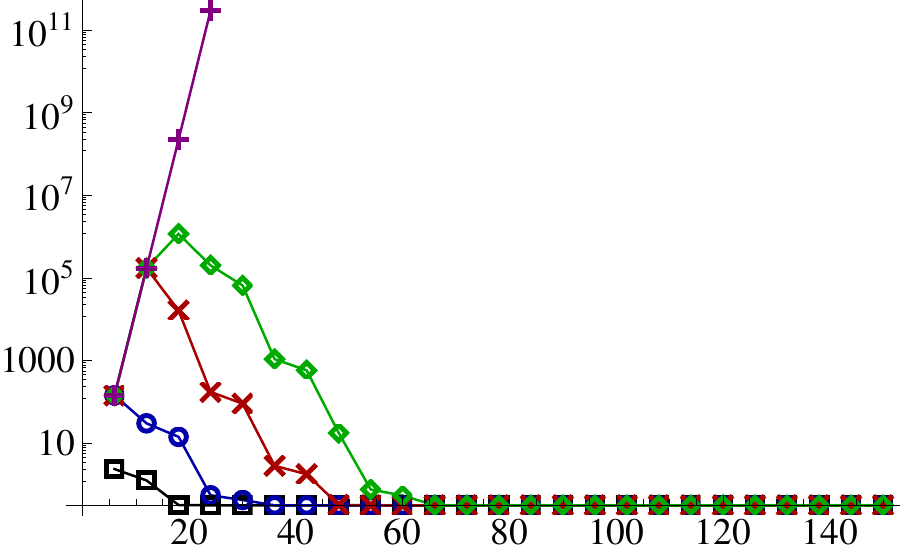} \\
\includegraphics[width=6.25cm]{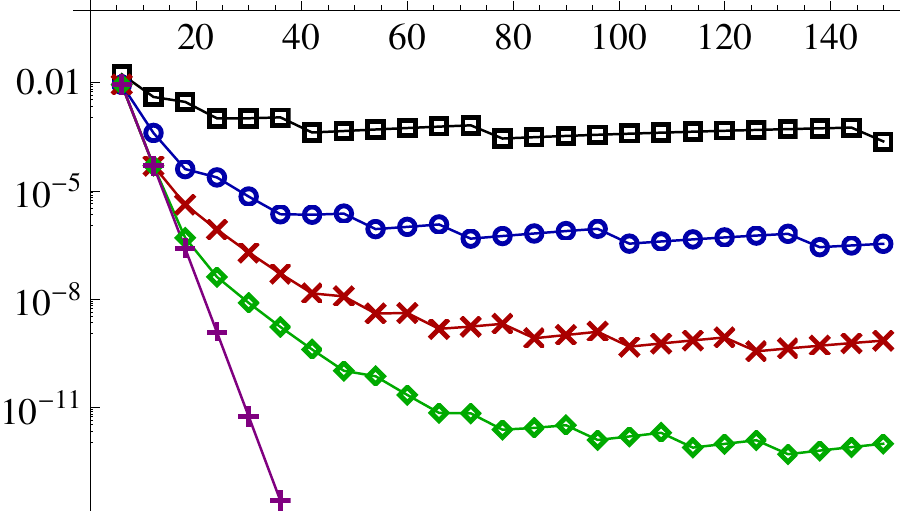} & \hspace{1pc} & \includegraphics[width=6.25cm]{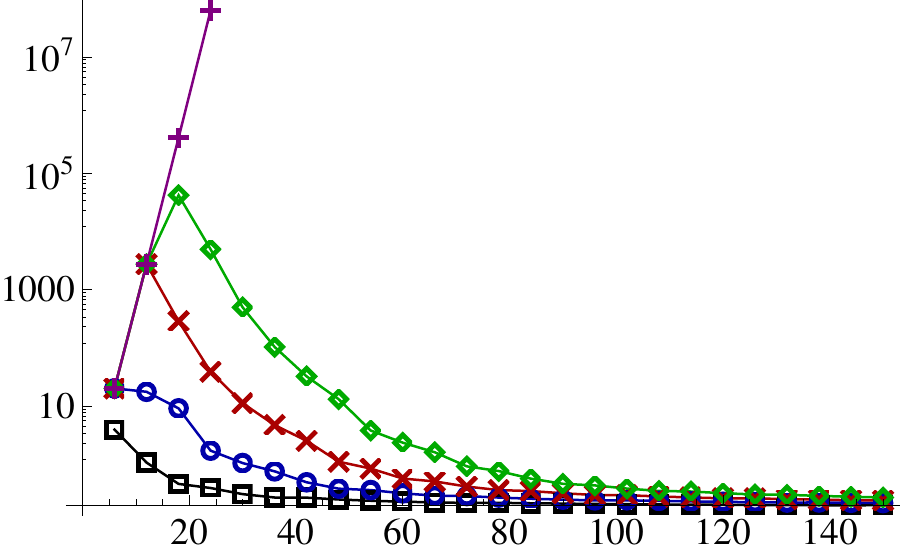} \\
\includegraphics[width=6.25cm]{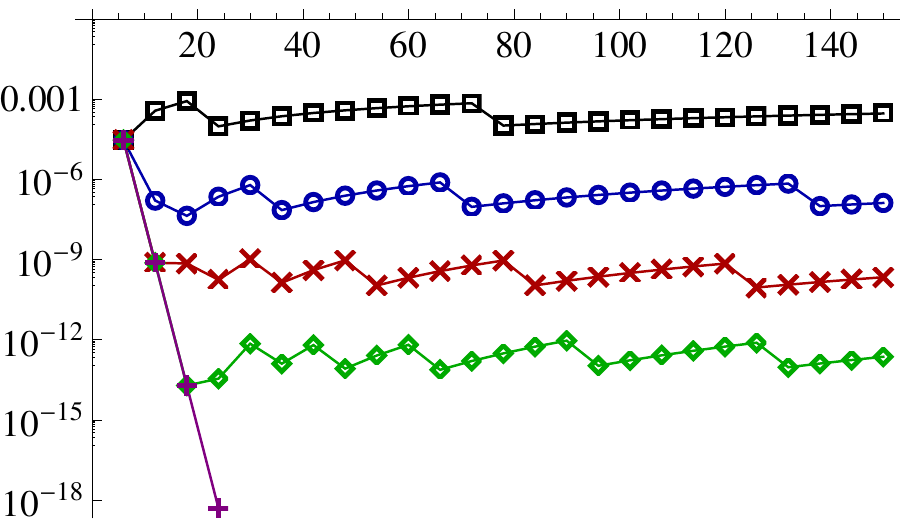} & \hspace{1pc} & \includegraphics[width=6.25cm]{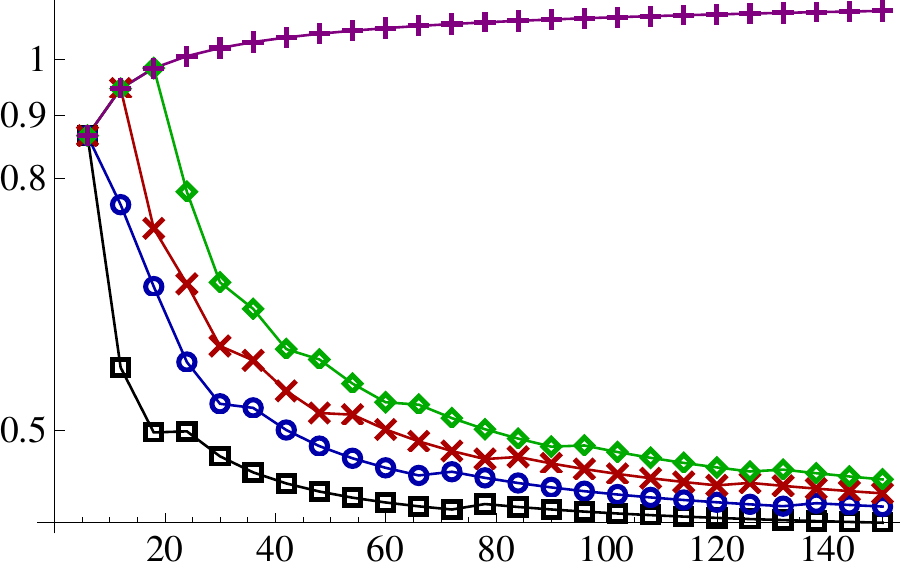}
\end{array}$
\caption{\small Error (left) and coefficient norm (right) against
$N$ for the continuous FE with $T=2$, where $f(x) = \frac{1}{1+16 x^2}$ (top row), $f(x)= \frac{1}{8-7 x} $ (middle row) and $f(x) = x$ (bottom row).  Squares, circles, crosses and diamonds correspond to the truncated SVD extension $H_{N,\epsilon}(f)$ with
$\epsilon = 10^{-6},10^{-12},10^{-18},10^{-24}$ respectively, and pluses correspond to the exact extension $F_N(f)$.}
\label{f:ManSVDexact}
\end{center}
\end{figure}

This analysis is corroborated in Figure \ref{f:ManSVDexact}, where we plot the error and coefficient norm for the truncated SVD extension for various test functions.  Note that the maximal achievable accuracy in all cases is order $\sqrt{\epsilon}$, consistently with our analysis.  Moreover, for the meromorphic functions $f(x) = \frac{1}{1+16 x^2}$ and $f(x) = \frac{1}{8-7 x}$ we see initial geometric convergence followed by slower convergence after $N_0$, again as our analysis predicts.  The qualitative difference in convergence for these functions in the regime $N > N_0$ is due to the contrasting behaviour of their derivatives (recall the discussion in \S \ref{ss:magnitude}).  On the other hand, the convergence effectively stops at $N_0$ for $f(x) = x$, since this function has small constant $c_f$ and is therefore already resolved down to order $\sqrt{\epsilon}$ when $N = N_0$.

Since $N_0(10^{-6},2) \approx 4$, $N_0(10^{-12},2) \approx 8$, $N_0(10^{-18},2) \approx 12,$ and $N_0(10^{-24},2) \approx 16$, Figure \ref{f:ManSVDexact} also confirms the expression \R{exactbreakpt} for the breakpoint in convergence.  In particular, the breakpoint is independent of the function being approximated.  This latter observation is unsurprising.  As noted, $N_0(\epsilon,T)$ is the largest value of $N$ for which $H_{N,\epsilon}(f)$ coincides with $F_N(f)$.  Beyond this point, $H_{N,\epsilon}(f)$ will not typically agree with $F_N(f)$, and thus we cannot expect further geometric convergence in general.  Note that our analysis does not rule out geometric convergence for $N > N_0$.  There may well be certain functions for which this occurs.  However, extensive numerical tests suggest that in most cases, one sees only superalgebraic convergence in this regime, and indeed, this is all that we have proved.

\rem{
\label{r:singularfunctions}
At first sight, it may appear counterintuitive that one can still obtain good accuracy when excluding all singular values below a certain tolerance.  However, recall that we are not interested in the accuracy of computing $a$, but rather the accuracy of $F_N(f)$ on the domain $[-1,1]$.  Since the $n^{\rth}$ singular value $\sigma_n$ is equal to $\| \Phi_n \|^2 / \| \Phi_n \|^2_{[-T,T]}$, the functions $\Phi_n$ excluded from $H_{N,\epsilon}(f)$ are precisely those for which $\| \Phi_n \|^2 < \epsilon \| \Phi_n \|^2_{[-T,T]}$.  In other words, they have little effect on $F_{N}(f)$ in $[-1,1]$.

In Figure \ref{f:SingFnPlot} we plot the functions $\Phi_n$ for several $n$.  Note that these functions are precisely the discrete prolate spheroidal wavefunctions of Slepian \cite{SlepianV}.  As predicted, when $n$ is small, the function $\Phi_n$ is large in $[-1,1]$ and small in $[-T,T] \backslash [-1,1]$.  When $n$ is in the transition region ($n \approx 2N/T$, see \S \ref{ss:ASVD}), the function $\Phi_n$ is roughly of equal magnitude in both regions, and for $n \approx 2N$, $\Phi_n$ is much smaller in $[-1,1]$ than on $[-T,T]$.  Note also that $\Phi_n$ is increasingly oscillatory in $[-1,1]$ as $n$ increases, and decreasingly oscillatory in $[-T,T] \backslash [-1,1]$.  This follows from the fact that $\Phi_n$ has precisely $n$ zeroes in $[-1,1]$ and $2N-n$ zeroes in $[-T,T] \backslash [-1,1]$ \cite{SlepianV}.  Such behaviour also implies that any `nice' function will eventually be well approximated by functions $\Phi_n$ corresponding to `nice' eigenvalues, as expected.

\begin{figure}
\begin{center}
$\begin{array}{ccc}
\includegraphics[width=4.75cm]{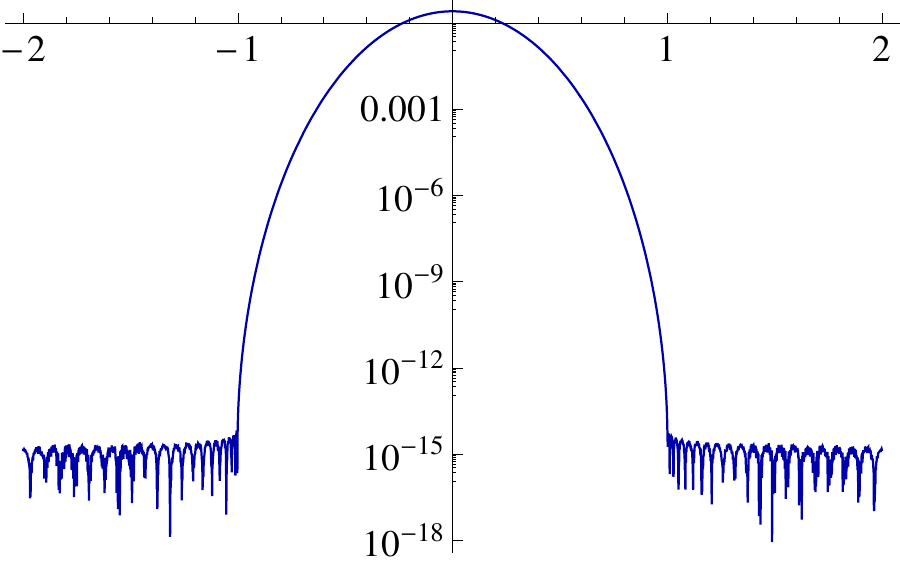}  & \includegraphics[width=4.75cm]{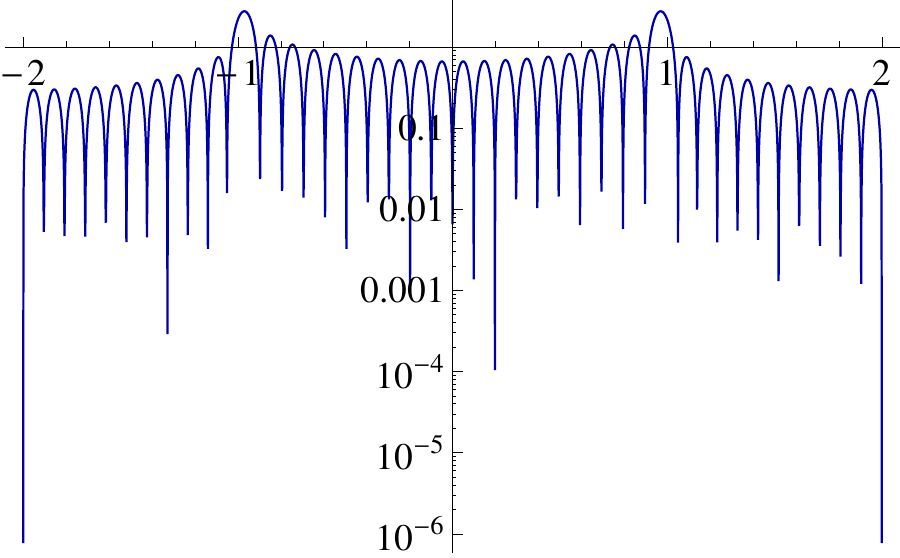} &
 \includegraphics[width=4.75cm]{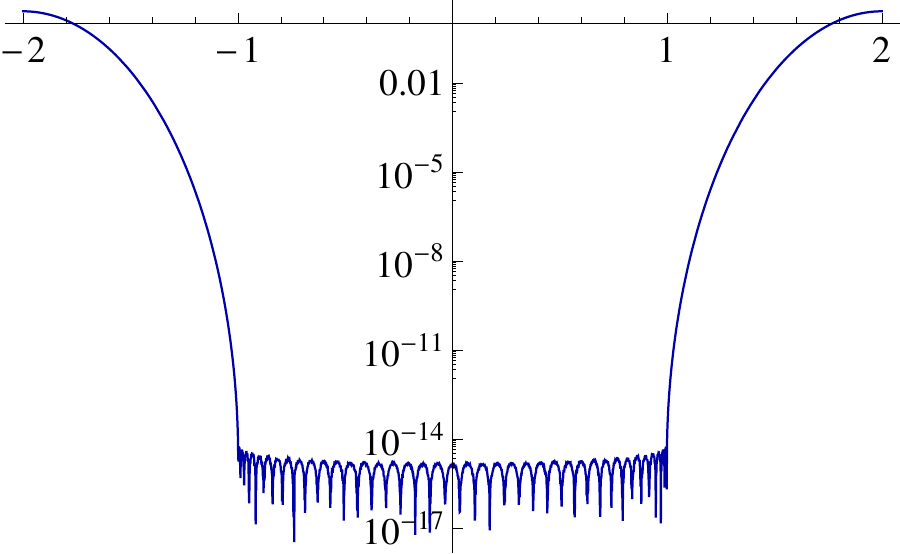}
\end{array}$
\caption{\small The SVD functions $|\Phi_n(x)|$ for $n=0$,
$n=20$ and $n=40$, where $N = 20$ and $T=2$.} \label{f:SingFnPlot}
\end{center}
\end{figure}
}

\subsubsection{The discrete Fourier extension}\label{sss:discSVDanalysis}
In this case, we have $( \Phi_n , \Phi_m)_N = \sigma^2_n \delta_{n,m}$, where
\bes{
(f,g)_N = \frac{\pi}{N+1} \sum^{N}_{n=-N-1}  f(x_n) \overline{g(x_n)},
}
is the discrete inner product corresponding to the quadrature nodes $\{ x_n \}^{N}_{n=-N-1}$.  Therefore
\be{
\label{discTSVD}
\tilde{H}_{N,\epsilon}(f) = \sum_{n : \sigma_n > \epsilon} \frac{1}{\sigma^2_n} (f,\Phi_n)_N \Phi_n \in \cG'_{N,\epsilon} : = \spn \left \{ \Phi_n : \sigma_n > \epsilon \right \},
}
is the orthogonal projection of $f$ onto $\cG'_{N,\epsilon}$ with respect to the discrete inner product $(\cdot,\cdot)_{N}$.

\thm{
Let $f \in \rL^\infty(-1,1)$ and $\tilde{H}_{N,\epsilon}(f)$ be given by \R{discTSVD}.  Then\be{
\label{discTSVDerr}
\| f - \tilde{H}_{N,\epsilon}(f) \|_{W} \leq \| f - \phi \|_W + \sqrt{2 \pi Q(N;\epsilon)} \| f - \phi \|_{\infty} + \epsilon \| \phi \|_{[-T,T]},\quad \forall \phi \in \cG_N,
}
and
\be{
\label{discTSVDcoeff}
\| a_{\epsilon} \|  = \| \tilde{H}_{N,\epsilon}(f) \|_{[-T,T]} \leq \frac{1}{\epsilon} \sqrt{2 \pi Q(N;\epsilon)} \| f - \phi \|_{\infty} + \| \phi \|_{[-T,T]},\quad \forall \phi \in \cG_N,
}
where $Q(N;\epsilon) = | \{ n : \sigma_n > \epsilon \} | \leq 2(N+1)$ and $W$ is the weight function of Lemma \ref{l:weightedOpt}.
}
\prf{
By the triangle inequality,
\bes{
\| f - \tilde{H}_{N,\epsilon}(f) \|_W \leq \| f - \phi \|_W + \| \phi - \tilde{H}_{N,\epsilon}(\phi) \|_W + \| \tilde{H}_{N,\epsilon}(f-\phi) \|_W,\quad \forall \phi \in \cG'_N.
}
Consider the second term.  Since $\phi \in \cG'_{N}$, and the quadrature is exact on $\cG'_N$, we have
\bes{
\| \phi - \tilde{H}_{N,\epsilon}(\phi) \|^2_W = (\phi - \tilde{H}_{N,\epsilon}(\phi),\phi - \tilde{H}_{N,\epsilon}(\phi))_N = \sum_{n : \sigma_n < \epsilon} \sigma^2_n |\ip{\phi}{\Phi_n}_{[-T,T]}|^2 \leq \epsilon^2 \| \phi \|^2_{[-T,T]}.
}
For the third term, let $g$ be arbitrary.  Then $ (\tilde{H}_{N,\epsilon}(g),\tilde{H}_{N,\epsilon}(g))_N = \sum_{n : \sigma_n > \epsilon} \frac{1}{\sigma^2_n} | (g,\Phi_n)_N |^2$.  Hence
\be{
\label{step1}
\| \tilde{H}_{N,\epsilon}(g) \|^2_{W} =(\tilde{H}_{N,\epsilon}(g),\tilde{H}_{N,\epsilon}(g))_N \leq (g,g)_N \sum_{n : \sigma_n > \epsilon} \frac{1}{\sigma^2_n} (\Phi_n,\Phi_n)_N =(g,g)_N Q(N;\epsilon),
}
since $(\Phi_n,\Phi_n)_N  = \sigma^2_n$.  It is straightforward to show that $(g,g)_N \leq 2 \pi \| g \|^2_{\infty}$.  Setting $g = f - \phi$ now gives the corresponding term in \R{discTSVDerr}, and completes its proof.  For \R{discTSVDcoeff}, we proceed as in the proof of Theorem \ref{t:exactSVD}.  Note that
\be{
\label{step2}
\| \tilde{H}_{N,\epsilon}(g) \|^2_{[-T,T]} = \sum_{n : \sigma_n > \epsilon} \frac{1}{\sigma^4_n} | (g,\Phi_n)_N |^2 \leq \frac{1}{\epsilon^2} \| \tilde{H}_{N,\epsilon}(g) \|^2_W,
}
for any $g \in  \rL^\infty(-1,1)$.  Also,
\be{
\label{step3}
\| \tilde{H}_{N,\epsilon}(\phi) \|_{[-T,T]} \leq \| \phi \|_{[-T,T]}, \quad \phi \in \cG_N.
}
The result now follows by writing $\| \tilde{H}_{N,\epsilon}(f) \|_{[-T,T]} \leq \| \tilde{H}_{N,\epsilon}(f-\phi) \|_{[-T,T]} +\| \tilde{H}_{N,\epsilon}(\phi) \|_{[-T,T]}$ and using \R{step1}--\R{step3}.
}

\begin{figure}[t]
\begin{center}
$\begin{array}{ccc}
\includegraphics[width=6.25cm]{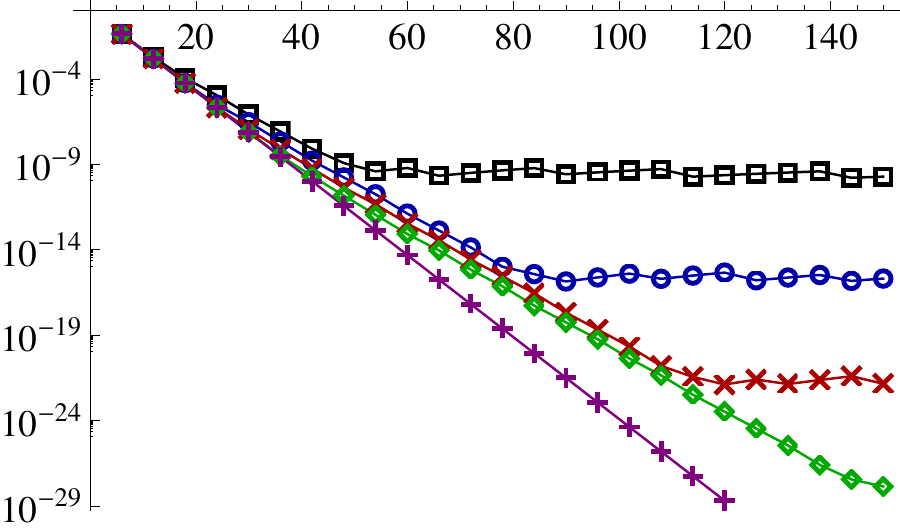} & \hspace{1pc} & \includegraphics[width=6.25cm]{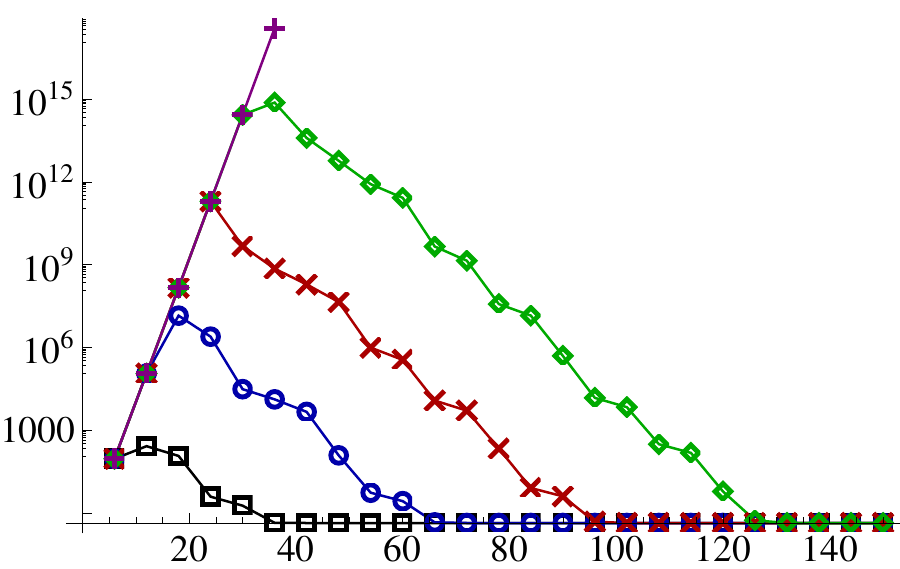} \\
\includegraphics[width=6.25cm]{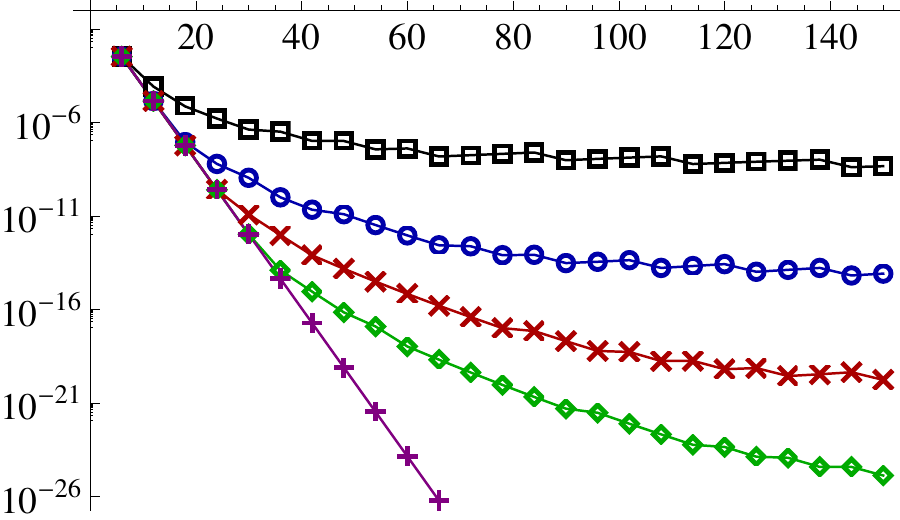} & \hspace{1pc} & \includegraphics[width=6.25cm]{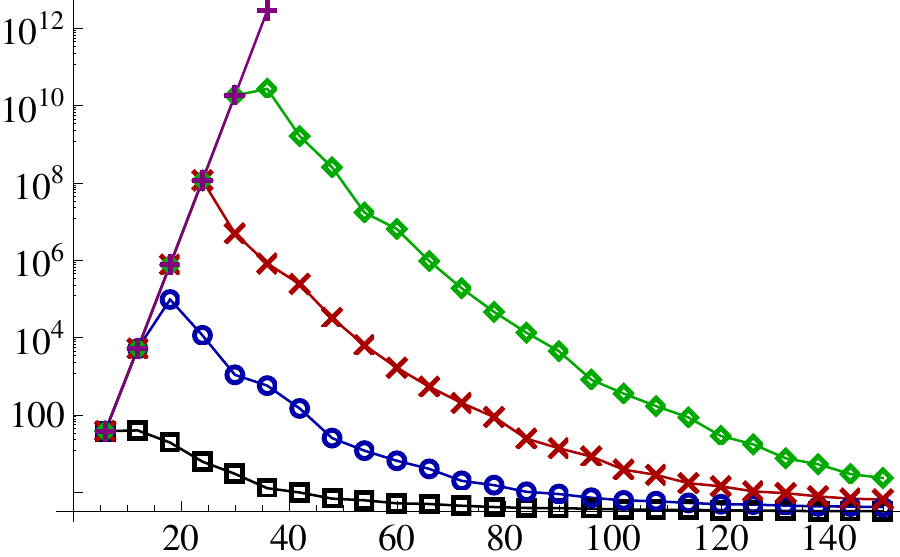} \\
\includegraphics[width=6.25cm]{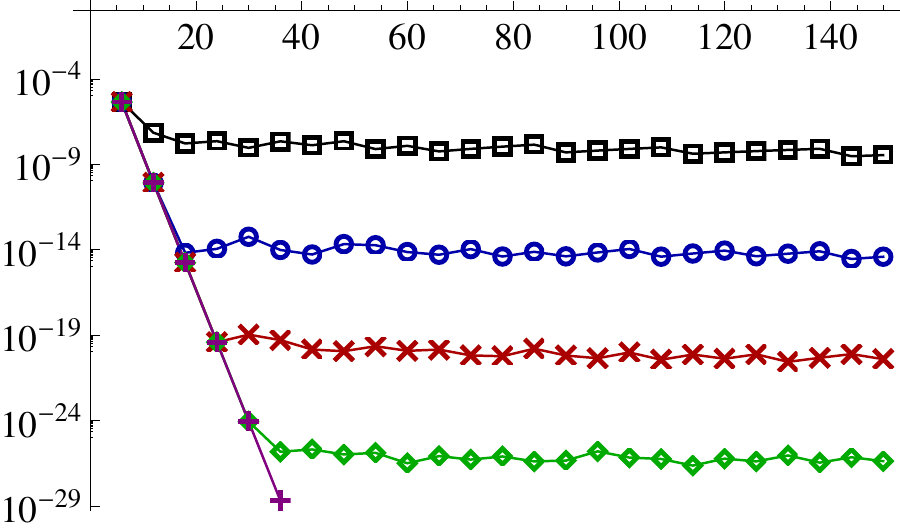} & \hspace{1pc} & \includegraphics[width=6.25cm]{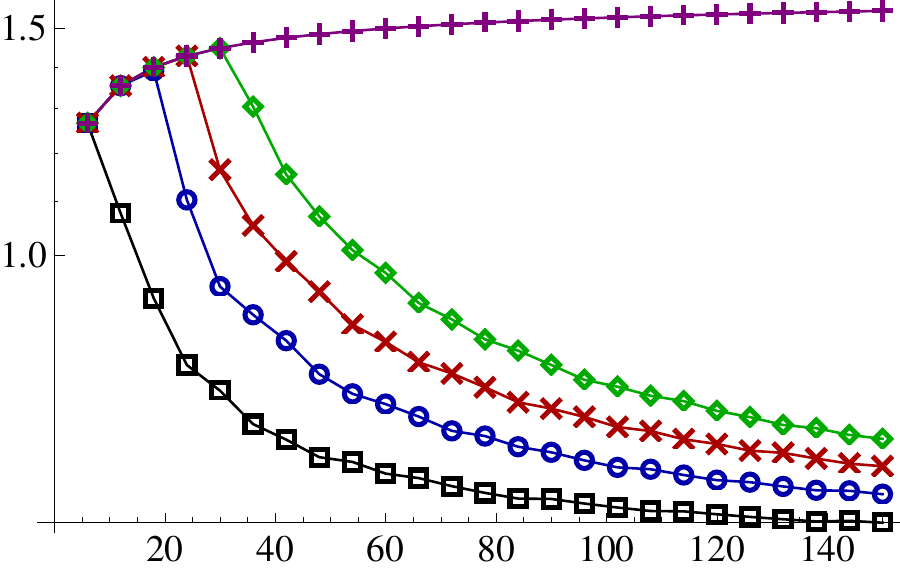}
\end{array}$
\caption{\small  Error (left) and coefficient norm (right) against
$N$ for the discrete FE with $T=2$, where $f(x) = \frac{1}{1+16 x^2}$ (top row), $f(x) = \frac{1}{8-7 x}$ (middle row) and $f(x) = x$ (bottom row).  Squares, circles, crosses and diamonds correspond to the truncated SVD extension $H_{N,\epsilon}(f)$ with
$\epsilon = 10^{-6},10^{-12},10^{-18},10^{-24}$ respectively, and pluses correspond to the exact extension $F_N(f)$.} \label{f:ManSVDdisc}
\end{center}
\end{figure}

As with the continuous FE, this theorem allows us to
analyze the numerical discrete extension $\tilde{G}_N(f)$.  Once more we deduce geometric
convergence in $N$ up to the function-independent breakpoint
\be{
\label{N1def}
N_1(T;\epsilon) : = - \frac{\log \epsilon}{\log E(T)} \equiv 2 N_0(T;\epsilon),
}
with superalgebraic convergence beyond this point.    These conclusions are confirmed in Figure \ref{f:ManSVDdisc}.  Note, however,
two key differences between the continuous and discrete FE.  First, the bound \R{discTSVDerr} involves
$\epsilon$, as opposed to $\sqrt{\epsilon}$, meaning that we
expect convergence of $\tilde{G}_N(f)$ down to close to machine precision.
Second, the breakpoint $N_{1}(T;\epsilon)$ is precisely twice
$N_0(T;\epsilon)$.  Hence, the regime of geometric convergence
of $\tilde{G}_N(f)$ is exactly twice as large as that of the
continuous FE.   These observations are in close agreement with the behaviour seen in the numerical examples in \S \ref{ss:magnitude}.

\subsection{The condition numbers of the numerical continuous and discrete FEs}\label{ss:stabfn}
Having analyzed the convergence of the numerical FE---and in particular, established 5.\ of \S \ref{s:introduction}---we next address its condition number.  Once more, we do this by considering the extensions $H_{N,\epsilon}$ and $\tilde{H}_{N,\epsilon}$:

\thm{ \label{t:stabfn}
Let $H_{N,\epsilon}$ be the continuous truncated SVD FE given by \R{exactTSVD}.  Then
 \bes{
 \kappa(H_{N,\epsilon})=1/ \min \{ \sqrt{\sigma_n} : \sigma_n > \epsilon \} \leq
\min \left \{ 1/\sqrt{\epsilon} , c(T) N^{\frac32} E(T)^N
\right \},\quad N \in \bbN,\  \epsilon > 0,
}
where $c(T)$ is a positive constant independent of $N$.  Conversely, if $\tilde{H}_{N,\epsilon}$ is the discrete extension \R{discTSVD}, then $\kappa( \tilde{H}_{N,\epsilon})=
1$ for all $N \in \bbN$ and $\epsilon > 0$.
}
\prf{
The proof of the equalities is similar to that of Lemma \ref{l:ExactCondNumb} with $A$ and $\tilde{A}$ replaced by their truncated SVD versions.  The upper bound for $\kappa(H_{N,\epsilon})$ is a consequence of Theorem \ref{t:exactCondNumb}.
}

This theorem, which establishes 3.\ of \S \ref{s:introduction}, has some interesting consequences.  First, the discrete FE is perfectly stable.  On the other hand, the numerical continuous FE is far from stable.  The condition number grows exponentially fast at rate $E(T)$ until it reaches $1/\sqrt{\epsilon}$, where $\epsilon$ is the truncation parameter in the SVD.  Thus, with the continuous FE, we may see perturbations being magnified by a factor of $1/\sqrt{\epsilon_{\mathrm{mach}}} \approx 10^8$ in practice.  

Note that $G_N$ and $\tilde{G}_N$ are both substantially better conditioned than the corresponding coefficient mappings.  The explanation for this difference comes from Remark \ref{r:singularfunctions}.  A perturbation $\eta$ in the input vector $b$ gives large errors in the FE coefficients if $\eta$ has a significant component in the direction of a singular vector $v_n$ associated with a small singular value $\sigma_n$.  However, since the corresponding function $\Phi_n$ is small on $[-1,1]$, this error is substantially reduced (in the case of the
continuous FE) or cancelled out altogether (for the discrete FE) in the resulting extension.

Another implication of Theorem \ref{t:stabfn} is the following: varying $T$ has no substantial effect on stability.  Although the condition number of the FE matrices depends on $T$ (recall Theorems \ref{t:exactCondNumb} and \ref{t:discCondNumb}), as does the condition number of the exact continuous FE (see Lemma \ref{l:ExactCondNumb}), the condition numbers of the numerical mappings $\tilde{G}_N$ and, for all large $N$, $G_N$ are actually independent of this parameter.

It is important to confirm that the results of this theorem on the condition number of the truncated SVD extensions predict the behaviour of the numerical extensions $G_N$ and $\tilde{G}_N$.  It is easiest to do this by computing upper bounds for $\kappa(G_N)$ and $\kappa(\tilde{G}_N)$.  Let $\{ e_n \}^{2N+1}_{n=1}$ be the standard basis for $\bbC^{2N+1}$.  Then a simple argument gives that
\be{
\| G_N(b) \| \leq \| b \| \sqrt{\sum^{2N+1}_{n=1} \| G_N(e_n) \|^2 },\quad \forall b \in \bbC^{2N+1}, 
}
and therefore
\be{
\label{KGN_def}
\kappa(G_N) \leq K(G_N) : = \sqrt{\sum^{2N+1}_{n=1} \| G_N(e_n) \|^2 }.
}
We define the upper bound $K(\tilde{G}_N)$ in a similar manner:
\bes{
\kappa(\tilde{G}_N) \leq K(\tilde{G}_N) : = \sqrt{\sum^{2N+2}_{n=1} \| \tilde{G}_N(e_n) \|^2_W}.
}
In Table \ref{t:FEStabFn} we show $K(G_N)$ and $K(\tilde{G}_N)$ for various choices of $N$.  As we see, the discrete FE is extremely stable: not only is there no blowup in $N$, but the value of $K(\tilde{G}_N)$ is also close to one in magnitude.  For the continuous extension, we see that $K(G_N) \approx 5 \times 10^{6} =1/\sqrt{\epsilon}$, where $\epsilon = 2.5 \times 10^{-13}$.  This behaviour is in good agreement with Theorem \ref{t:stabfn}.

\begin{table}
\begin{center}
\begin{tabular}{|c||c|c|c|c|c|}
\hline
$N$ & 40 & 80 & 120 & 160 & 200\\
\hline
$K(G_N)$ & $4.93 \times10^6$ & $4.22 \times 10^6$ & $3.30 \times 10^6$ & $3.82 \times 10^6$ & $5.28\times 10^6$ \\ \hline
$K(\tilde{G}_N)$ & $8.00\times10^0$ & $1.04\times10^1$ & $1.23\times10^1$ & $1.39\times10^1$ & $1.53\times10^1$ \\ \hline
\end{tabular}
\caption{The functions $K(G_N)$ and $K(\tilde{G}_N)$ for $T=2$.}\label{t:FEStabFn}
\end{center}
\end{table}

\begin{figure}
\begin{center}
$\begin{array}{ccc}
\includegraphics[width=6.25cm]{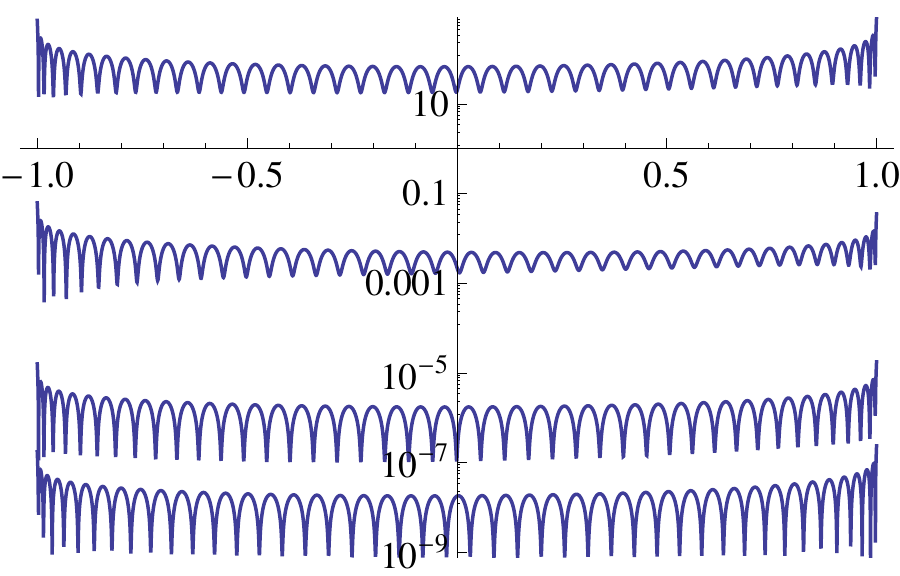}  & \hspace{1pc} & \includegraphics[width=6.25cm]{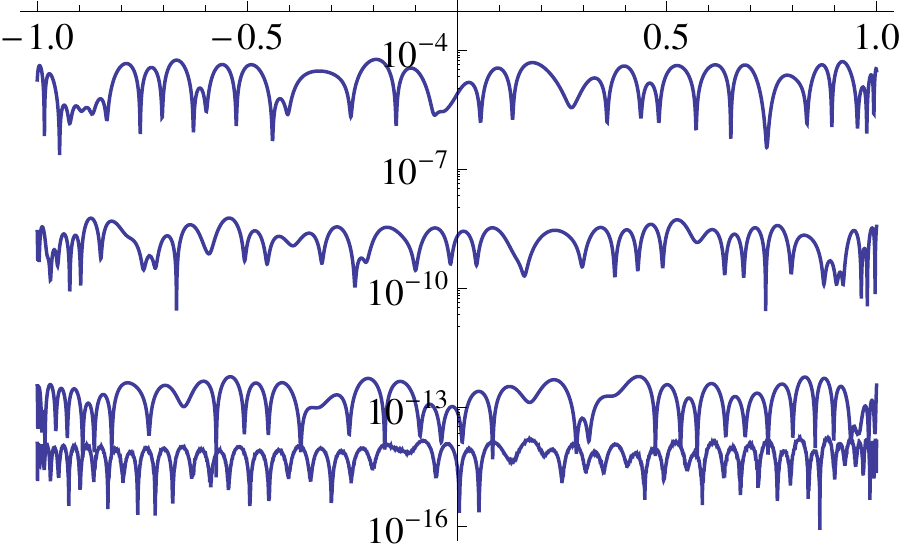}
\end{array}$
\caption{\small The error $| f(x) - f_N(x) |$ against $x$, where
$f_N = G_N(f)$ (left) or $f_N = \tilde{G}_N(f)$ (right), for $N =
30$, $T=2$ and $f(x) = \E^x$, with noise at amplitudes $\delta =
10^{-4},10^{-8},10^{-12},0$.}  \label{f:Noisy}
\end{center}
\end{figure}

The difference in stability between the continuous and discrete FEs is highlighted in Figure \ref{f:Noisy}.  Here we perturbed the right-hand side $b$ of the function $f(x) = \E^x$ by noise of magnitude $\delta$, and then computed its FE.  As is evident, the discrete extension approximates $f$ to an error of magnitude roughly $\delta$, whereas for the continuous extension the error is of magnitude $\approx 10^{6} \delta$, as predicted by Table \ref{t:FEStabFn}.

\section{Fourier extensions from equispaced data}\label{s:equispacedI}
We now turn our attention to the problem of computing FEs when only equispaced data is prescribed.  As discussed in \S \ref{s:introduction}, a theorem of Platte, Trefethen \& Kuijlaars \cite{TrefPlatteIllCond} states that any exponentially-convergent method for this problem must also be exponentially ill-conditioned (see \S \ref{ss:PTKrelation} for the precise result).  However, as we show in this section, FEs give rise to a method, the so-called \textit{equispaced Fourier extension}, that allows this barrier to be circumvented to a substantial extent.  Namely, it achieves rapid convergence in a numerically stable manner.  

\subsection{The equispaced Fourier extension}\label{ss:equiFEdef}
Let \be{ \label{equinodes} x_{n} = \frac{n}{M},\quad
n=-M,\ldots,M, } be a set of $2M+1$ equispaced points in $[-1,1]$,
where $M \geq N$.  We define the \textit{equispaced Fourier
extension} of a function $f \in \rL^\infty[-1,1]$ by \be{
\label{equiL2} F_{N,M}(f) : =  \underset{\phi \in
\cG_N}{\operatorname{argmin}}  \sum_{|n| \leq M} | f(x_n) -
\phi(x_n) |^2. }
If $F_{N,M}(f) = \sum_{|n| \leq N}
a_n \phi_n$, then the vector $a = (a_{-N},\ldots,a_N)^{\top}$ is
the least squares solution to $\bar{A} a \approx \bar{b}$,
where $\bar{A} \in \bbC^{(2M+1) \times (2N+1)}$ has $(n,m)^{\rth}$
entry $\frac{1}{\sqrt{M+1/2}} \phi_m(x_n)$ and $\bar{b}$ has $n^{\rth}$ entry $\frac{1}{\sqrt{M+1/2}} f(x_n)$.

Note that $F_{N,M}(f)$, as defined by \R{equiL2}, is (up to minor changes of parameters/notation) identical to the extensions considered in the previous papers \cite{BoydFourCont,BoydRunge,brunoFEP,LyonFESVD,LyonFast} on equispaced FEs.

\subsection{The exact equispaced Fourier extension}\label{ss:equiFEthy}
Consider first the case $M=N$.  Then $F_{N,N}(f)$ is equivalent to polynomial interpolation in $z$:

\prop{ Let $F_{N,N}(f)= f_N = f_{e,N}+f_{o,N} \in \cG_N$ be
defined by \R{equiL2} with $N=M$ and let $h_{i,N}(z)$
be given by \R{E:hNdef}.  Then $h_{i,N}(z)$, $i=1,2$ is the
$(N+1-i)^{\rth}$ degree polynomial interpolant of $h_i(z)$ at the
nodes $\{ z_n \}^{N}_{n=i-1} \subseteq [-1,1]$, where \be{
\label{unodes} z_n = m(x_{n}) = 2 \frac{\cos \left ( \tfrac{n
\pi}{NT} \right ) - c(T)}{1-c(T)} -1,\quad n=0,\ldots,N. }} 
This proposition allows us to analyze the theoretical
convergence/divergence of $F_{N,N}(f)$ using standard results on
polynomial interpolation.  Recall that associated with a set of
nodes $\{ z_n \}^{N}_{n=0}$ is a \textit{node density
function} $\mu(z)$, i.e.\ a function such that (i) $\int^{1}_{-1}
\mu(z) \D z = 1$ and (ii) each small interval $[z,z+h]$ contains a
total of $N \mu(z) h$ nodes for large $N$ \cite{Fps}.  In the case
of \R{unodes} we have

\lem{
The nodes \R{unodes} have node density function $\mu(z) = T /(\pi \sqrt{(1-z)(z-m(T))})$.
}

\prf{
Note first that $\int^{1}_{-1} \mu(z) \D z = 1$.  Now let $I = [z,z+h] \subseteq [-1,1]$ be an interval.  Then the node $z_n \in I$ if and only if $m^{-1}(z+h) \leq x_n \leq  m^{-1}(z)$.  Therefore, as $N \rightarrow \infty$, the proportion of nodes lying in $I$ tends to $m^{-1}(z)   - m^{-1}(z+h)$.  Now suppose that $h \rightarrow 0$.  Then
\eas{
m^{-1}(z+h)&=\frac{T}{\pi} \arccos \left [ c(T) + \frac{1-c(T)}{2} (z+h+1) \right ] = m^{-1}(z) - \mu(z) h + \ord{h^2}.
}
Thus $m^{-1}(z)   - m^{-1}(z+h)=  \mu(z) h + \ord{h^2}$, as required.
}

It is useful to consider the behaviour of $\mu(z)$.  When $z \rightarrow 1^{-}$,
$\mu(z) \sim T /(\pi \sqrt{1-z})$.  On the other hand, $\mu$ is
continuous at $z = -1$ with $\mu(-1) = \frac{T}{2 \pi} \tan \left ( \frac{\pi}{2T} \right )$.  Hence the nodes $\{ z_n
\}^{N}_{n=0}$ cluster quadratically near $z = 1$ and are linearly
distributed near $z=-1$.  It is well known that to avoid the Runge
phenomenon in a polynomial interpolation scheme, it is essentially
necessary for the nodes to cluster quadratically near both endpoints
(as is the case with Chebyshev nodes) \cite{Fps}.  If this is not
the case, one expects the Runge phenomenon: that is, divergence
(at a geometric rate) of the interpolant for any function having a singularity in a certain complex
region containing $[-1,1]$ (the \textit{Runge region} for the
interpolation scheme).  Since the nodes \R{unodes} do not exhibit
the correct clustering at the endpoint $z=-1$, we consequently expect this behaviour in the
equispaced FE $F_{N,N}(f)$.

As it transpires, the corresponding Runge region $\cR =
\cR(T)$ for $F_{N,N}$ can be defined in terms of the
potential function $\phi(t) = - \int^{1}_{-1} \mu(z) \log | t - z
| \D z + c$.  Here $c$ is an arbitrary constant.  Standard
polynomial interpolation theory \cite{Fps} then gives that \bes{
\cR(T)= \left \{ x \in \bbC : \phi(m(x)) = \phi(-1) \right
\}, } (observe that this is a subset of the complex $x$-plane). We
note also that the convergence/divergence of $F_{N,N}(f)$ at a
point $x$ will be exponential at a rate $\E^{\phi(m(x_0)) -
\phi(m(x))}$, where $x_0$ is the limiting singularity of $f$.  This
follows from a general result on polynomial interpolation
\cite{Fps}.  In particular, if $f$ has a singularity in
$\cR(T)$, then there will be some points $x \in [-1,1]$ for
which $F_{N,N}(f)$ diverges.

We next discuss two approaches for overcoming the Runge phenomenon in $F_{N,N}(f)$.

\subsubsection{Overcoming the Runge phenomenon I: varying $T$}\label{sss:oversamp}
One way to attempt to overcome (or, at least, mitigate) the Runge
phenomenon observed above is to vary the parameter $T$.  Note that: \lem{
The Runge region $\cR(T)$ satisfies $\cR(T)
\rightarrow [-1,1]$ as $T \rightarrow 1^+$, and $\cR(T)
\rightarrow \cR$  as $T \rightarrow \infty$, where $\cR$ is the
Runge region for equispaced polynomial interpolation.
}
\prf{
Suppose first that $T \rightarrow 1^+$. Since $m(T) \sim -1$, we have $\mu(z) \sim  \frac{1}{\pi \sqrt{1-z^2}}$.  The right-hand side is the potential function for Chebyshev interpolation, and thus the first result follows.

For the second result, we first recall that $\phi(m(x)) = - \int^{1}_{-1} \mu(z) \log | m(x) - z | \D z$.  Define the change of variable $z = m(s)$.  Since $m'(s) = - 1/\mu(m(s))$ we have
\bes{
\phi(m(x)) =- \int^{1}_{0} \log | m(x) - m(s) | \D s .
}
Note that
\bes{
m(x) - m(s) = \frac{\cos \frac{\pi x}{T} - \cos \frac{\pi s}{T}}{\sin^2 \frac{\pi}{2T}} =- \frac{2 \sin \tfrac{\pi(x-s)}{2T} \sin \tfrac{\pi(x+s)}{2T} }{\sin^2 \frac{\pi}{2T}} \sim - 2 (x-s)(x+s),\quad T \rightarrow \infty.
}
Therefore
\eas{
 \phi(m(x)) \sim  -  \int^{1}_{-1} \log \left | x - s\right |  \D s + C,\quad T \rightarrow \infty,
}
which is the potential function of equispaced polynomial interpolation, as required.
}

This lemma comes as no surprise.  As $T \rightarrow 1^+$ for fixed $N$, the system $\{ \E^{\I \frac{n \pi}{T} \cdot } \}_{|n| \leq N}$ tends to the standard Fourier basis on $[-1,1]$.  The problem of equispaced interpolation with trigonometric polynomials is well-conditioned and convergent.  On the other hand, when $T \rightarrow \infty$, the subspaces $\cC_{N}$ and $\cS_N$ both resemble spaces of algebraic polynomials in $x$.  Thus, in the large $T$ limit, $F_{N,N}(f)$ is an algebraic polynomial interpolant of $f$ at equispaced nodes.

Since the Runge region $\cR(T)$ can be made arbitrarily small by letting $T \rightarrow 1^+$, one way to overcome the Runge phenomenon is to vary $T$ in the way described in \S \ref{ss:Tchoice} and set $T = T(N;\epsilon)$.  One could also take $T \approx 1$ fixed.  However, this will always lead to a nontrivial Runge region, and consequently divergence of $F_{N,N}$ for some nonempty class of analytic functions.

\subsubsection{Overcoming the Runge phenomenon II: oversampling}
An alternative means to overcome the Runge phenomenon in $F_{N,M}(f)$ is to allow $M \geq N$.  Oversampling is known to defeat the Runge phenomenon in equispaced polynomial interpolation \cite{boyd2009divergence,BoydRunge,TrefPlatteIllCond}, and the same is true in this context (see \cite{BoydFourCont,brunoFEP} for previous discussions on oversampling for equispaced FEs).

It is now useful to introduce some notation.  For nodes $\{ x_n \}_{|n| \leq M}$ given by \R{equinodes}, let $(\cdot,\cdot)_M$ be the discrete bilinear form $(g,h)_M = \frac{1}{M+\frac12} \sum_{|n| \leq M} g(x_n) \overline{h(x_n)}$, and denote the corresponding discrete semi-norm by $\nm{\cdot}_{M}$.  Much as before, we define the condition number of $F_{N,M}$ by
\be{
\label{FNM_cond}
\kappa(F_{N,M} ) = \sup \left \{ \| F_{N,M}(b) \| : b \in \bbC^{2M+1},\ \nm{b} = 1 \right \}.
}
We now have:

\thm{
\label{t:oversamp}
Let $F_{N,M}(f)$ be given by \R{equiL2}, and suppose that
\be{
\label{DNM}
D(N,M) = \sup \left \{ \|\phi \| : \phi \in \cG_N,\  \| \phi \|_M = 1 \right \},
}
then
\bes{
\| f - F_{N,M}(f) \| \leq \sqrt{2}\left ( 1+D(N,M) \right ) \inf_{\phi \in \cG_N} \| f - \phi \|_{\infty}.
}
Moreover, the condition number $\kappa(F_{N,M}) = D(N,M)$.
}

\prf{
For the sake of brevity, we omit the first part of the proof (a very similar argument is given in \cite{boyd2009divergence} for the case of polynomial interpolation).  For the second part, we first notice that
\bes{
\kappa(F_{N,M}) = \sup \left \{ \| F_{N,M}(f) \| : f \in \rL^\infty(-1,1),\ \nm{f}_M = 1 \right \}.
}
Since $F_{N,M}(\phi) = \phi$ for $\phi \in \cG_N$ we have $\kappa(F_{N,M}) \geq D(N,M)$.  Conversely, since $F_{N,M}(f) \in \cG_N$, and since $F_{N,M}$ is an orthogonal projection with respect to the bilinear form $(\cdot,\cdot)_M$, we have $\| F_{N,M}(f) \| \leq D(N,M) \| F_{N,M}(f) \|_M \leq D(N,M) \| f \|_M$.  Hence $\kappa(F_{N,M}) \leq D(N,M)$, and we get the result.
}

This theorem implies that $F_{N,M}(f)$ will converge, regardless of the analyticity of $f$, provided $M$ is chosen such that $D(N,M)$ is bounded.  Note that this is always possible: $D(N,M) \rightarrow 1$ as $M \rightarrow \infty$ for fixed $N$ since $\nm{\cdot}_M$ is a Riemann sum approximation to $\nm{\cdot}$ and $\cG_N$ is finite-dimensional.  Up to small algebraic factors in $M$ and $N$, the quantity $D(N,M)$ is equivalent to
\be{
\label{tilDNM}
\tilde D(N,M) = \sup \left \{ \| p \|_{\infty} : p \in \bbP_N,\  | p(z_n) | \leq 1,\ n=0,\ldots,M \right \}.
}
Note the meaning of $\tilde D(N,M)$: it informs us how large a polynomial of degree $N$ can be on $[-1,1]$ if that polynomial is bounded at the $M$ points $z_n$.  Unfortunately, numerical evidence suggests that
\be{
\label{coppriv}
\alpha^{\frac{N^2}{M}} \leq \tilde D(N,M) \leq \beta^{\frac{N^2}{M}},
}
for constants $\beta \geq \alpha > 1$.  Thus one requires $M = \ord{N^2}$  nodes for boundedness of $D(N,M)$.  This is clearly less than ideal:  it means that we require many more samples of $f$ to compute its $N$-term equispaced FE.  In particular, the exact equispaced FE $F_{N,M}(f)$ of an analytic function $f$ will converge only root-exponentially fast in the number $M$ of equispaced grid values.

Had the nodes $\{ z_n \}^{M}_{n=0}$ clustered quadratically
near $z= \pm 1$, then $M = \ord{N}$ would be sufficient to ensure
boundedness of $\tilde D(N,M)$.  Note that when $N=M$, $\tilde
D(N,M)$ is precisely the Lebesgue constant of polynomial
interpolation.  On the other hand, if $\{ z_n \}^{M}_{n=0}$
were equispaced nodes on $[-1,1]$ then \R{coppriv} would coincide
with a well-known result of Coppersmith \& Rivlin
\cite{copprivlinpolygrowth}.  The intuition for a bound of the
form \R{coppriv} for the nodes \R{unodes} comes from the fact that
these nodes are linearly distributed near $z=-1$.  Thus, at least
near $z=-1$ they behave like equispaced nodes.

We remark that it is straightforward to show that the scaling $M = \ord{N^2}$ is sufficient for boundedness of $\tilde D(N,M)$.  This is based on Markov's inequality for polynomials.  Necessity of this condition would follow directly from the lower bound in \R{coppriv}, provided \R{coppriv} were shown to hold.  It may be possible to adapt the proof of \cite{copprivlinpolygrowth} to establish this result.

Since the scaling $M = \ord{N^2}$ is undesirable, one can ask what
happens when $M = \gamma N$ for some fixed oversampling parameter
$\gamma \geq 1$.  Using potential theory arguments, one can show
that $\tilde D(N, \gamma N)$  grows exponentially in $N$
(with the constant of this growth becoming smaller as $\gamma$
increases), as predicted by the conjectured bound \R{coppriv}. In other words,
\be{
\label{cPotentialDef}
N^{-1} \log D(N , \gamma N ) \sim \log c(\gamma;T),\quad N \rightarrow \infty,
}
for some $c(\gamma ; T) >1$\footnote{The constant
of growth was obtained in private communication with A.\ Kuijlaars.
A closed expression (up to several integrals involving the
potential function $\phi$ for the nodes $z_n$) can be found
for $c(\gamma;T)$. We omit the full argument as it is rather
lengthy, but note that it is based on standard results in
potential theory. A general reference is
\cite{ransford1995potentialtheory}.}.   In view of this behaviour, Theorem \ref{t:oversamp} guarantees
convergence of the FE \R{equiL2}, provided $\rho \geq c(\gamma ;
T)$, where $\rho$ is as in Theorem \ref{t:expconv}.  In other
words, $f$ needs to be analytic in the region
$\cD(c(\gamma;T))$ (recall $\cD$ from Theorem
\ref{t:expconv}) to ensure convergence.  Therefore, one expects
a Runge phenomenon whenever $f$ has a complex singularity lying in
the corresponding Runge region $\cR(\gamma ; T) =
\cD(c(\gamma;T))$. Naturally, a larger value of $\gamma$ leads to
a smaller (but still nontrivial) Runge region. However, regardless
of the choice of $\gamma$, there will always be analytic functions
for which one expects divergence of $F_{N,\gamma N}(f)$ (see
\cite{boyd2009divergence} for a related discussion in the case of
equispaced polynomial interpolation).  Moreover, the mapping $f \mapsto F_{N,\gamma N}$ will always be exponentially ill-conditioned for any fixed $\gamma$, since the condition number is precisely $D(N,\gamma N)$ (Theorem \ref{t:oversamp}).

Primarily for later use, we now note that it is also possible to study the condition number of the equispaced FE matrix $\bar{A}$ in a similar way.  Straightforward arguments show that
\be{
\label{minsingBNM}
1/\sigma_{\min}(\bar A)= B(N,M),\qquad B(N,M) =  \sup \left \{ \| \phi \|_{[-T,T]} : \phi \in \cG_N, \| \phi \|_M=1  \right \}.
}
Using the fact that $1/\sigma_{\min}(A) =  \sup \left \{ \| \phi \|_{[-T,T]} : \phi \in \cG_N, \| \phi \|=1  \right \}$, where $A$ is the matrix of the continuous FE, one can show that
\bes{
1/\sigma_{\min}(A) \lesssim B(N,M) \leq D(N,M) / \sigma_{\min}(A),
}
(here we use $\lesssim$ to mean up to possible algebraic factors in $N$).  Theorem \ref{t:exactCondNumb} now shows that $\bar{A}$ is always exponentially ill-conditioned in $N$, regardless of $M \geq N$.

Much like the case of $D(N,M)$ and $\tilde D(N,M)$, one can show that the quantity $B(N,M)$ is, up to algebraic factors, equivalent to
\be{
\label{tilBNM} \tilde
B(N,M) = \sup \left \{ \| p \|_{\infty,[m(T),1]} : p \in \bbP_N, | p(z_n) | \leq 1, n=0,\ldots,M \right \}.
}
Potential theory can be used  once more to determine the exact behaviour of $\tilde B(N,\gamma N)$.  In particular,
\be{
\label{dPotentialDef}
N^{-1} \log B(N,\gamma N) \sim d(\gamma;T),\quad N \rightarrow \infty,
}
for some constant $d(\gamma;T) \geq c(\gamma;T) > 1$.

\subsubsection{Numerical examples}\label{sss:numexampEqui}

In the previous section we established (up to the conjecture \R{coppriv}) 6., 7.\ and 8.\ of \S
\ref{s:introduction}. The main conclusion is as follows.  In
order to obtain a convergent FE in exact arithmetic using equispaced data one either needs to oversample quadratically (and thereby
reduce the convergence rate to only root-exponential), or scale
the extension parameter $T$ suitably with $N$ or both.  However,
recall from \S \ref{s:stability} that a FE obtained from a
finite precision computation may differ quite dramatically from
the corresponding infinite-precision extension.  Is it therefore
possible that the unpleasant effects described in the previous
section may not be witnessed in finite precision?  The answer
transpires to be yes, and consequently FEs can safely be used for
equispaced data, even in situations where divergence
is expected in exact arithmetic.

\begin{figure}
\begin{center}
$\begin{array}{ccc}
\includegraphics[width=6.25cm]{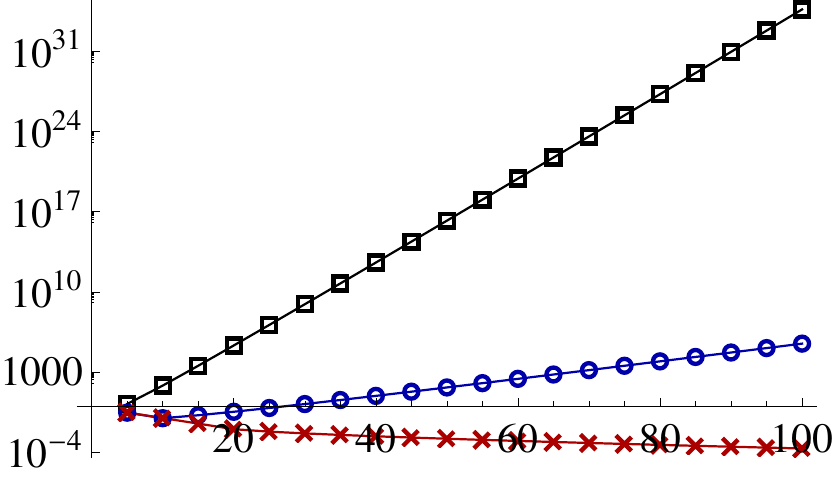} & \hspace{1pc}  & \includegraphics[width=6.25cm]{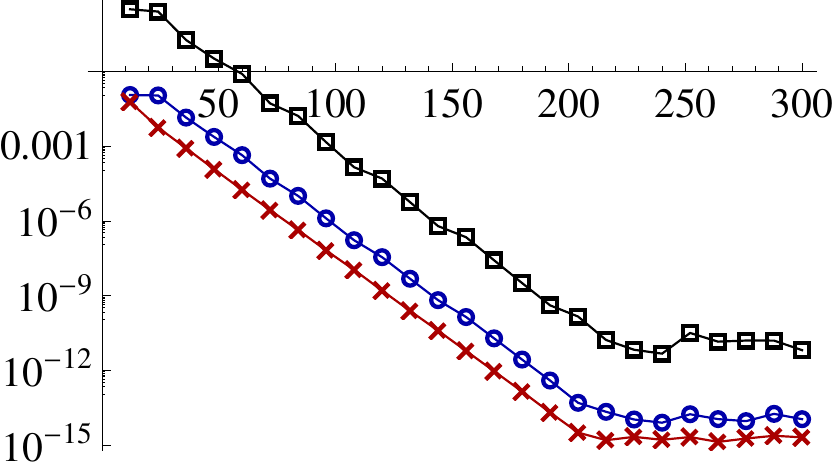}
\end{array}$
\caption{\small The error $\| f - f_N \|_{\infty}$ against $N$ for the equispaced FEs $f_N = F_{N,\gamma N}(f)$ (left) and $f_{N} = G_{N,\gamma N}(f)$ (right) of $f(x) = \frac{1}{1+100 x^2}$ with oversampling factor $\gamma = 1,2,4$ (squares, circles and crosses) and $T=2$.}  \label{f:EquiFnApp}
\end{center}
\end{figure}

To illustrate, consider the function $f(x) =
\frac{1}{1+100 x^2}$.  When $T=2$, this function has a singularity
lying in the Runge region $\cR(1;2)$.  The predicted divergence of its exact (i.e.\ infinite-precision) equispaced FE is shown in Figure \ref{f:EquiFnApp}.
Note that double oversampling also gives divergence, whilst with
quadruple oversampling the singularity of $f$ no longer lies in
$\cR(\gamma ; T)$.  We therefore witness geometric
convergence, albeit at a very slow rate. This behaviour is typical.  Given a function $f$ it is always possible to select the oversampling parameter $\gamma$ in such a way that $F_{N,\gamma N}(f)$ converges geometrically.  However, such a $\gamma$ depends on $f$ in a nontrivial manner (i.e.\ the location of the nearest complex singularity of $f$) and therefore cannot in practice be determined from the given data.  Note that this phenomenon---namely, the fact that careful tuning of a particular parameter in a function-dependent way allows geometric convergence to be restored---is also seen in other methods for approximating functions to high accuracy, such as the Gegenbauer reconstruction technique \cite{GottGibbsRev,GottGibbs1} (see Boyd \cite{BoydGegen} for a description of the phenomenon) and polynomial least squares \cite{boyd2009divergence}.

Fortunately, and unlike for these other methods, the situation changes completely for Fourier extensions when we carry out computations in finite precision.  This is shown in Figure \ref{f:EquiFnApp}.  For all choices of $\gamma$ used, the finite precision FE, which we denote $G_{N,\gamma N}(f)$, converges geometrically fast, and there is no drift in the error once the best achievable accuracy is attained.  Note that oversampling by a constant factor improves the approximation, but in all cases we still witness convergence.  In particular, no careful selection of $\gamma$, such as that discussed above, appears to be necessary in finite precision.

\subsection{The numerical equispaced Fourier extension}\label{ss:equianalysis}
We now explain these results by analysing the numerical equispaced FE.  Proceeding as in \S \ref{ss:numsolnanalysis} we shall consider the truncated SVD approximation, which we denote $H_{N,M,\epsilon}(f)$.  Note that a similar analysis has also recently been presented in \cite{LyonFESVD}; see Remark \ref{r:Lyon} for further details.

Let $\Phi_n \in \cG_N$ be the function corresponding to the right singular vector $v_n$ of the matrix $\bar{A}$.  Write $\cG_{N,M,\epsilon} = \spn \left \{ \Phi_n : \sigma_n > \epsilon \right \}$ and $\cG^{\perp}_{N,M,\epsilon} = \spn \{ \Phi_n : \sigma_n\leq\epsilon \}$, and note that $H_{N,M,\epsilon}$ is the orthogonal projection onto $\cG_{N,M,\epsilon}$ with respect to $(\cdot,\cdot)_{M}$.  Since $(\Phi_n,\Phi_m)_M = \sigma^2_n \delta_{n,m}$, we have
\be{
\label{equiTSVD}
H_{N,M,\epsilon}(f) = \sum_{n: \sigma_n > \epsilon} \frac{1}{\sigma^2_n} (f,\Phi_n)_M \Phi_n.
}
Our main result is as follows:

\thm{
\label{t:equiTSVD}
Let $f \in \rL^\infty(-1,1)$ and $H_{N,M,\epsilon}(f)$ be given by \R{equiTSVD}.  Then
\be{
\label{equiTSVDerr}
\| f - H_{N,M,\epsilon}(f) \| \leq  \sqrt{2} \left ( 1+ C_1(N,M;T,\epsilon) \right ) \| f - \phi \|_{\infty} + C_2(N,M;T,\epsilon) \| \phi \|_{[-T,T]},\quad \forall \phi \in \cG_N,
}
and
\be{
\label{equiTSVDcoeff}
\nm{a_{\epsilon}} = \| H_{N,M,\epsilon}(f) \|_{[-T,T]} \leq \frac{\sqrt{2}}{\epsilon} \| f - \phi \|_{\infty} + \| \phi \|_{[-T,T]},\quad \forall \phi \in \cG_N,
}
where
\ea{
\label{C1}
C_{1}(N,M;T,\epsilon) = \sup_{\substack{\phi \in \cG_{N,M,\epsilon} \\ \phi \neq 0}} \left \{ \frac{\| \phi \|}{\| \phi \|_M} \right \},
\quad
C_{2}(N,M;T,\epsilon) = \sup_{\substack{\phi \in \cG^{\perp}_{N,M,\epsilon} \\ \phi \neq 0}} \left \{ \frac{\| \phi \|}{\| \phi \|_{[-T,T]}}  \right \}.
}
}
\prf{
Let $\phi \in \cG_{N}$.  Then
\be{
\label{eqstep1}
\| f - H_{N,M,\epsilon}(f) \| \leq \| f - \phi \| + \| H_{N,M,\epsilon}(f-\phi)  \| + \| \phi - H_{N,M,\epsilon}(\phi) \|.
}
Consider the second term.  By definition of $C_1(N,M;T,\epsilon)$,
\bes{
\| H_{N,M,\epsilon}(f - \phi) \|  \leq C_{1}(N,M,\epsilon) \| H_{N,M,\epsilon }(f- \phi) \|_M \leq C_{1}(N,M,\epsilon) \| f - \phi \|_M,
}
where the second inequality follows from the fact that $H_{N,M,\epsilon}$ is an orthogonal projection with respect to $(\cdot,\cdot)_M$.  Noting that $\| g \|,\| g \|_M \leq \sqrt{2} \| g \|_{\infty}$ for any function $g \in \rL^\infty(-1,1)$ now gives the corresponding term in \R{equiTSVDerr}.  The bound for the third term of \R{eqstep1} follows immediately from the definition of $C_2(N,M;T,\epsilon)$ and the inequality $\| \phi - H_{N,M,\epsilon}(\phi) \|_{[-T,T]} \leq \| \phi \|_{[-T,T]}$.

For \R{equiTSVDcoeff}, we first write $\| H_{N,M,\epsilon}(f) \|_{[-T,T]} \leq \| H_{N,M,\epsilon}(f - \phi) \|_{[-T,T]} + \| H_{N,M,\epsilon}(\phi) \|_{[-T,T]}$.  Observe that for any $g \in \rL^\infty(-1,1)$ we have
\bes{
\| H_{N,M,\epsilon}(g) \|^2_{[-T,T]} = \sum_{n:\sigma_n > \epsilon} \frac{1}{\sigma^4_n} | (g,\Phi_n)_M |^2 \leq \frac{1}{\epsilon^2} \| H_{N,M,\epsilon}(g) \|^2_M \leq  \frac{1}{\epsilon^2} \| g \|^2_M \leq \frac{2}{\epsilon^2} \| g \|^2_{\infty}.
}
Also, $\| H_{N,M,\epsilon}(\phi) \|_{[-T,T]} \leq \| \phi \|_{[-T,T]}$ for $\phi \in \cG_N$.  Setting $g=f-\phi$ and combining these two bounds now gives \R{equiTSVDcoeff}.
}

\cor{
If $f \in \rL^\infty(-1,1)$ then $\| H_{N,M,\epsilon}(f) \| \leq \sqrt{2} / \epsilon \| f \|_{\infty}$, $\forall N \in \bbN$, $M \geq N$.  Moreover, if $f \in \rH^1(-1,1)$, $\bbT=[-T,T)$ is the $T$-torus and $c_1(T) > 0$ is as in Lemma \ref{l:smoothextension}, then
\bes{
\limsup_{\substack{N,M \rightarrow \infty \\ M \geq N}} \| H_{N,M,\epsilon}(f) \| \leq \inf \left \{ \| \tilde f \|_{[-T,T]} : \tilde{f} \in \rH^1(\bbT),\ \tilde{f} |_{[-1,1]} = f \right \} \leq c_1(T) \| f \|_{\rH^1(-1,1)}.
}
}
\prf{
By \R{equiTSVDcoeff}, we have
\be{
\label{normbd}
\| H_{N,M,\epsilon}(f) \| \leq \| H_{N,M,\epsilon}(f) \|_{[-T,T]} \leq \frac{\sqrt{2}}{\epsilon} \| f - \phi \|_{\infty} + \| \phi \|_{[-T,T]},\quad \forall \phi \in \cG_N.
}
Setting $\phi = 0$ gives the first result.  For the second, we let $\phi$ be the $N$-term Fourier series of $\tilde{f}$ on $\bbT$, so that $\| f - \phi \|_{\infty} \rightarrow 0$ as $N \rightarrow \infty$.  The final inequality follows from Lemma \ref{l:smoothextension}.
}
This corollary shows that the equispaced FE cannot suffer from a Runge phenomenon in finite precision, since it is bounded in $N$ and $M$.  This should come as no surprise.  Divergence of $H_{N,M,\epsilon}(f)$ would imply unboundedness of the coefficients $a_{\epsilon}$, a behaviour which is prohibited by truncating the singular values of $\bar{A}$ at level $\epsilon$.  Note that this corollary actually shows a much stronger result, namely that $H_{N,M,\epsilon}(f)$ is bounded on the extended domain $[-T,T]$, not just on $[-1,1]$.

Although this corollary demonstrates lack of divergence of $H_{N,M,\epsilon}(f)$, it says littles about its convergence besides the observation that $\|H_{N,M,\epsilon}(f)\|$ is asymptotically bounded by $\| f \|_{\rH^1(-1,1)}$.  To study convergence we shall use \R{equiTSVDerr}.  For this we first need to understand the constants $C_{i}(N,M;T,\epsilon)$.

\subsubsection{Behaviour of $C_{i}(N,M;T,\epsilon)$}\label{Cibehaviour}
\begin{figure}
\begin{center}
$\begin{array}{ccc}
  \includegraphics[width=4.75cm]{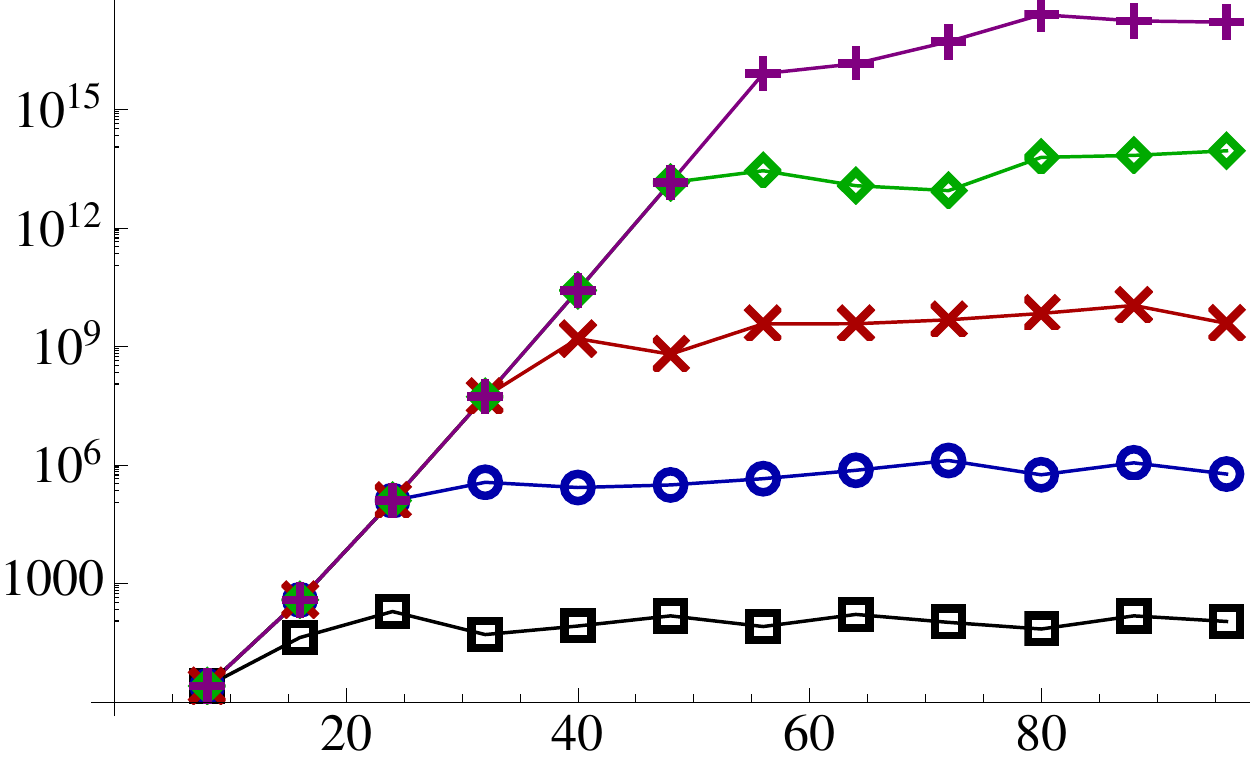} &  \includegraphics[width=4.75cm]{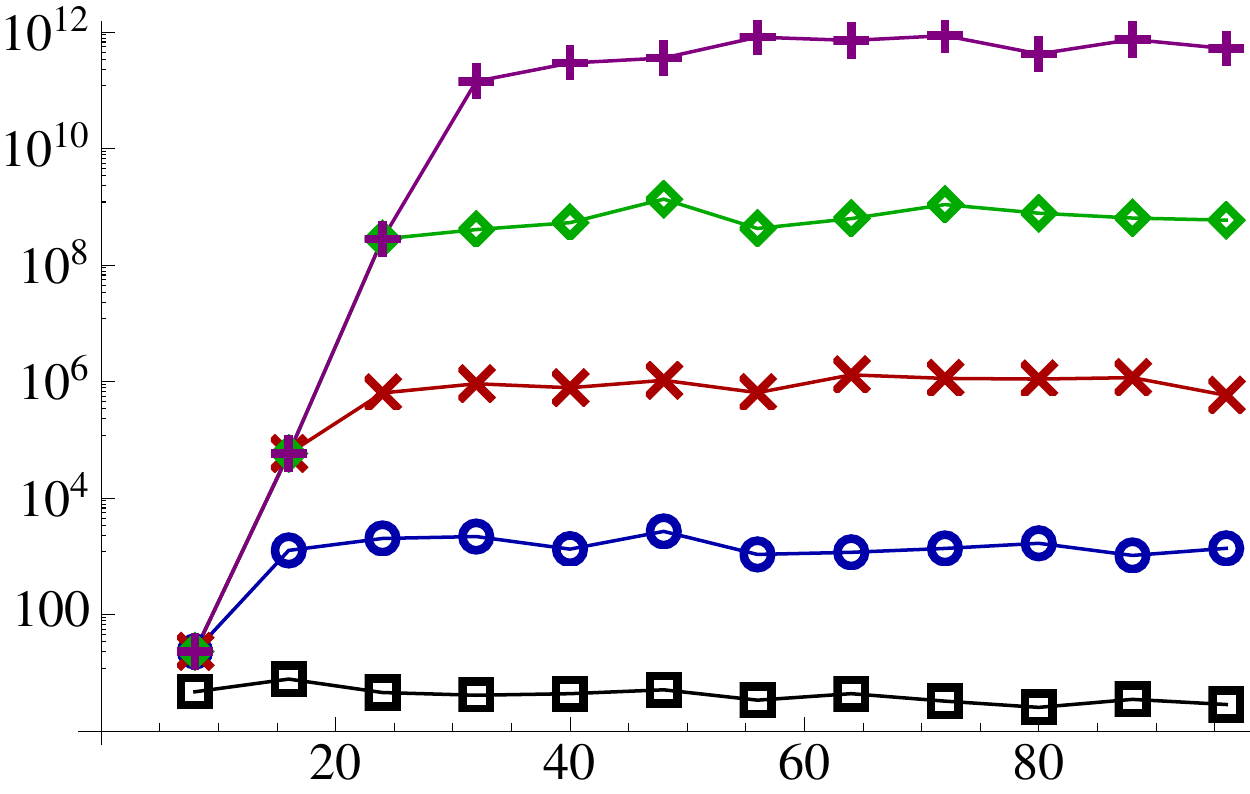}  & \includegraphics[width=4.7cm]{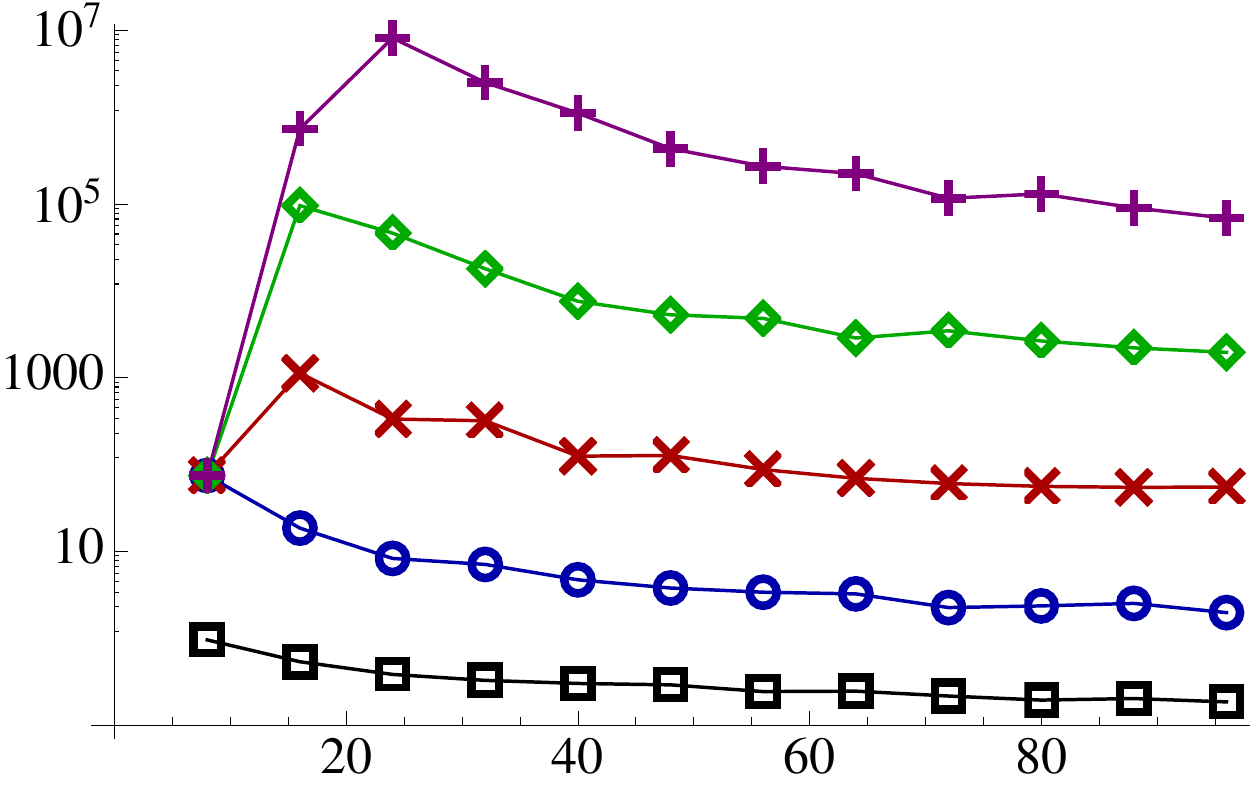}    \\    \includegraphics[width=4.75cm]{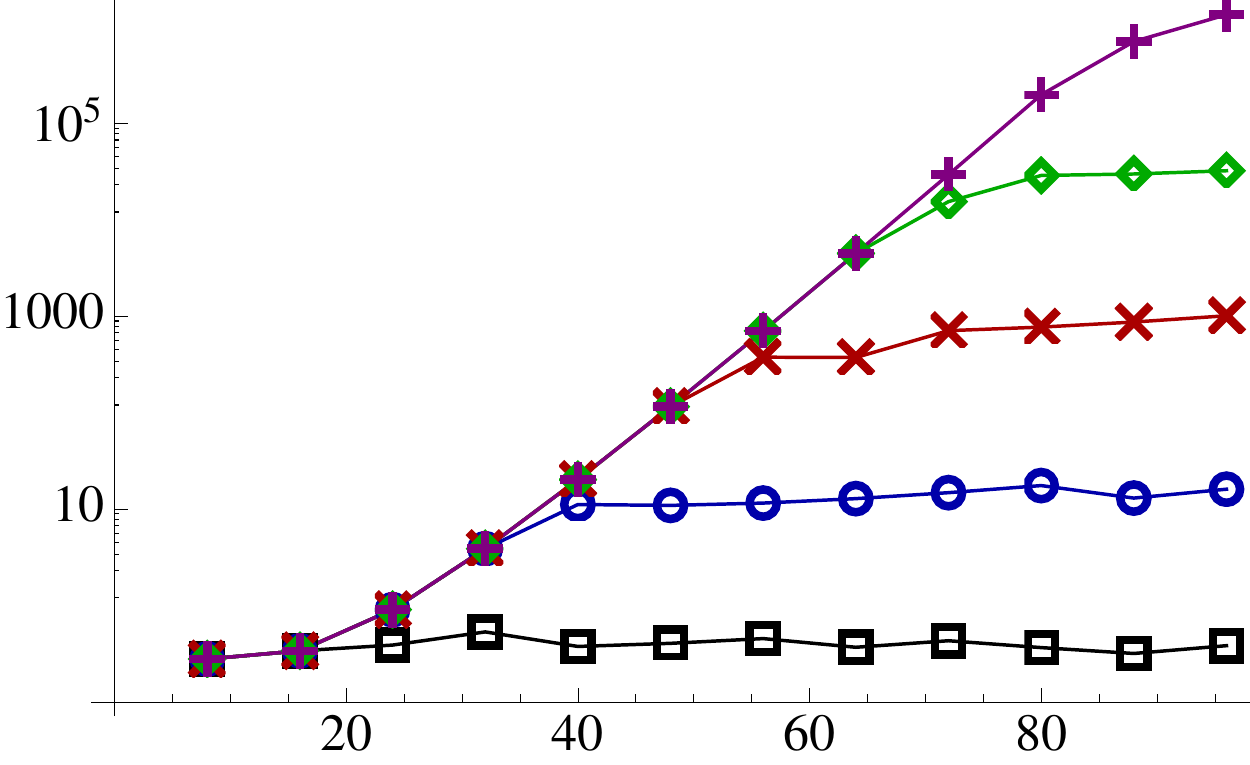} &  \includegraphics[width=4.75cm]{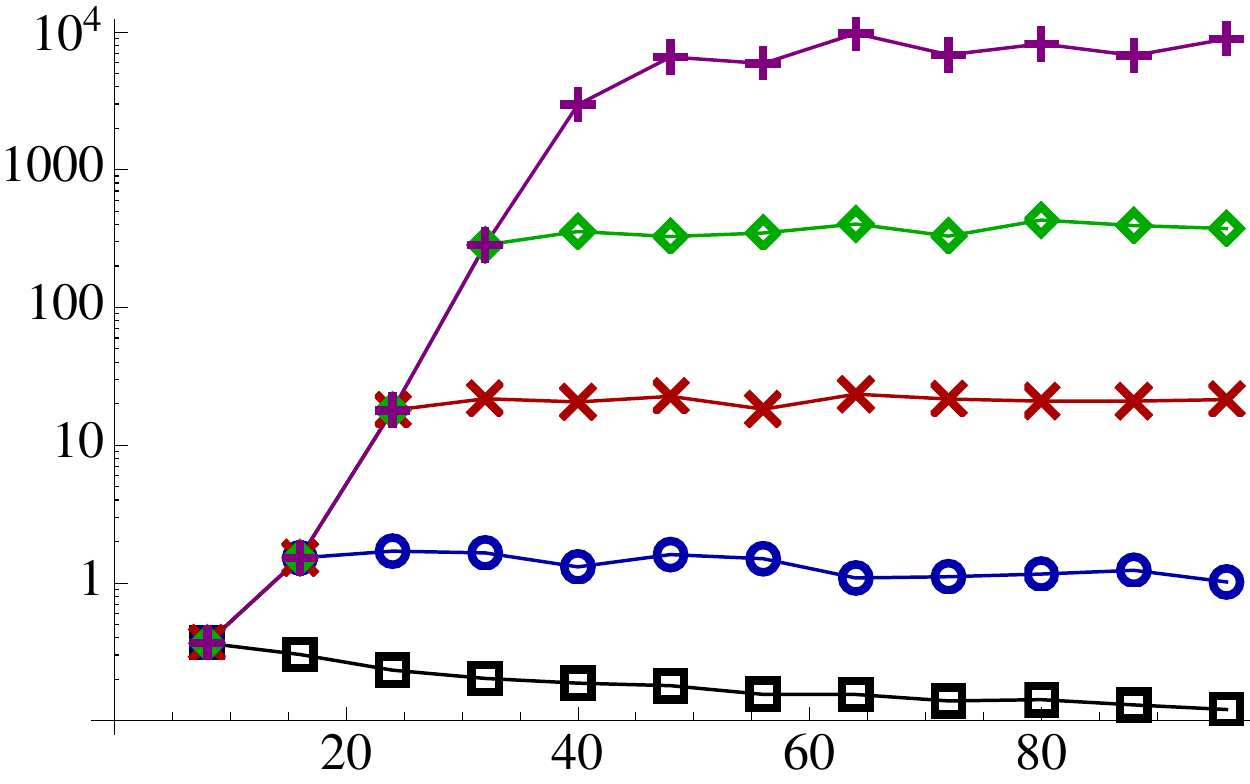}  & \includegraphics[width=4.7cm]{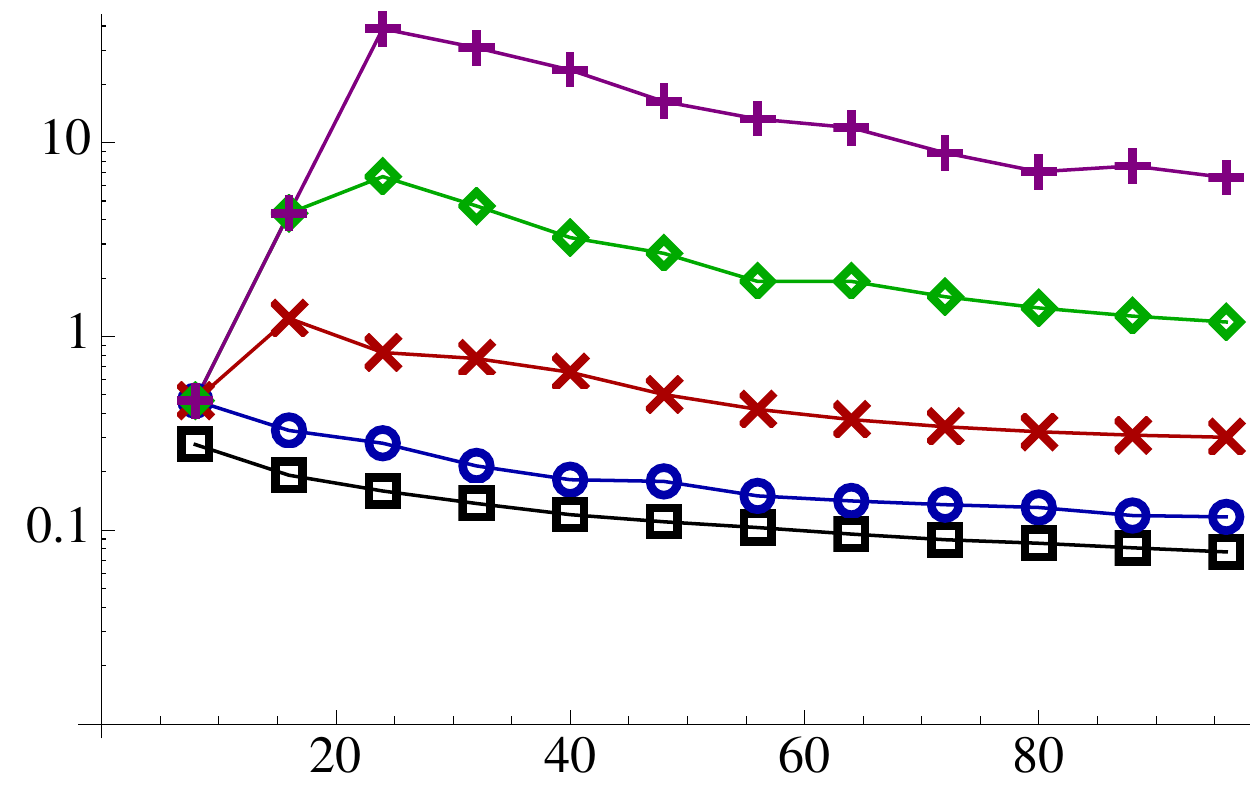}
   \\
 T=\frac43 & T=2 & T=4
\end{array}$
\caption{\small The quantity $C_1(N,\gamma N;T,\epsilon)$ against $N$ for $\gamma=1$ (top row) or $\gamma = 2$ (bottom row) and $\epsilon = 10^{-6},10^{-12},10^{-18},10^{-24},10^{-30}$ (squares, circles, crosses, diamonds and dashes respectively).}  \label{f:C1}
\end{center}
\end{figure}

Although Theorem \ref{t:equiTSVD} holds for arbitrary $M \geq N$, we now focus on the case of linear oversampling, i.e.\ $M = \gamma N$ for some $\gamma \geq 1$.

Let $N_2(\gamma,T,\epsilon)$ be the largest $N$ such that all the
singular values of $\bar{A}$ are at least
$\epsilon$ in magnitude:
\bes{
N_2(\gamma,T,\epsilon) = \max \left \{ N : \sigma_{\min}(\bar{A}) > \epsilon \right \}.
}
For $N \leq N_2(\gamma,T,\epsilon)$ we have
$\cG_{N, \gamma N,\epsilon} = \cG_N$ and therefore $C_1(N,\gamma N;T,\epsilon) =
D(N,\gamma N)$, where $D(N,M)$ is given by \R{DNM}.  Thus we witness exponential
divergence of $C_1(N,\gamma N;T,\epsilon)$ at rate $c(\gamma;T)$, where  $c(\gamma;T)$ is as in \R{cPotentialDef}.    This is shown numerically in Figure
\ref{f:C1}.

However, once $N > N_2(\gamma,T,\epsilon)$ the numerical
results in Figure \ref{f:C1} indicate a completely different behaviour: namely, $C_1(N,\gamma N;T,\epsilon)$ appears to be bounded.
Although we have no proof of this fact, these results strongly suggest the following conjecture:
 \be{ \label{conjecture1} C_1(N,\gamma N;T,\epsilon) \lesssim
C_1(N_2,\gamma N_2;T,\epsilon) \sim c(\gamma;T)^{N_2},\quad
\forall N \in \bbN. } In other words, $C_1(N,\gamma N;T,\epsilon)$ achieves its
maximal value in $N$ at $N \approx N_2$.  Recalling \R{minsingBNM} and \R{dPotentialDef}, we note that
\be{
\label{N2}
N_2(\gamma ,T,\epsilon) \approx -
\frac{\log \epsilon}{\log d(\gamma ; T)}.
}
Thus,
substituting this into bound \R{conjecture1} results in the following conjecture:
\be{ \label{C1bd}
C_1(N,\gamma N;T,\epsilon) \lesssim \min \left \{ c(\gamma;T)^N, \epsilon^{-\frac{\log
c(\gamma;T)}{\log d(\gamma;T)}} \right \},\quad \forall N \in \bbN.
}
In particular, $C_1(N,\gamma N;T,\epsilon)$ is bounded for all $N$ by some power
of $\epsilon^{-1}$.  Importantly, this power cannot be too large.
Note that  $c(\gamma;T) \leq d(\gamma;T)$, $\forall T > 1$, since the maximum of a  polynomial on $[m(T),1]$ is at least as large as its maximum on the smaller interval $[-1,1]$---compare \R{tilBNM} to \R{tilDNM}. Therefore the ratio $\frac{\log c(\gamma;T)}{\log
d(\gamma;T)}$ is at most one.  Moreover, by varying either $\gamma$ or $T$ we may decrease this ratio to arbitrarily close to $1$.  We discuss this further in the next section.

\begin{figure}
\begin{center}
$\begin{array}{ccc}
  \includegraphics[width=4.75cm]{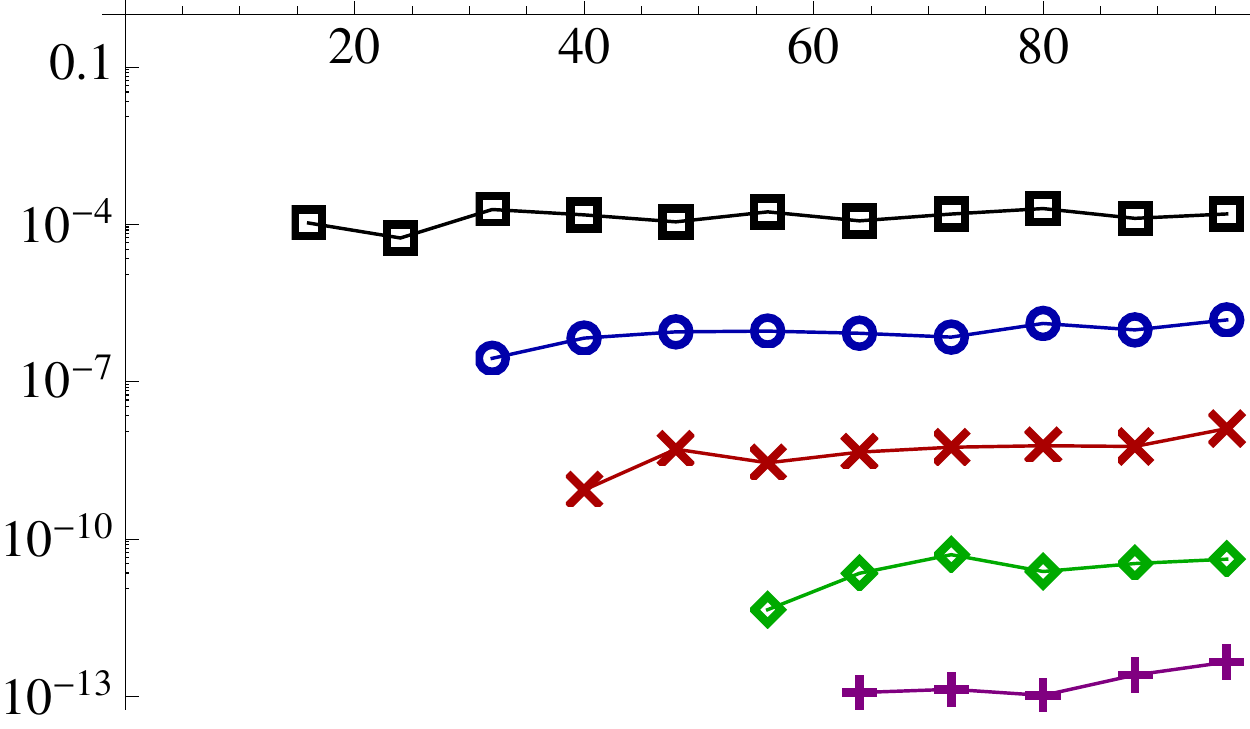} &  \includegraphics[width=4.75cm]{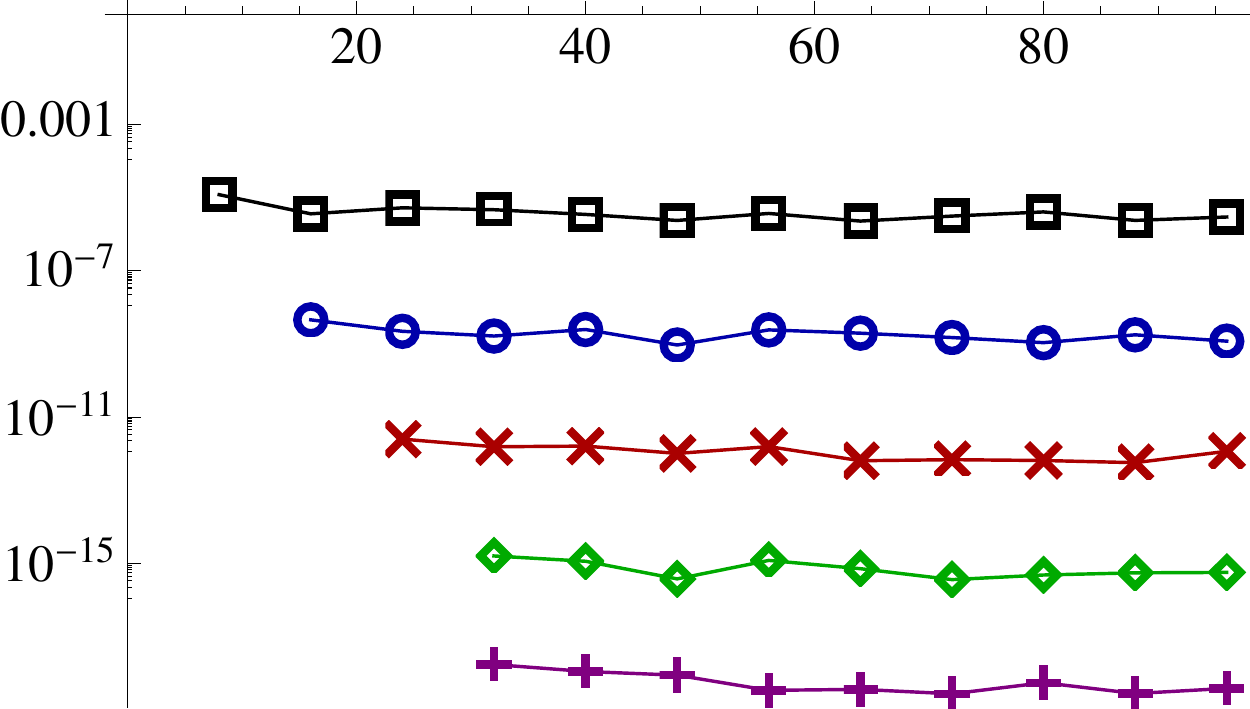}  & \includegraphics[width=4.7cm]{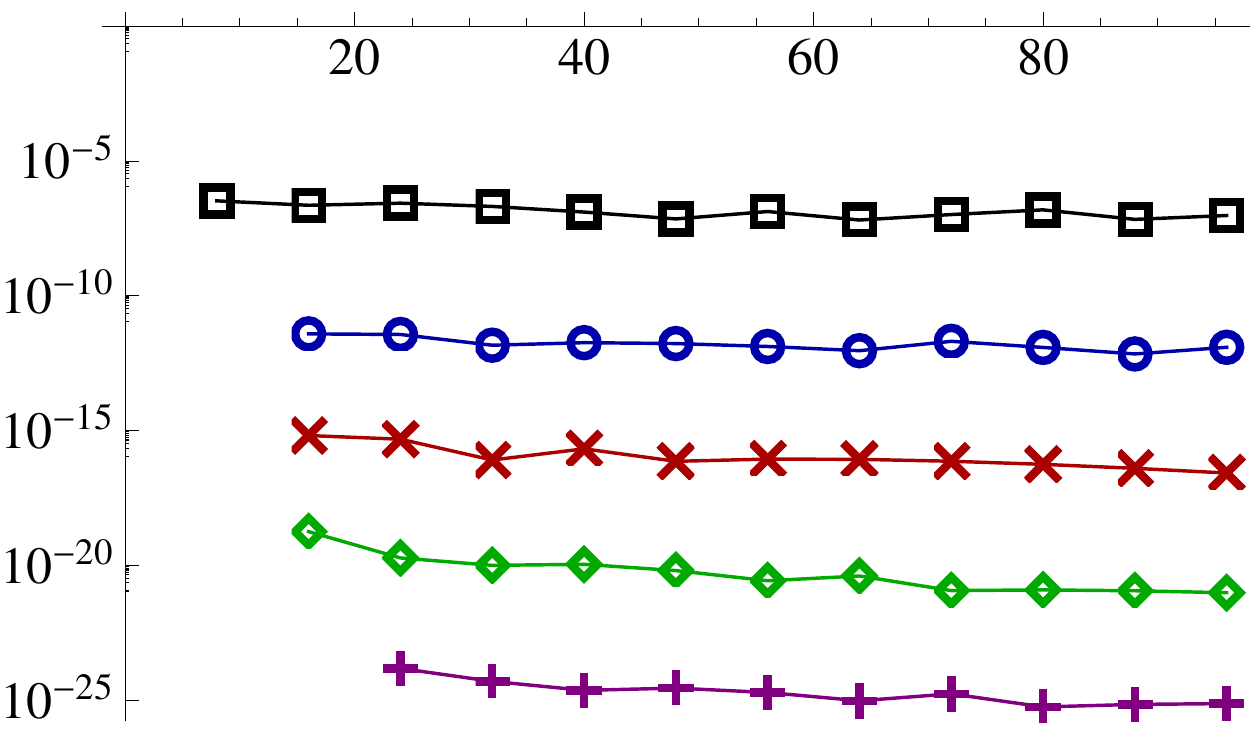}    \\    \includegraphics[width=4.75cm]{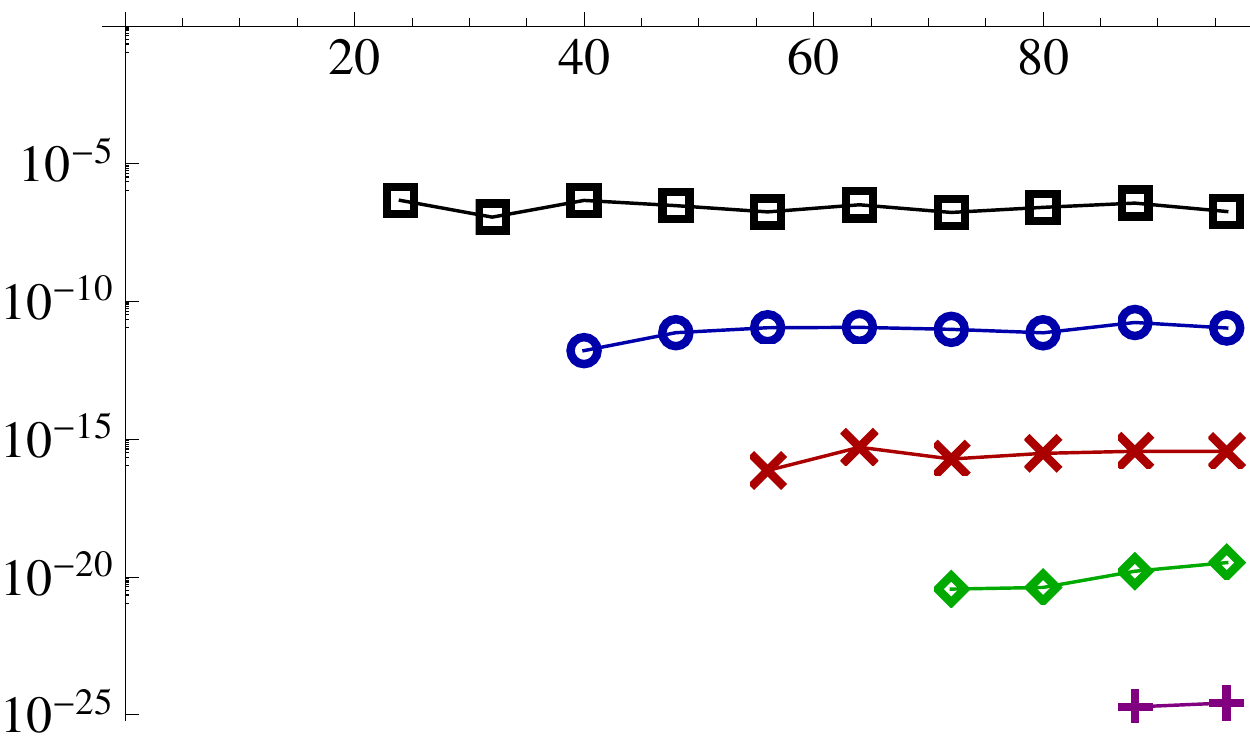} &  \includegraphics[width=4.75cm]{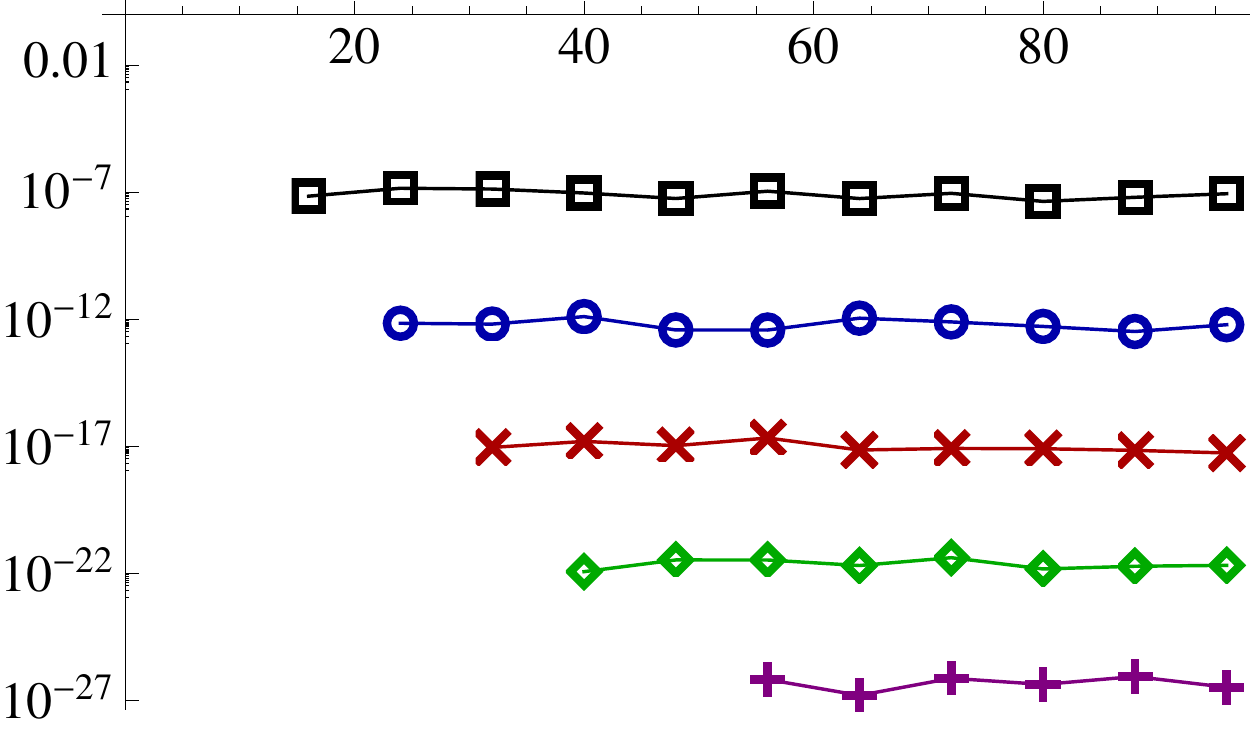}  & \includegraphics[width=4.7cm]{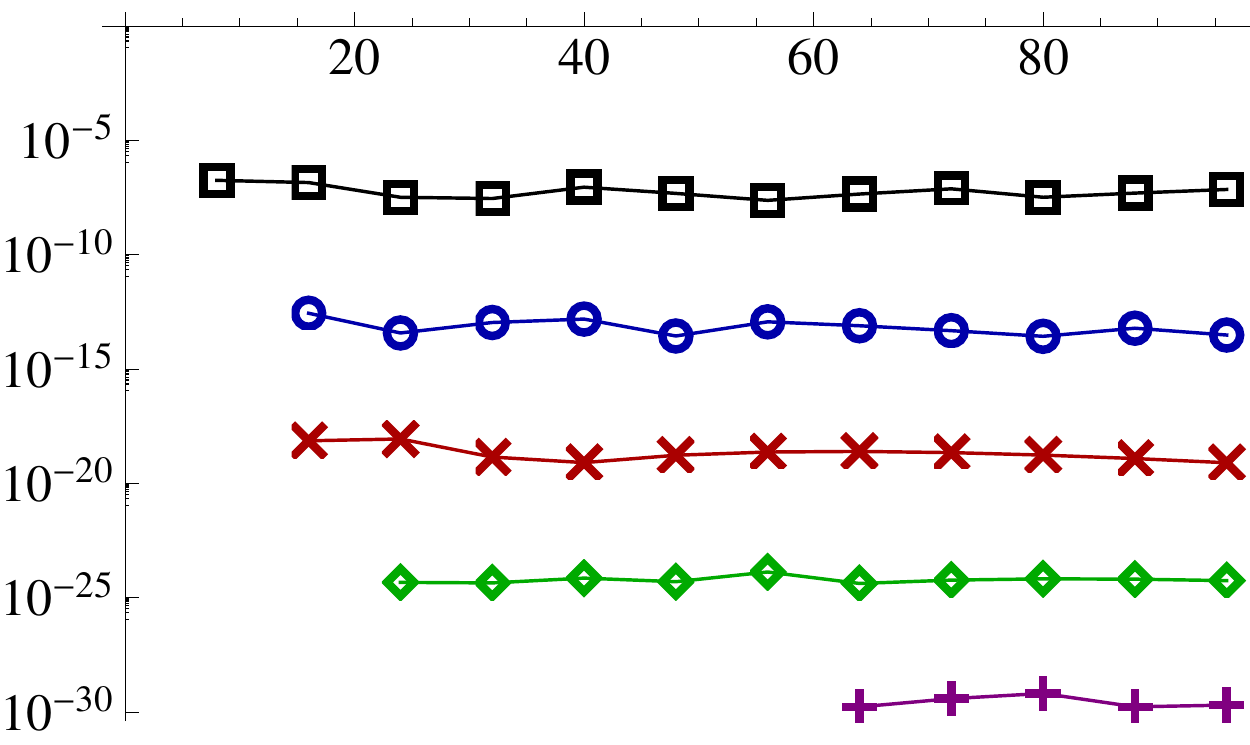}
   \\
 T=\frac43 & T=2 & T=4
\end{array}$
\caption{\small The quantity $C_2(N,\gamma N;T,\epsilon)$ against $N$ for $\gamma=1$ (top row) or $\gamma = 2$ (bottom row) and $\epsilon = 10^{-6},10^{-12},10^{-18},10^{-24},10^{-30}$ (squares, circles, crosses, diamonds and dashes respectively).}  \label{f:C2}
\end{center}
\end{figure}

The quantity $C_2(N,M;T,\epsilon)$ is harder to analyze, although clearly we have $C_2(N,M;T,\epsilon) = 0$ when $N < N_2$.  Figure \ref{f:C2} demonstrates that $C_2(N,\gamma N,\epsilon)$ is also bounded in $N$.  Moreover, closer comparison with Figure \ref{f:C1} suggests the existence of a bound of the form
\be{
\label{conjecture2}
C_2(N,\gamma N;T,\epsilon) \lesssim \epsilon C_1(N,\gamma N;T,\epsilon).
}
Once more, we have no proof of this observation.

\rem{
The quantities $C_1(N,M;T,\epsilon)$ and $C_2(N,M;T,\epsilon)$ have the explicit expressions
\bes{
C_1(N,M;T,\epsilon) = \sqrt{\| (S^{\epsilon})^{\dag} V^* A V (S^{\epsilon})^{\dag} \|},\quad C_2(N,M;T,\epsilon) = \sqrt{\| (V^{\epsilon})^* A V^{\epsilon} \|},
}
where $A$ is the continuous FE matrix, $U S V^*$ is the singular value decomposition of the equispaced FE matrix $\bar{A}$, $S^{\epsilon}$ is formed by replacing the $n^{\rth}$ column of $S$ by the zero vector whenever $\sigma_n \leq \epsilon$, and $V^{\epsilon}$ is formed by doing the same for columns of $V$ corresponding to indices $n$ with $\sigma_n > \epsilon$.  These expressions were used to obtain the numerical results in Figures \ref{f:C1} and \ref{f:C2}.  Computations were carried out with additional precision to avoid effects due to round-off.
}

\subsubsection{Behaviour of the truncated SVD Fourier extension}
Combining the analysis of the previous section with Theorem \ref{t:equiTSVD}, we now conjecture the bound
\be{
\label{exactnumerr}
\| f -H_{N,\gamma N,\epsilon}(f) \| \leq C(\gamma,T,\epsilon) \left ( \| f - \phi \|_{\infty} + \epsilon \| \phi \|_{[-T,T]} \right ),\quad \forall \phi \in \cG_N,
}
where $C(\gamma,T,\epsilon)$ is proportional to $\epsilon^{-a(\gamma;T)}$ and $a(\gamma;T)$ is given by
\be{
\label{agammaT}
a(\gamma;T) = \frac{\log c(\gamma;T)}{\log d(\gamma;T)}.
}
This estimate allows us to understand the behaviour of the numerical equispaced FE $G_{N,\gamma N}(f)$.  When $N < N_2$ we have $G_{N,\gamma N}(f) = F_{N,\gamma N}(f)$ and therefore $G_{N,\gamma N}(f)$ will diverge geometrically fast in $N$
whenever $f$ has a singularity in the Runge region
$\cR(\gamma;T)$ (see \S \ref{sss:oversamp}).  However, once
$N$ exceeds $N_2$, one obtains convergence.  Indeed,
setting $\phi = F_N(f)$ in \R{exactnumerr}, we find that the convergence is geometric up to the breakpoint $N_1$ (see \R{N1def}), and then, much as before, at least superalgebraic beyond that point.  Note that the maximal achievable accuracy of order $C(\gamma,T,\epsilon) \epsilon \approx \epsilon^{ 1-a(\gamma;T) }$.

In summary, we have now identified the following convergence behaviour for $H_{N,\gamma N,\epsilon}(f)$:
\begin{enumerate}
\item[(i)] $N < N_2(\gamma,T,\epsilon) \approx - \frac{\log \epsilon}{\log d(\gamma;T)}$.  Geometric divergence/convergence of $H_{N,\gamma N,\epsilon}(f)$ at a rate of, at worst, $c(\gamma;T) / \rho$, where $\rho$ is as in Theorem \ref{t:expconv} and $c(\gamma;T)$ is given by \R{cPotentialDef}.
\item[(ii)] $N_2(\gamma,T,\epsilon) \leq N < N_1(T,\epsilon) \approx - \frac{\log \epsilon}{\log E(T)}$.  Geometric convergence at a rate of at least $\rho$.
\item[(iii)]  $N= N_1(T,\epsilon)$.  The error
\bes{
\| f - H_{N,\gamma N,\epsilon}(f) \| \approx c_f \epsilon^{d_f-a(\gamma;T)},
}
where $a(\gamma;T)$ is as in \R{agammaT} and $d_f = \frac{\log \rho}{\log E(T)} \in (0,1]$.
\item[(iv)] $N \geq N_1(\gamma,T)$.  Superalgebraic convergence of $H_{N,\gamma N,\epsilon}(f)$ down to a maximal achievable accuracy proportional to $\epsilon^{1-a(\gamma;T)}$.
\end{enumerate}
(This establishes 10.\ of \S \ref{s:introduction}).  Much as in the case of the discrete FE, we see that if $f$ is analytic in $\cD(E(T))$, and if $c_f$ is not too large, then convergence stops at $N = N_1$ with maximal accuracy of order $c_f \epsilon^{1-a(\gamma;T)}$.  Otherwise, we have a further regime of at least superalgebraic convergence before this accuracy is reached.  

An important question is the role of the oversampling parameter $\gamma$ in this convergence.  We note:
\lem{
\label{l:a_behaviour}
Let $a(\gamma;T)$ be given by \R{agammaT}.  Then $a(\gamma;T)$ satisfies $0 \leq a(\gamma;T) \leq 1$ for all $\gamma$ and $T$.  Moreover, $a(\gamma;T) \rightarrow 0$ as $\gamma \rightarrow \infty$ for fixed $T$, and $a(\gamma;T) \rightarrow 0$ as $T \rightarrow \infty$ for fixed $\gamma$.
}
\prf{
Note that $c(\gamma;T) \leq d(\gamma;T)$.  Also $c(\gamma;T) \rightarrow 1$ and $d(\gamma;T) \rightarrow E(T)$ as $\gamma \rightarrow \infty$ for fixed $T$,  and $d(\gamma;T) \rightarrow \infty$ as $T \rightarrow \infty$ for fixed $\gamma$, whereas $c(\gamma;T)$ is bounded.  
}
This lemma suggests that increasing $\gamma$ will lead to a smaller constant $C(\gamma,T,\epsilon)$ in \R{exactnumerr}.  In fact, numerical results (Figures \ref{f:C1} and \ref{f:C2}) indicate that using $T=2$ and $\gamma =2$ gives a bound of a little over $1$ in magnitude for $\epsilon = 10^{-14}$.  Note that the effect of even just double oversampling is quite dramatic.  Without oversampling (i.e.\ $\gamma = 1$), the constant $C(\gamma,T,\epsilon)$ is approximately $10^4$ in magnitude when $\epsilon = 10^{-14}$ (see Figures \ref{f:C1} and \ref{f:C2}).

\begin{figure}
\begin{center}
$\begin{array}{ccc}
\includegraphics[width=6.25cm]{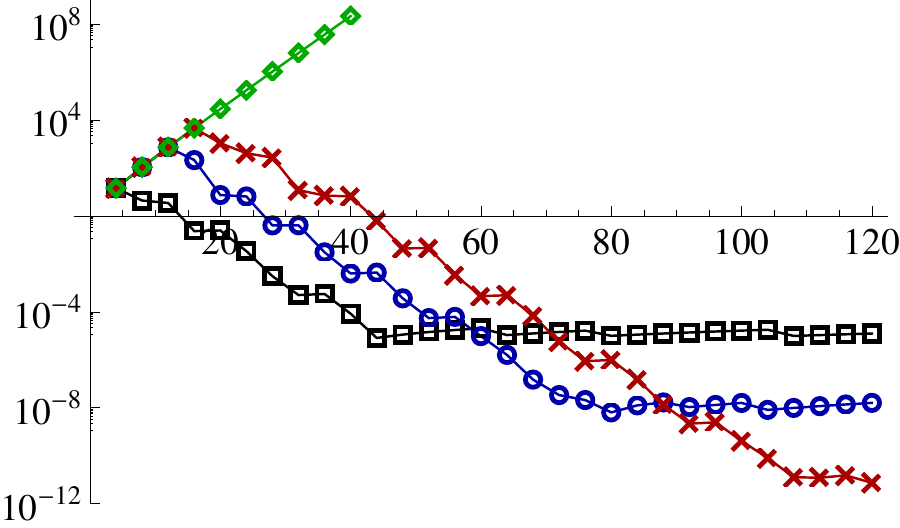} & \hspace{1pc} & \includegraphics[width=6.25cm]{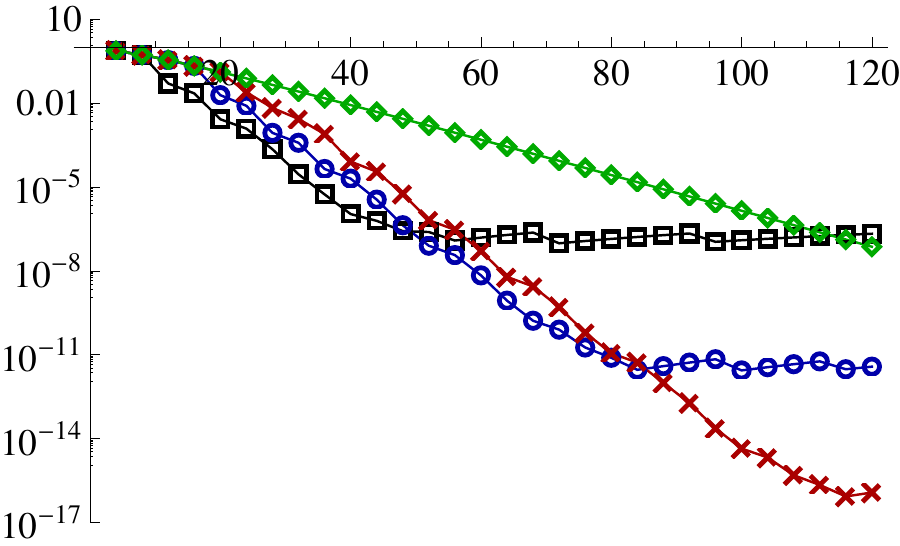}
\end{array}$
\caption{\small Error against $N$ for $H_{N,\gamma N , \epsilon}(f)$, where $f(x) = \frac{1}{1+16 x^2}$, $T=2$, $\gamma = 1$ (left) or $\gamma = 2$ (right) and $\epsilon = 10^{-6} , 10^{-12}, 10^{-18}$ (squares, circles, crosses).  Diamonds correspond to the exact equispaced FE $F_{N,\gamma N}(f)$.}  \label{f:ManSVDequi}
\end{center}
\end{figure}

Let us make several further remarks.  First, in practice the regime $N < N_1$ is typically very small---recall that $N_1$ is around $20$ for $T=2$ (see \S \ref{sss:discSVDanalysis})---and therefore one usually does not witness all three types of behaviour in numerical examples.  Second, as $\gamma \rightarrow \infty$, we have $N_2 \rightarrow N_1$ (recall that $d(\gamma;T) \rightarrow E(T)$ as $\gamma \rightarrow \infty$).  Thus, with a sufficient amount of oversampling, the regime (ii) will be arbitrarily small.  On the other hand, oversampling decreases $c(\gamma;T)$, and therefore the rate of divergence in the regime (i) is also lessened by taking $\gamma > 1$.  Indeed, the numerical results in Figure \ref{f:ManSVDequi}, as well as in \S \ref{ss:NumExp} later, indicate that oversampling by a factor of $2$ is typically sufficient in practice to mitigate the effects of divergence for most reasonable functions.

Figure \ref{f:ManSVDequi} confirms these observations for the function $f(x) = \frac{1}{1+16 x^2}$.  For $\gamma = 1$ the initial exponential divergence is quite noticeable.  However, this effect largely vanishes when $\gamma =2$.  Notice that a larger cutoff $\epsilon$ actually gives a smaller error initially, since there is a smaller regime of divergence.  However, the maximal achievable accuracy is correspondingly lessened.  We note also that maximal achievable accuracies for $\epsilon = 10^{-6},10^{-12},10^{-18}$ are roughly $10^{-4}$, $10^{-8}$ and $10^{-12}$ respectively when $\gamma = 1$ and $10^{-7}$, $10^{-12}$ and $10^{-16}$ when $\gamma = 2$.  These are in close agreement with the corresponding numerical values of $C_2(N,\gamma N; T , \epsilon)$ (see Figure \ref{f:C2}), as predicted by Theorem \ref{t:equiTSVD}.

\rem{\label{r:scaling}
A central conclusion of this section is that one requires a lower asymptotic scaling of $M$ with $N$ for the numerical equispaced FE than for its exact counterpart.  Since $\cG_{N,M,\epsilon}$ is a subset of $\cG_N$, we clearly have $C_{1}(N,M;T,\epsilon) \leq D(N,M)$,  where $D(N,M)$ is given by \R{DNM}.  Hence quadratic scaling $M = \ord{N^2}$ is sufficient (see \S \ref{sss:oversamp}) to ensure boundedness of $C_{1}(N,M;T,\epsilon)$, and one can make a similar argument for $C_{2}(N,M;T,\epsilon)$.  However, Figures \ref{f:C1} and \ref{f:C2} indicate that this condition is not necessary, and that one can get away with the much reduced scaling $M=\ord{N}$ in practice.

This difference can be understood in terms of the singular values of $\bar{A}$.  Recall that small singular values correspond to functions $\phi \in \cG_N$ with $\| \phi \|_{[-T,T]} \gg \| \phi \|_M$.  Now consider an arbitrary $\phi \in \cG_N$.  If the ratio $\| \phi \| / \| \phi \|_M$ is large, it suggests that $\phi$ lies approximately in the space $\cG^{\perp}_{N,M,\epsilon}$ corresponding to small singular values.  Hence, $\| \phi \| / \| \phi \|_M$ cannot be too large over $\phi \in \cG_{N,M,\epsilon}$, and thus we see boundedness of $C_{1}(N,M,\epsilon)$, even when $D(N,M)$---the supremum of this ratio over the whole of $\cG_N$---is unbounded.
}

\rem{
\label{r:Lyon}
A similar analysis of the equispaced FE, also based on truncated SVDs, was recently presented by M. Lyon in \cite{LyonFESVD}.  In particular, our expressions \R{equiTSVDerr} and \R{exactnumerr} are similar to equations (30) and (31) of \cite{LyonFESVD}.  Lyon also provides extensive numerical results for his analogues of the quantities $C_1(N,M;T,\epsilon)$ and $C_2(N,M;T,\epsilon)$, and describes a bound which is somewhat easier to use in computations.  The main contributions of our analysis are the conjectured scaling of the constant $C(\gamma,T,\epsilon)$ in terms of $\epsilon$, $\gamma$ and $T$, the description and analysis of the breakpoints $N_2$ and $N_1$, and the differing convergence/divergence in the corresponding regions.
}

\subsubsection{The condition number of the numerical equispaced FE}\label{ss:equistabfn}

\begin{table}
\begin{center}
\begin{tabular}{|c||c|c|c|c|c|}
\hline  N & 40 & 80 & 120 & 160 & 200
\\ \hline
$\gamma = 1$ &  $2.37 \times 10^4$    &  $3.50 \times 10^4$  & $2.24 \times 10^4$ &  $2.47 \times 10^4$ &  $1.93 \times 10^4$
\\ \hline
$\gamma = 2$ &  $2.18 \times 10^1$    &  $2.66 \times 10^1$  & $2.40  \times 10^1$ &  $2.56  \times 10^1$ &  $2.47  \times 10^1$
\\ \hline
$\gamma = 4$ &  $8.03 \times 10^0$    &  $1.05  \times 10^1$  & $1.23  \times 10^1$  &  $1.39  \times 10^1$ &  $1.54  \times 10^1$
\\ \hline
\end{tabular}
\caption{ The function $K(G_{N,\gamma N})$ against $N$ with $T=2$ and $\gamma = 1, 2, 4$.}\label{t:FEStabFnEquispaced}
\end{center}
\end{table}

We now consider the condition number $\kappa(G_{N,M})$ (defined as in \R{FNM_cond}) of the numerical equispaced extension.  In Table \ref{t:FEStabFnEquispaced} we plot $K(G_{N,\gamma N})$ against $N$, where $K(G_{N,M})$ is an upper bound for $\kappa(G_{N,M})$ defined analogously to \R{KGN_def}.  The results indicate numerical stability, and, as we expect, improved stability with more oversampling.

\begin{figure}
\begin{center}
$\begin{array}{c}
\includegraphics[width=6.25cm]{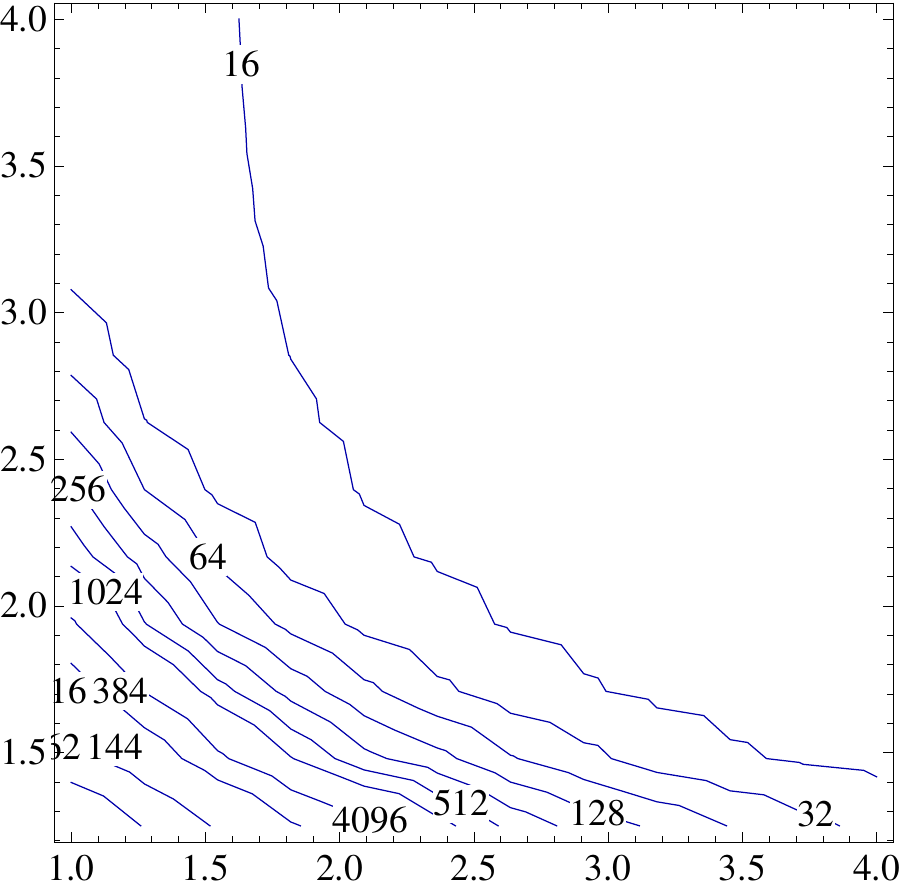}
\end{array}$
\caption{\small Contour plot of the quantity $K(G_{N,\gamma N})$ against $1\leq \gamma \leq 4$ and $1 < T \leq 4$ for $N=200$.} \label{f:contours}
\end{center}
\end{figure}

Besides oversampling it is also possible to improve stability by varying the extension parameter $T$.  In Figure \ref{f:contours} we give a contour plot of $K(G_{N,\gamma N})$ in the $(\gamma,T)$-plane.  Evidently, increasing $T$ improves stability.  Recall, however, that a larger $T$ corresponds to worse resolution power (see \S \ref{ss:Tchoice}).  Conversely, increasing $\gamma$ also leads to worse resolution when measured in terms of the total number $M = \gamma N$ of equispaced function values required.  Hence a balance must be struck between the two quantities.  Figure \ref{f:contours} suggests that $\gamma = T = 2$ is a reasonable choice in practice.  Recall also that the choice $T=2$ allows for fast computation of the equispaced FE (Remark \ref{sss:T2}), and hence is desirable to use in computations.

The behaviour of the condition number can be investigated with the following theorem (the proof is similar to that of Theorem \ref{t:stabfn} and hence omitted):

\thm{
\label{t:equistabfn}
The condition number $\kappa(H_{N,M,\epsilon})$ of the truncated SVD equispaced FE $H_{N,M,\epsilon}$ satisfies $\kappa(H_{N,M,\epsilon}) =C_1(N,M;T,\epsilon)$, where $C_{1}(N,M;T,\epsilon)$ is given by \R{C1}.
}

From the analysis of \S \ref{Cibehaviour} we conclude that $\kappa(H_{N,\gamma N,\epsilon}) \lesssim \epsilon^{-a(\gamma;T)}$, where $a(\gamma;T)$ is as in \R{agammaT}.  Lemma \ref{l:a_behaviour} therefore shows that $\kappa(H_{N,\gamma N,\epsilon}) \lesssim 1$ as $\gamma \rightarrow \infty$ for fixed $T$, and $\kappa(H_{N,\gamma N,\epsilon}) \lesssim 1$ as $T \rightarrow \infty$ for fixed $\gamma$.  This confirms the behaviour described above.

\subsubsection{Numerical examples}\label{ss:NumExp}
\begin{figure}
\begin{center}
$\begin{array}{ccc}
  \includegraphics[width=6.25cm]{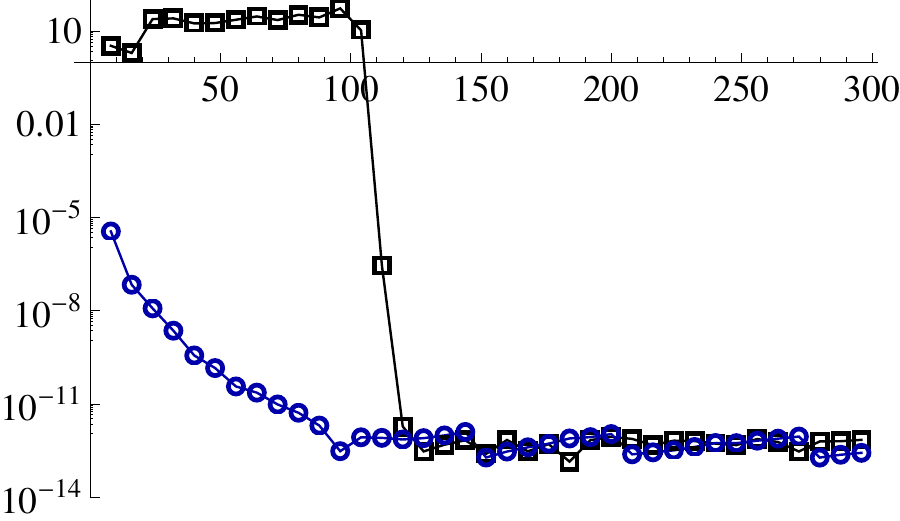} & \hspace{1pc} & \includegraphics[width=6.25cm]{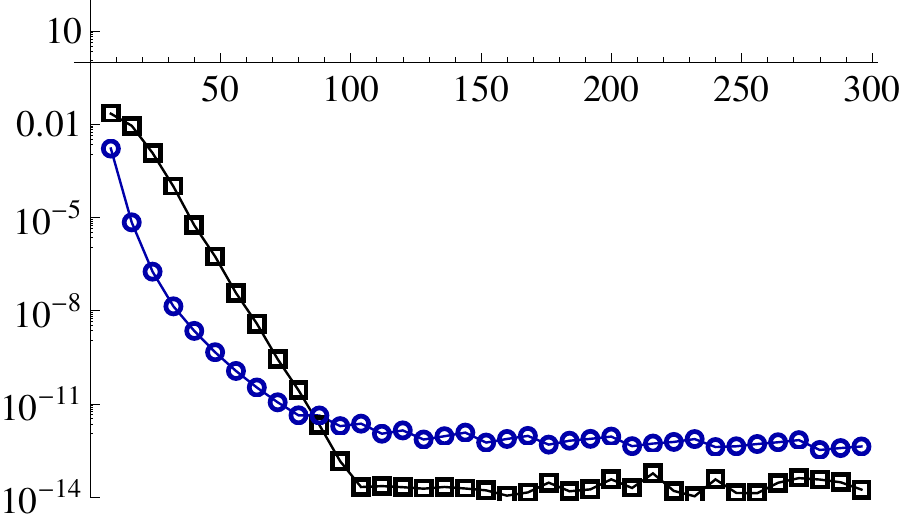} \end{array}$
\caption{\small The error $\| f - G_{N,\gamma N}(f) \|_{\infty}$,
where $\gamma = 2$ and $T=2$.  Left: $f(x) = \E^{25 \sqrt{5} \pi \I x}$ (squares), $f(x) = | x |^7$ (circles).  Right: $f(x) = \frac{1}{1+25 x^2}$ (squares), $f(x) = \frac{1}{8-7 x}$ (circles).  }
\label{f:EquiFnApp2}
\end{center}
\end{figure}

\begin{figure}
\begin{center}
$\begin{array}{ccc}
\includegraphics[width=6.25cm]{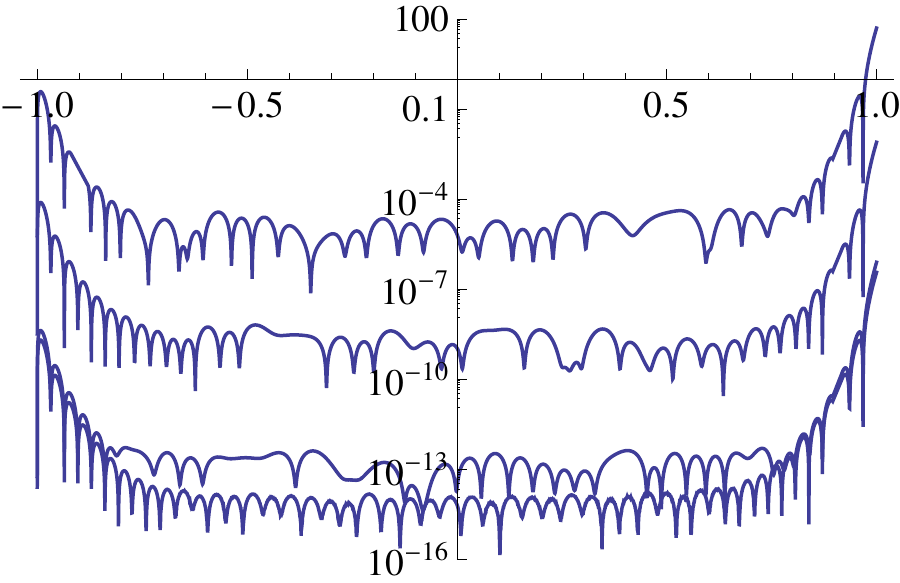} & \hspace{1pc} & \includegraphics[width=6.25cm]{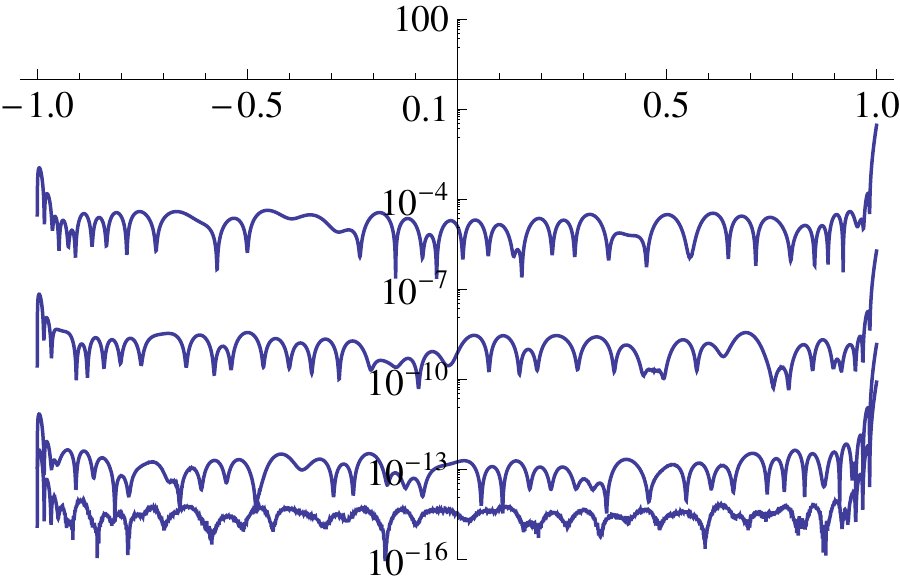}
\end{array}$
\caption{\small The error $| f(x) - G_{N,\gamma N}(f)(x) |$ against $x$, where  $\gamma =1$ (left) or $\gamma = 2$ (right), for $N = 30$, $T=2$ and $f(x) = \E^x$, with noise at amplitudes $\delta = 10^{-4},10^{-6},10^{-8},10^{-10},0$.}  \label{f:NoisyEqui}
\end{center}
\end{figure}

In Figure \ref{f:EquiFnApp2} we consider the equispaced FE for four test functions.  In all cases we use $\gamma =2$ and $T=2$.  As is evident, all choices of $T$ give good, stable numerical results, with the best achievable accuracy being at least $10^{-12}$.  Robustness in the presence of noise is shown in Figure \ref{f:NoisyEqui}.  Observe that when $\gamma = 1$, noise of amplitude $\delta$ is magnified by around $10^{5}$, in a manner consistent with Theorem \ref{t:equistabfn}.  Conversely, with double oversampling, this factor drops to less than $10^{2}$, again in agreement with Theorem \ref{t:equistabfn}.

\subsection{Relation to the theorem of Platte, Trefethen \& Kuijlaars}\label{ss:PTKrelation}
We are now in a position to explain how FEs relate to the impossibility theorem of Platte, Trefethen \& Kuijlaars \cite{TrefPlatteIllCond}.  A restatement of this theorem (with minor modifications to notation) is as follows:
\thm{[\cite{TrefPlatteIllCond}]
\label{t:impossibility}
Let $F_M$, $M \in \bbN$, be a sequence of approximations such that $F_M(f)$ depends only on the values of $f$ on an equispaced grid of $M$ points.  Let $\cE \subseteq \bbC$ be compact and suppose that there exists $C < \infty$, $\alpha>1$ and $\tau \in (\frac12,1]$ such that
\be{
\label{FMerror}
\| f - F_M(f) \|_{\infty} \leq C c_f \alpha^{-M^{\tau}},\qquad c_f = \max_{x \in \cE} | f(x) |,
}
for all $M \in \bbN$ and all $f$ that are continuous on $\cE$ and analytic in its interior.  Then there exists a $\beta > 1$ such that the condition numbers $\kappa(F_M) \geq \beta^{M^{2 \tau -1}}$ for all sufficiently large $M$.
}
This theorem has two important consequences. First, any exponentially convergent method is also exponentially ill-conditioned.  Second,  the best possible convergence for a stable method is root-exponential in $M$.  Note that the theorem is valid for all methods, both linear and nonlinear, that satisfy \R{FMerror}.

Consider now the exact equispaced Fourier extension $F_{N,M}$.  As shown in \S \ref{ss:equiFEthy}, when $N = \ordu{\sqrt{M}}$ this method is stable and root-exponentially convergent.  Hence equispaced FEs in infinite precision \textit{attain the maximal possible convergence rate for stable methods} satisfying the conditions of the theorem.

Now consider the numerical equispaced FE $G_{\eta M,M}$, where $0 < \eta \leq 1$ is the reciprocal of the oversampling parameter $\gamma$ used in the previous sections.  We have shown that this approximation is stable, so at least one condition in Theorem \ref{t:impossibility} must be violated.  Suppose that we take $\cE = \cD(E(T))$, for example.  Then \R{exactnumerr} shows that
\be{
\label{conv_nice}
\| f - G_{\eta M,M}(f) \|_{\infty} \lesssim   c_f \epsilon^{-a(\eta^{-1} ; T)}\left ( \left ( E(T)^{\eta} \right )^{-M} + \epsilon \right ).
}
The finite term $\epsilon$ in the brackets means that this approximation does not satisfy \R{FMerror}, and hence Theorem \ref{t:impossibility} does not apply.  Recall that if $c_f$ is small then \R{conv_nice} describes the full convergence behaviour for all $M$.  On the other hand, if $c_f$ is large, or if $f \in \cD(\rho)$ with $\rho < E(T)$, then the convergence is, after initial geometric convergence, at least superalgebraic down to the maximal achievable accuracy $\epsilon^{1-a(\eta^{-1};T)}$.  This is also not in contradiction with the conditions of Theorem \ref{t:impossibility}.

To summarize, equispaced FEs, when implemented in finite precision, possess both numerical stability and rapid convergence, and hence allow one to circumvent the impossibility theorem to an extent.  In particular, for all functions $f \in \cD(E(T))$ possessing small constants $c_f$, the approximations converge geometrically fast down to a maximal accuracy of order $\epsilon^{1-a(\eta^{-1};T)}$.  In all other cases, the convergence is at least superalgebraic down to the same accuracy.

\section{Conclusions and challenges}\label{s:conclusions}
We conclude by making the following remark.  Extensive numerical
experiments \cite{BoydFourCont,BoydRunge,brunoFEP,LyonFESVD,LyonFast}
have shown the effectiveness of FEs in approximating even badly
behaved functions to high accuracy in a stable fashion.  The purpose of this paper has been to provide analysis to explain these results.  In particular, we have shown numerical stability for all three types of extensions considered, and analyzed their convergence.  The reason for this robustness, despite the presence of exponentially ill-conditioned matrices, is due to the fact that the FE is a frame
approximation and that for all functions $f$, even those with
oscillations or large derivatives, there eventually exist
coefficient vectors with small norms which approximate $f$ to high
accuracy.

The main outstanding theoretical challenge is to understand the constants $C_i(N,M;T,\epsilon)$ of the equispaced FE.  In particular, we wish to show that linear scaling $M = \gamma N$ is sufficient to ensure boundedness of these constants in $N$, with a  larger $\gamma$ corresponding to a smaller bound.  Note that the analysis of \S \ref{sss:oversamp} implies the suboptimal result that $M = \ord{N^2}$ is sufficient (Remark \ref{r:scaling}).  It is also a relatively straightforward exercise to show that if $M = c N / \epsilon$ for suitable $c > 0$, then $C_i(N,M;T,\epsilon)$ is bounded.  This is based on making rigorous the arguments given in Remark \ref{r:scaling}---we do not report it here for brevity's sake.  Unfortunately, although this estimate gives the correct scaling $M = \ord{N}$, it is wildly pessimistic.  It implies that $M$ should scale like $\approx 10^{16} N$, whereas the numerics in \S \ref{Cibehaviour} indicate that $M= \gamma N$ is sufficient for \textit{any} $\gamma \geq 1$.

One approach towards establishing a more satisfactory result is to perform a closer analysis of the singular values of the matrix $\bar{A}$.  Some preliminary insight into this problem was given in \cite{EdelmanFuture}.  Therein it was proved that (whenever $M=N$ and $2T \in \bbN$) the singular values cluster near zero and one, and the transition region is $\ord{\log N}$ in width, much like for the prolate matrix $A$.  Unfortunately, little is known outside of this result.  There is no existing analysis for $\bar A$ akin to that of Slepian's for the prolate matrix---see  \cite{EdelmanFuture} for a discussion.  Note, however, that the normal form $B = \bar{A}^*\bar{A}$ has entries $B_{n,m} = \frac{\sin \frac{(n-m) \pi}{T}}{M T \sin \frac{(n-m) \pi}{M T}}$, and can therefore be viewed as a discretized version of the prolate matrix $A$.  Indeed, $B \rightarrow A$ as $M \rightarrow \infty$ for fixed $N$.  Given the similarities between the two matrices, there is potential for Slepian's analysis to be extended to this case.  However,
 this remains an open problem.

Another issue is that of understanding how to choose the parameters $T$ and $\gamma$ in the case of the equispaced extension.  As discussed in \S \ref{ss:Tchoice}, the choice of $T$ is reasonably clear for the continuous and discrete FEs (where there is no $\gamma$).  If resolution of oscillatory functions is a concern, one should choose a small value of $T$ (in particular, \R{Tkte}).  Otherwise, a good choice appears to be $T=2$.  However, for the equispaced FE, small $T$ adversely affects stability (see \S \ref{ss:equistabfn}).  Hence it must be balanced by taking a larger value of the oversampling parameter $\gamma$, which has the effect of reducing the effective resolution power.  In practice, however, a reasonable choice appears to be $T = \gamma = 2$.  Investigating whether or not this is optimal is a topic for further investigation.

\section*{Acknowledgements}
The authors would like to thank John Boyd, Doug Cochran, Laurent Demanet, Anne Gelb, Anders Hansen, Arieh Iserles, Arno Kuijlaars, Mark Lyon, Nilima Nigam, Sheehan Olver, Rodrigo Platte, Jie Shen and Nick Trefethen for useful discussions and comments.  They would also like to thank the anonymous referees for their constructive and helpful remarks.

\bibliographystyle{abbrv}
\small
\bibliography{FEStabilityRefs}

\newpage
\section*{Symbols}

\small

\renewcommand{\arraystretch}{1.2}
\begin{longtable}{|c|c|p{9.0cm}|}
 \hline Symbol &  Section &
Description
\\ \hline
$T$ & \ref{ss:FE} & Extension parameter
\\ \hline
$N$ & \ref{ss:FE} & Truncation parameter
\\ \hline
$M$, $\gamma$ & \ref{ss:intro_summary_equispaced} & Number of equispaced nodes of the equispaced
FE, and the oversampling parameter $\gamma = M / N$
 \\ \hline
 $ \phi_n(x) $ & \ref{ss:FE} & The exponential $ \frac{1}{\sqrt{2T}} \E^{\I \frac{n \pi}{T}
x } $
\\ \hline
$ \cG_{N} $, $\cS_{N}$, $\cC_{N}$  &  \ref{ss:FE}, \ref{ss:FEinterp} &
Finite-dimensional spaces of exponentials, sines and
cosines
\\ \hline
$F_N$, $\tilde{F}_N(f)$, $F_{N,M}(f)$  & \ref{ss:FE}, \ref{ss:intro_summary_equispaced} & Exact continuous, discrete and equispaced FEs
\\ \hline
$G_N$, $\tilde{G}_N(f)$, $G_{N,M}(f)$  & \ref{ss:intro_summary}, \ref{ss:intro_summary_equispaced} & Numerical
 continuous, discrete and equispaced FEs
\\ \hline
$a$ & \ref{ss:intro_summary}, \ref{s:FEtypes} & Vector of coefficients of an FE
\\ \hline
$A$, $\tilde{A}$, $\bar{A}$  & \ref{ss:intro_summary}, \ref{s:FEtypes} & Matrices of
the continuous, discrete and equispaced FE`s
\\ \hline
$b$, $\tilde{b}$, $\bar{b}$ & \ref{s:FEtypes}, \ref{ss:equiFEdef}  & Data
vectors for the continuous, discrete and equispaced FEs
\\ \hline
$x$, $y$, $z$ & \ref{ss:FEinterp} & Physical domain variable $x \in [-1,1]$, and the mapped variables $y\in [c(T),1]$
 and $z \in [-1,1]$
 \\ \hline
$f_e(x)$, $f_o(x)$ & \ref{ss:FEinterp} & Even and odd parts of the function
$f(x)$
\\ \hline
$ g_1(y) $, $ g_2(y) $, $g_{1,N}(y)$, $g_{2,N}(y)$ & \ref{ss:FEinterp} & Images
of $ f_e(x)$ and $ f_o(x) / \sin \frac{\pi}{T} x$ in the
$y$-domain and their polynomial approximations
 \\ \hline
$ h_i(z) $, $ h_{i,N}(z) $ & \ref{ss:FEinterp} & Images of $g_i$ and $g_{i,N}$
in the $z$-domain
\\ \hline
$m(x)$ & \ref{ss:FEinterp} &   The mapping $x \mapsto z$
\\ \hline
$c(T)$, $E(T)$ & \ref{ss:intro_summary}, \ref{ss:FEinterp} & FE constants $\cos \frac{\pi}{T}$ and
$\cot^2 \left ( \frac{\pi}{4 T} \right )$.
\\ \hline
$\cB(\rho) $, $\cD(\rho)$ & \ref{s:FEconv} & Bernstein ellipse  in the
$z$-domain and its image in the $x$-domain
\\ \hline
$\kappa(F)$ & \ref{ss:cond_numb_exact} & Condition number of a mapping $F$
\\ \hline
$N_0$, $N_1$, $N_2$ & \ref{ss:numsolnanalysis}, \ref{ss:equianalysis} & Breakpoints in convergence
\\ \hline
  $\{ u_n , \sigma_n, v_n \}$ & \ref{ss:numsolnanalysis}, \ref{ss:equianalysis} & Singular system of $A$, $\tilde{A}$ or $\bar{A}$
\\ \hline
$ \Phi_n $ & \ref{ss:numsolnanalysis} & Fourier series corresponding to $v_n$
 \\ \hline
 $ \cG_{N,\epsilon} $, $ \cG'_{N,\epsilon} $,
$\cG_{N,M,\epsilon}$    & \ref{ss:numsolnanalysis}, \ref{ss:equianalysis}    & The subspace $ \spn \{ \Phi_n :
\sigma_n > \epsilon \}$ 
 \\ \hline
$H_{N,\epsilon}(f)$, $\tilde{H}_{N,\epsilon}(f)$,
$H_{N,M,\epsilon}(f)$ & \ref{ss:numsolnanalysis}, \ref{ss:equianalysis} & Truncated SVD FEs
corresponding to the continuous, discrete and equispaced cases
\\ \hline
$a(\gamma ; T )$ & \ref{ss:equianalysis} & Quantity determining the maximal achievable accuracy of the equispaced FE
\\ \hline
$\rL^2(I)$, $\ip{\cdot}{\cdot}_I$, $\nm{\cdot}_I$ & N/A &
Space of square-integral functions on a domain $I$ and
corresponding inner product and norm
\\ \hline
$\ip{\cdot}{\cdot}$, $\nm{\cdot}$ & N/A  & Inner product and
norm on $\rL^2(-1,1)$
\\ \hline
$\rL^2_w(I)$, $\ip{\cdot}{\cdot}_{w,I}$, $\nm{\cdot}_{w,I}$ &
 N/A & Space of square integrable functions with respect to a
weight function $w$ and corresponding inner product and norm
\\ \hline
$\nm{\cdot}_{\infty,I}$, $\nm{\cdot}_{\infty}$ &  N/A  &
Uniform norms on an arbitrary domain $I$ and the interval $[-1,1]$
respectively
\\ \hline
\end{longtable}

\end{document}